\documentclass[12pt]{article}
\usepackage[margin=35mm]{geometry}
\usepackage{amsmath}
\usepackage{amsfonts}
\usepackage{amssymb}
\usepackage{amsthm}
\usepackage{graphicx}
\usepackage{mathrsfs}
\usepackage{verbatim}
\immediate\write18{texcount -tex -sum  \jobname.tex > \jobname.wordcount.tex}

\newtheorem*{claim*}{Claim}
\newtheorem{claim}{Claim}
\newtheorem{theorem}{Theorem}
\newtheorem{lemma}{Lemma}

\newtheorem{corollary}{Corollary}
\newtheorem{conjecture}{Problem}
\usepackage[hidelinks]{hyperref}

\providecommand{\keywords}[1]
{
  \small	
  \textbf{\textit{Keywords---}} #1
}

\title{The Gamma-Theta Conjecture holds for planar graphs}
\author{Dmitrii Taletskii$^{1}$\\
\small
$^1$ National Research University Higher School of Economics, \\
\small
Bolshaya Pechyorskaya ul. 25/12, Nizhny Novgorod, 603155 Russia\\%
    \small
    e-mail: dmitalmail@gmail.com
}
\date{} % Comment this line to show today's date

\begin{document}
\maketitle

\begin{abstract}
The Gamma-Theta Conjecture states that if the domination number of a graph is equal to its eternal domination number, then it is also equal to its clique covering number. This conjecture is known to be true for several graph classes, such as outerplanar graphs, subcubic graphs and $C_k$-free graphs, where $k \in \{3,4\}$. In this paper, we prove the Conjecture for the class of planar graphs.
\end{abstract} \hspace{10pt}

\keywords{dominating set, eternal dominating set, planar graph}

\section{Introduction}

In this paper, all graphs are finite, simple and undirected.  Let $G$ be a graph with a vertex set $V(G)$ and an edge set $E(G)$. A set $D \subseteq V(G)$ is called \textit{dominating} if every vertex not from $D$ has a neighbor in $D$. \textit{An eternal dominating set} of $G$ is a dominating set that can defend any infinite series of vertex attacks,
where an attack is defended by moving one guard along an edge from its current
position to the attacked vertex. The cardinality of a minimum dominating  set (eternal dominating set) of $G$ is called \textit{the domination number} (\textit{the eternal domination number}) and denoted by $\gamma(G)$ ($\gamma^{\infty}(G)$). A set $I \subseteq V(G)$ is called \textit{independent}, if its vertices are pairwise nonadjacent. The cardinality of a maximum independent set of $G$ is called \textit{the independence number} and denoted by $\alpha(G)$.  \textit{A clique partition} (CP) of a graph $G$ is a partition of $V(G)$ into cliques. The cardinality of a minimum clique partition (MCP) of $G$ (that is, the minimum number of cliques needed to cover $V(G)$) is called \textit{the clique covering number}, and denoted by $\theta(G)$.

% In other words, a dominating set $D \subseteq V(G)$ is eternal dominating in $G$, if for every sequence  $v_1, \ldots, v_s$ of $s \geq 1$ vertices from $V(G)$, there exist dominating sets $D_0 = D,D_1, \ldots,D_{s}$ and vertices $u_0 \in D_0 , u_1 \in D_1, \ldots, u_{s-1} \in D_{s-1}$ such that the vertex $u_i$ is adjacent to $v_{i+1}$ and $D_{i+1} = (D_{i} \setminus \{u_{i}\}) \cup \{v_{i+1}\}$ for all  $i \in [0;s-1]$.

The study of the eternal domination number started in~\cite{BCGMVW04}, where the trivial relation $\gamma(G) \leq \alpha(G) \leq \gamma^{\infty}(G) \leq \theta(G)$ was first mentioned. The inequality $\gamma^{\infty}(G) \leq \binom{\alpha(G)}{2}$ was proved in~\cite{KM07} for all graphs $G$. In~\cite{GK08} graphs $G'$ were found such that $\gamma^{\infty}(G') = \binom{\alpha(G')}{2}$. In~\cite{DKKM20}, the inequality $\gamma^{\infty}(G \boxtimes H) \geq \alpha(G) \cdot \gamma^{\infty}(H)$ was obtained (here $G \boxtimes H$ denotes the strong product of graphs $G$ and $H$). 

In this paper, we study the eternal domination model known as \textit{one guard moves} (as mentioned earlier, in this model only one guard is allowed to move to a neighboring vertex after each attack). There is another well-known model \textit{all guards move}, which was introduced in~\cite{GHH05} (in this model all guards are allowed to move after each attack). There are also some generalizations like oriented~\cite{BJK21} and fractional~\cite{DKKMS23}  eternal domination. See a survey on eternal domination in~\cite{KM20}.

Following~\cite{ABBCVY07}, we call a graph $G$  \textit{maximum-demand}, if $\gamma^{\infty}(G) = \theta(G)$. The following open problems are of interest:

\begin{conjecture}
Characterize the class of maximum-demand graphs.
\end{conjecture} 

\begin{conjecture}[The Gamma-Theta Conjecture]
For every graph $G$, if $\gamma(G) = \gamma^{\infty}(G)$, then $\gamma^{\infty}(G) = \theta(G)$.
\end{conjecture} 

We summarize the known results on these problems in the following theorems.

\begin{theorem}
All graphs from the following classes are maximum-demand:

(a) \cite{BCGMVW04}~Perfect graphs;

(b) \cite{BCGMVW04}~Graphs with the clique covering number at most 3;

(c) \cite{ABBCVY07}~$K_4$-minor-free graphs;

(d) \cite{ABBCVY07}~The cartesian products $C_m \square P_m$, $C_m \square C_n$ for all $m,n \geq 1$.

(d) \cite{R07}~Circular-arc graphs;

(f) \cite{MMV22}~Graphs with at most 9 vertices;

(g) \cite{MMV22}~Triangle-free graphs with at most 12 vertices;

(i) \cite{MMV22}~Planar graphs with at most 11 vertices;

(j) \cite{MMV22}~Cubic graphs with at most 16 vertices.
\end{theorem}

\begin{theorem}
The Gamma-Theta conjecture holds for the following classes of graphs:

(a) Maximum-demand graphs;

(b) Subcubic graphs~\cite{KM15};

(c) Triangle-free~\cite{KM15} and $C_4$-free~\cite{KKM18} graphs.
\end{theorem}

 It was proved in~\cite{GHH05} that not all graphs are maximum-demand and an example of a graph on 11 vertices that is not maximum-demand was given. Later, in~\cite{MMV22}, it was shown that the smallest such graphs have 10 vertices, the eternal domination number~3 and the clique covering number~4 (there are two such graphs up to isomorphism). It was also shown in~\cite{MMV22} that the smallest triangle-free graphs that are not maximum-demand have 13 vertices (there are 13 such graphs up to isomorphism).  To date, there are no known examples of planar or cubic graphs that are not maximum-demand. Note that all outerplanar graphs are maximum-demand, since every such graph is $K_4$-minor-free.  It was proved in~\cite{KSV20} that there are infinitely many graphs $G$ such that $G$ is maximum-demand, but its prism (that is, the cartesian product $G \square K_2$) is not.

This paper proves the following fact:

\begin{theorem}\label{thm3}
For every planar graph $G$, if $\gamma(G) = \gamma^{\infty}(G)$, then $\gamma^{\infty}(G) = \theta(G)$.
\end{theorem}

We introduce a new technique based on \textit{$\theta$-independent sets}, that is, vertex subsets with no two vertices from the same clique of a given minimum clique partition of a graph. Supposing for a contradiction that a minimal (by the number of vertices) planar counterexample $G$ exists, we consider its minimum degree vertex $v$ and the subgraph $G[N_1(v)]$ induced by the vertices from the open neighborhood of $v$. The number of such subgraphs (up to isomorphism) is relatively small, since the minimum vertex degree of a planar graph is at most 5. This allows us to consider all possible subgraphs in Lemmas~\ref{l18}--\ref{l27} and~\ref{l30}--\ref{l40} and conclude that a minimal counterexample does not exist. Therefore, Theorem~\ref{thm3} follows instantly from the 21 lemmas mentioned above.

This paper is rather long. There are two main reasons for this. First, we need to obtain a large number of auxiliary facts and properties in order to prove the main result. Some facts about the structure of a minimal counterexample (Lemma~\ref{l15} and especially Lemma~\ref{l17}) are difficult to prove. Second, we have not found a universal approach to check the configurations studied in Lemmas~\ref{l18}--\ref{l27} and~\ref{l30}--\ref{l40}. Hence we have to use slightly different approaches for different configurations, also leading to a lengthy proof.

\section{Terminology}

\subsection{Basics}

Let $G$ be a graph. Denote by $\delta(G)$ ($\Delta(G)$) the minimum (maximum) vertex degree of~$G$. A vertex with degree 0 is called \textit{isolated}; a vertex with degree 1 is called \textit{pendant}. Denote by $G[W]$ the induced subgraph of $G$ with the vertex set $W \subseteq V(G)$. For a vertex $v \in V(G)$ and an integer $k \geq 1$, denote by $N_k(v)$ ($N_k[v]$) the set of all vertices on distance exactly $k$ (at most $k$) from $v$. For a nonempty set $U \subseteq V(G)$, let $N_k(U) = \cup_{u \in U}N_k(u)$ and $N_k[U] = \cup_{u \in U}N_k[u]$. Denote by $G_U$ the induced subgraph $G[V(G) \setminus N_1[U]]$. If $U = \{u\}$, we write $G_u$ instead of $G_{\{u\}}$; if $U = \{u,v\}$, we write $G_{u,v}$ instead of $G_{\{u,v\}}$. 

Let $A,B \subseteq V(G)$. We say that $A$ \textit{dominates} $B$, if every vertex from $B$ belongs to $A$ or has a neighbor in $A$. If $A = \{a\}$, we say that $a$ dominates $B$. Let $v \in V(G)$ and $W \subseteq V(G)$ (here $v \notin W$). We say that $v$ is \textit{adjacent} (\textit{not adjacent}) to $W$ if $v$ has a neighbor (does not have a neighbor) in $W$.  A graph $G$ is called \textit{$H$-free} if it does not contain a graph $H$ as a subgraph.

For every planar graph that we consider throughtout the paper, we assume that it has a fixed embedding in the plane. Let $G$ be a plane graph and $C$ be its cycle. Denote by $int(C)$ ($ext(C)$) the set of vertices inside (outside) $C$. If both $int(C)$ and $ext(C)$ are nonempty, $C$ is said to be \textit{separating}. We call an edge $ab \in E(G)$ \textit{separating}, if the induced subgraph $G[V(G) \setminus \{a,b\}]$ is disconnected. 

Denote by $IS(G), DS(G), MDS(G), EDS(G), MEDS(G)$ the families of all independent, dominating, minimum dominating, eternal dominating and minimum eternal dominating sets of a graph $G$, respectively. Likewise, denote by $CP(G)$ and $MCP(G)$ the families of all clique partitions and minimum clique partitions of $G$, respectively. 

Call a graph $G$ \textit{critical}, if it is the smallest by the number of vertices planar graph such that $\gamma(G) = \gamma^{\infty}(G) < \theta(G)$.

\subsection{$\theta$-independent sets} 

Let $G$ be a graph and $\mathcal{C} \in MCP(G)$. Call a nonempty set $S \subseteq V(G)$ \textit{$\theta$-independent in $\mathcal{C}$} if no two vertices of $S$ belong to the same member of $\mathcal{C}$. We say that $S$ is \textit{$\theta$-independent in $G$} if there exists a family $\mathcal{C}' \in MCP(G)$ such that $S$ is $\theta$-independent in $\mathcal{C}'$. Note that every independent set $J \subseteq V(G)$ is $\theta$-independent in $G$, but the converse is not true.

Consider the graph $H$ on a Fig.~1 and its MCP $\mathcal{C} = \{\{u_1,u_4\},\{u_2,u_5\},\{u_3,u_6\}\}$. The set $\{u_1,u_3,u_5\}$ is independent.  The set $\{u_1,u_2,u_3\}$ is not independent, but it is $\theta$-independent in $\mathcal{C}$. The set $\{u_1,u_3,u_4\}$ is not $\theta$-independent in $\mathcal{C}$ but it is $\theta$-independent in $H$, since there exists another MCP $\mathcal{C}' =\{\{u_1,u_2\}, \{u_4,u_5\}, \{u_3,u_6\}\}$, such that $\{u_1,u_3,u_4\}$ is $\theta$-independent in $\mathcal{C}'$. Finally, the set $\{u_1,u_2,u_4\}$ is not $\theta$-independent in~$H$, since for every family $\mathcal{C}^* \in MCP(H)$ either $\{u_1,u_2\} \in \mathcal{C}^*$ or $\{u_1,u_4\} \in \mathcal{C}^*$.

\begin{center}
\includegraphics[scale=0.45]{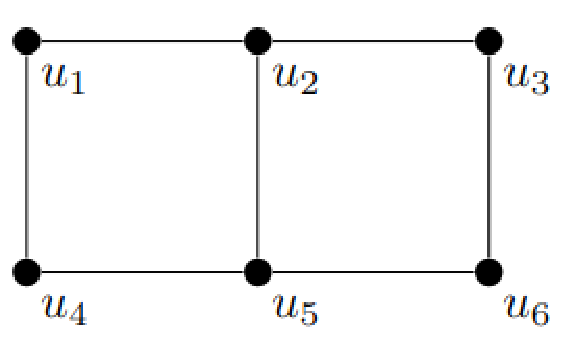}

Fig.~1. The graph $H$
\end{center}

\subsection{Strategies}

Let $G$ be a graph and $D \in DS(G)$. Call \textit{an attack sequence} (AS) a finite sequence $\mathfrak{A}$ of vertices from $V(G)$. A triple $(G,D,\mathfrak{A})$ is called \textit{a strategy}. A strategy $(G,D,\mathfrak{A})$ is said to be \textit{losing}, if the guards from the set $D$ can protect the vertices of $G$ from a sequence of attacks $\mathfrak{A}$, and \textit{winning} otherwise. More formally, a strategy $(G,D,v_1 \ldots v_k)$ is losing, if for all $i \in [1;k]$ there exist sets $D_i \in DS(G)$ and vertices $u_i \in D_{i-1} \cap N_1[v_i]$ such that  $D_i = (D_{i-1} \setminus \{u_i\}) \cup \{v_i\}$. Here $D_0 = D$.

Consider the graph $H$ on a Fig.~1 and its dominating set $D = \{u_2,u_5\}$. It is easy to see that the strategy $(H,D,u_1u_2u_3)$ is losing. However, the strategy $(H,D,u_1u_3)$ is winning, hence $D$ is not eternal dominating in $G$.

\section{Preliminary results}

\subsection{Planar graphs}

\begin{lemma}\label{l1}
Let $G$ be a planar graph and $uv \in E(G)$. Moreover, let $N_1[u] \setminus N_1[v] = \{u_1,\ldots,u_p\}$ and $N_1[v] \setminus N_1[u] = \{v_1,\ldots,v_q\}$. Then the following holds:

(a) If $p,q \geq 2$, then there exist integers $i \in [1;p]$ and $j \in [1;q]$ such that $u_iv_j \notin E(G)$.

(b) If $p,q \geq 3$, then there exist integers $i,i' \in [1;p]$ and $j,j' \in [1;q]$ such that $i < i'$, $j < j'$ and $u_iv_j,u_{i'}v_{j'} \notin E(G)$.
\end{lemma}

\begin{proof}
The first statement of the lemma is obvious, since $G$ does not contain $K_{3,3}$ as a subgraph. We now prove the second statement. Assume that $u_1v_1 \notin E(G)$. Then, by the first statement, there exist integers $i' \in [2;p]$ and $j' \in [2;q]$ such that $u_{i'}v_{j'} \notin E(G)$, as required.
\end{proof}

\begin{lemma}\label{l2}
Let $G$ be a plane graph with a separating cycle $C$ such that every vertex of $G[int(C)]$ has a neighbor in $C$. If $|int(C)| \geq 2$ $(|int(C)| \geq 4)$, then $G[int(C)]$ has two vertices (two nonadjacent vertices) with at least $\delta(G) - 2$ neighbors in~$C$.
\end{lemma}

\begin{proof} The subgraph $G[int(C)]$ is outerplanar, since every its vertex has a neighbor in~$C$. It is well known that every outerplanar graph with at least 2 vertices (with at least 4 vertices) has two vertices (two nonadjacent vertices) of degree 2 or less. Therefore, $G[int(C)]$ has two vertices (two nonadjacent vertices) with at least $\delta(G) - 2$ neighbors in $V(G) \setminus int(C)$, as required.
\end{proof}

\subsection{Eternal dominating sets}

First, we mention several known facts.

\begin{lemma}[\cite{BCGMVW04}]\label{l3}
For every graph $G$ we have $\gamma(G) \leq \alpha(G) \leq \gamma^{\infty}(G) \leq \theta(G).$
\end{lemma}

\begin{lemma}[\cite{ABBCVY07}, Preposition~11]\label{l4} 
Let $G$ be a graph with a cutvertex $v$. Moreover, $G$ is obtained from graphs $G_1$ and $G_2$ by identifying vertices $v_1 \in V(G_1)$ and $v_2 \in V(G_2)$ with~$v$. If the graphs $G_1$, $G_2$, $G_1 \setminus v_1$, $G_2 \setminus v_2$ are maximum-demand, then $G$ is also maximum-demand.
\end{lemma}

\begin{lemma}[\cite{MMV22}, Observation~29]\label{l5} Every planar graph on 11 vertices or less is maximum-demand.
\end{lemma}

\begin{lemma}[\cite{KM09}]\label{l6}
Let $H$ be an induced subgraph of $G$, then $\gamma^{\infty}(H) \leq \gamma^{\infty}(G)$.
\end{lemma}

\textbf{Remark~1.} Clearly, if $H$ is an induced subgraph of $G$, then $\theta(H) \leq \theta(G)$. However, this property does not hold for the domination number. Consider the complete bipartite graph $K_{1,n}$ and its induced empty subgraph $\overline{K_n}$. Then $\gamma(\overline{K_n}) > \gamma(K_{1,n})$ for all $n \geq 1$.

We now prove a generalization of Lemma~\ref{l6}.

\begin{lemma}\label{l7}
Let $H$ be an induced subgraph of $G$, then the following holds:

(a) There exists a set $D_1 \in MEDS(G)$ such that $|D_1 \cap  V(H)| \geq \gamma^{\infty}(H)$.

(b) If $\gamma^{\infty}(H) = \gamma^{\infty}(G)$, then for every induced subgraph $H'$ of $H$, there exists a set $D_2 \in MEDS(G)$, such that $D_2 \subseteq V(H)$ and $\ |D_2 \cap V(H')| \geq \gamma^{\infty}(H')$.

\end{lemma}

\begin{proof}
(a) Consider a set $D \in MEDS(G)$ such that the cardinality $|D \cap V(H)|$ is maximum possible. Suppose that $|D \cap V(H)| < \gamma^{\infty}(H)$, then there exists an AS $\mathfrak{A}$ such that the strategy $(H,D \cap V(H),\mathfrak{A})$ is winning. Since no set from $MEDS(G)$ has more than $|D \cap V(H)|$ guards in $V(H)$, the strategy $(G,D,\mathfrak{A})$ is also winning, a contradiction.

(b) It follows from (a) that there exists a set $D_1 \in MEDS(G)$ such that $D_1 \subseteq V(H)$. We now apply (a) again to the graph $H$ and its induced subgraph $H'$ and obtain a required set $D_2$.
\end{proof}

\subsection{Graphs with $\gamma < \theta$}

In this subsection we obtain two sufficient conditions for a planar graph to have the dominating number less than the clique covering number. 

\begin{lemma}\label{l8}
Let $G$ be a graph such that $\gamma(G) = \theta(G)$. If a set $J \subseteq V(G)$ is $\theta$-independent in $G$, then there exists a set $J' \in MDS(G)$ such that $J \subseteq J'$.
\end{lemma}

\begin{proof}
Since $J$ is $\theta$-independent in $G$, there exists a family $\mathcal{C} \in MCP(G)$ such that $J$ is $\theta$-independent in $\mathcal{C}$. If $|J| = \theta(G)$, then $J$ has a common vertex with every member of $\mathcal{C}$, hence $J \in MDS(G)$ and we can choose $J' = J$. If $|J| < \theta(G)$, select a vertex from every member of $\mathcal{C}$ not intersecting with $J$ and obtain a set $J^*$. Then $J \cup J^*$ is $\theta$-independent set in $\mathcal{C}$ and $|J \cup J^*| = \theta(G)$. Clearly, $J \cup J^* \in MDS(G)$, thus we can choose $J' = J \cup J^*$. 
\end{proof}

\textbf{Remark~2.} Lemma~\ref{l8} implies that if $\gamma(G) = \theta(G)$, then for every vertex $v \in V(G)$ there exists a set $D \in MDS(G)$ such that $v \in D$.

\begin{lemma}\label{l9}
Let $G$ be a graph and $\mathcal{C} \in MCP(G)$. Moreover, let $J$ be a $\theta$-independent set in $\mathcal{C}$. If there exists a clique $W \in \mathcal{C}$ such that $J$ dominates $W$ and $J \cap W = \emptyset$, then $\gamma(G) < \theta(G)$.
\end{lemma}

\begin{proof}
Since $J \cap W = \emptyset$ and $J$ is $\theta$-independent in $\mathcal{C}$, we have $|J| \leq \theta(G) - 1$. If $|J| = \theta(G) - 1$, then every member of $\mathcal{C} \setminus \{W\}$ intersects with $J$, thus $J \in DS(G)$ and $\gamma(G) < \theta(G)$, as required. Suppose that $|J| < \theta(G) - 1$. If $\gamma(G \setminus W) < \theta(G \setminus W)$, then $$\gamma(G) \leq \gamma(G \setminus W) + 1 < \theta(G \setminus W) + 1 = \theta(G)$$ and we are done. Otherwise by the previous lemma there exists a set $J^* \in DS(G \setminus W)$ such that $J \subseteq J^*$ and $|J^*| = \theta(G) - 1$. Then $J^* \in DS(G)$ and $\gamma(G) < \theta(G)$, as required.
\end{proof}

\textbf{Remark~3.} Lemma~\ref{l9} implies that if $\gamma(G) = \theta(G)$ and $\mathcal{C} \in MCP(G)$, then all members of $\mathcal{C}$ are maximal by inclusion in $V(G)$. In particular, if $G$ has no isolated vertices, then $\mathcal{C}$ has no 1-vertex cliques.

\begin{lemma}\label{l10}
Let $G$ be a planar graph and $v \in V(G)$ be its vertex satisfying one of the properties (a)--(e). Then $\gamma(G) < \theta(G)$.

(a) $v$ is adjacent to at least two pendant vertices;

(b) $N_1(v) = \{u_1,u_2\}$, $u_1u_2 \notin E(G)$. Moreover, there exist vertices $u'_1 \in N_1(u_1) \setminus N_1(u_2)$ and $u'_2 \in N_1(u_2) \setminus N_1(u_1)$.

(c) $N_1(v) = \{u_1,u_2\}$, $u_1u_2 \notin E(G)$. Moreover, there exist nonadjacent vertices $x,y \in (N_1(u_1) \cap N_1(u_2)) \setminus \{v\}$.

(d) $N_1(v) = \{u_1,u_2\}$, $u_1u_2 \notin E(G)$ and $\min(deg(u_1),deg(u_2)) \geq 4$.

(e) $N_1(v) = \{u_1,\ldots,u_s\}$, where $s \geq 3$. Moreover, the vertices from $N_1(v)$ are not pendant and $N_1(v) \in IS(G)$.
\end{lemma}

\begin{proof}
(a) Suppose $v$ is adjacent to pendant vertices $u_1$ and $u_2$. Then $u_1u_2 \notin E(G)$ and there exists a set $D \in MEDS(G)$ such that $\{u_1,u_2\} \subseteq D$. However, $(D \setminus \{u_1,u_2\}) \cup \{v\} \in DS(G)$, thus $\gamma(G) < \gamma^{\infty}(G) = |D| \leq \theta(G)$, as required. 

(b) Consider a family $\mathcal{C} \in MCP(G)$ and its clique $W \ni v$. Assume by symmetry that $W = \{v,u_1\}$. The set $\{u'_1,u_2\} \in IS(G)$ is $\theta$-independent in $\mathcal{C}$ and dominates $\{v,u_1\}$, thus $\gamma(G) < \theta(G)$ by the previous lemma.

(c) Consider a family $\mathcal{C} \in MCP(G)$ and its clique $W \ni v$, assume that $W = \{v,u_1\}$. Clearly, one of the sets $\{u_2,x\}$ and $\{u_2,y\}$ is $\theta$-independent in $\mathcal{C}$ and dominates $\{v,u_1\}$, thus $\gamma(G) < \theta(G)$ by the previous lemma.

(d) It follows from (c) that the set $(N_1(u_1) \cap N_1(u_2))  \setminus \{v\}$ is a clique. Since $G$ is planar, $|(N_1(u_1) \cap N_1(u_2)) \setminus \{v\}| \leq 2$.  Then both sets $N_1(u_1) \setminus N_1(u_2)$ and $N_1(u_2) \setminus N_1(u_1)$ are nonempty; this contradicts~(b). 

(e) Consider a family $\mathcal{C} \in MCP(G)$ and its clique $W \ni v$. Assume by symmetry that $W = \{v,u_1\}$. Since $u_1$ is not pendant, there exists a vertex $u'_1 \in N_1[u_1] \setminus N_1[v]$. Since $u_2u_3 \notin E(G)$, one of the sets $\{u'_1,u_2\}$ and $\{u'_1,u_3\}$ is $\theta$-independent in $\mathcal{C}$ and dominates $\{v,u_1\}$. Thus $\gamma(G) < \theta(G)$ by the previous lemma. 
\end{proof}

\textbf{Remark~4.} The conditions of Lemma~9 and Lemma~10(a,b,c,e) also hold for nonplanar graphs (though we do not use this fact in what follows).

\section{Properties of a critical graph}

In this section we consider a critical graph $G$, that is, a minimal by the number of vertices planar graph such that $\gamma(G) = \gamma^{\infty}(G) < \theta(G)$. In subsection~\ref{s41} we obtain a few simple facts and outline the proof of the main result. In subsections~\ref{s42} and~\ref{s43} we prove more difficult structural properties.

\subsection{Induced subgraphs}\label{s41}

\begin{lemma}\label{l11}
Let $ W \subseteq V(G)$ be a nonempty clique, then $\gamma^{\infty}(G) = \gamma^{\infty}(G \setminus W)$.

\end{lemma}

\begin{proof}

Suppose for a contradiction that $\gamma^{\infty}(G) \neq \gamma^{\infty}(G \setminus W)$. Lemmas~\ref{l3} and~\ref{l6} imply that $$\gamma(G \setminus W) \leq \gamma^{\infty}(G \setminus W) < \gamma^{\infty}(G).$$ Consider a set $D \in DS(G \setminus W)$ and a vertex $w \in W$. Clearly, $D \cup \{w\} \in DS(G)$, hence $\gamma(G \setminus W) = \gamma^{\infty}(G \setminus W) = \gamma(G) - 1$. Consider a family $\mathcal{C} \in MCP(G \setminus W)$. Since $\mathcal{C} \cup\{W\} \in CP(G)$, we have $$\theta(G \setminus W) \geq \theta(G) - 1 > \gamma(G \setminus W) = \gamma^{\infty}(G \setminus W).$$
Therefore, $G$ is not critical, a contradiction.
\end{proof}

\begin{lemma}\label{l12}
For every edge $uv \in E(G)$ neither $N_1[u] \subseteq N_1[v]$ nor $N_1[v] \subseteq N_1[u]$. 
\end{lemma}

\begin{proof}
Suppose for a contradiction that $N_1[u] \subseteq N_1[v]$. Consider a set $D \in DS(G \setminus v)$. Clearly, $D$ dominates $v$, then $\gamma(G \setminus v) \geq \gamma(G)$. By Lemma~\ref{l6}, $\gamma^{\infty}(G \setminus v) \leq \gamma^{\infty}(G)$, hence $\gamma(G \setminus v) = \gamma^{\infty}(G \setminus v) = \gamma(G)$. We now prove that $\theta(G) = \theta(G \setminus v)$. Consider a family $\mathcal{C} \in MCP(G \setminus v)$ with a clique $W \ni u$. Clearly, the set $W' = W \cup \{v\}$ is also a clique, hence $(\mathcal{C} \setminus \{W\}) \cup \{W'\} \in CP(G)$ and $\theta(G \setminus v) = \theta(G)$. Therefore, $G$ is not critical, a contradiction.
\end{proof}

\textbf{Remark~5.} Lemma~\ref{l12} implies that for every vertex $x \in V(G)$ the set $N_1[x]$ is not a clique.

\begin{lemma}\label{l13}
For every nonempty set $I \subseteq IS(G)$ we have $$\gamma(G_{I}) = \theta(G_{I}) =  \gamma(G) - |I|.$$
\end{lemma}

\begin{proof}
Consider a set $D \in MDS(G_I)$. Clearly, $D \cup I \in DS(G)$. Therefore, 

\begin{equation}
\gamma(G) \leq \gamma(G_I) + |I|.
\end{equation}

Let $\mathfrak{A}$ be a sequence of all vertices from $I$ in some order. If there exists a set $D \in MEDS(G)$ such that $|N_1[I] \cap D| < |I|$, then the strategy $(G,D,\mathfrak{A})$ is winning, a contradiction. Then for every set $D' \in MEDS(G)$ we have $|V(G_I) \cap D'| \leq \gamma^{\infty}(G) - |I|$. Lemma~\ref{l7}(a) implies that 

\begin{equation}
\gamma^{\infty}(G_I) \leq \gamma^{\infty}(G) - |I|.
\end{equation}
 
Remind that $\gamma(G) = \gamma^{\infty}(G)$ and $\gamma(G_I) \leq \gamma^{\infty}(G_I)$. It follows from (1) and (2) that $\gamma(G_I) = \gamma^{\infty}(G_I) = \gamma(G) - |I|$. Since $G$ is critical,  $\gamma(G_{I}) = \theta(G_{I})$, as required.
\end{proof}

\textbf{Remark~6.} Lemma~\ref{l13} implies that for every vertex $a \in V(G)$ we have $$\gamma(G_{a}) = \theta(G_a) = \gamma(G) - 1.$$
Moreover, for every vertex $b \in V(G)$ such that $ab \notin E(G)$ we have
$$\gamma(G_{a,b}) = \theta(G_{a,b}) = \gamma(G) - 2.$$

\begin{lemma}\label{l14}
Let $ab \in E(G)$ and $Q = N_1[a] \cap N_1[b], \ A = N_1[a] \setminus Q, \ B = N_1[b] \setminus Q$. Then the following is true:

(a) if $G[Q] \cong K_m$, where $m \in [2;4]$, then 

\begin{equation}\label{eq3}
\gamma(G_{a,b}) = \theta(G_{a,b}) = \gamma(G) - 2.
\end{equation}

(b) if $|Q| \geq 4$ and the condition (\ref{eq3}) does not hold, then 

\begin{equation}\label{eq4}
\gamma(G \setminus Q) = \theta(G\setminus Q) = \gamma(G) - 1.
\end{equation}

\end{lemma}

\begin{proof}
Note that $A,B \neq \emptyset$ by Lemma~\ref{l12}. We now prove (a). By Lemmas~\ref{l6} and~\ref{l13}, $$\gamma^{\infty}(G_{a,b}) \leq \gamma^{\infty}(G_a) = \gamma(G) - 1.$$ 
Clearly, if $D \in MDS(G_{a,b})$, then $D \cup \{a,b\} \in DS(G)$, hence $\gamma(G_{a,b}) \geq \gamma(G) - 2$ and $\gamma^{\infty}(G_{a,b}) \in \{\gamma(G) - 2, \gamma(G) - 1\}$. 

\textbf{Case 1.} $\gamma^{\infty}(G_{a,b}) = \gamma(G) - 2$. Then $\gamma(G_{a,b}) = \gamma(G) - 2$. Since $G$ is critical, $\theta(G_{a,b}) = \gamma(G_{a,b})$ and we are done. 

\textbf{Case 2.} $\gamma^{\infty}(G_{a,b}) = \gamma^{\infty}(G) - 1$. By Lemma~\ref{l11}, $\gamma^{\infty}(G \setminus Q) = \gamma^{\infty}(G)$. Thus by Lemma~\ref{l7}(b) there exists a set $D' \in MEDS(G)$, such that 

\begin{equation}\label{eq5}
D' \cap Q = \emptyset, \ |D' \cap (A \cup B)| \leq 1.
\end{equation}
Since $A \cap B = \emptyset$, the set $D'$ does not dominate $\{a,b\}$, a contradiction.

We now prove (b). Again, it is easy to check, using Lemmas~\ref{l6} and~\ref{l13}, that $$\gamma(G) - 2 \leq \gamma(G_{a,b}) \leq \gamma^{\infty}(G_{a,b}) \leq \gamma^{\infty}(G_a) \leq \gamma(G) - 1.$$

\textbf{Case 1.} $\gamma^{\infty}(G_{a,b}) = \gamma(G) - 2$. Then $\gamma(G_{a,b}) = \gamma^{\infty}(G_{a,b})$. Since $G$ is critical, $\theta(G_{a,b}) = \gamma(G_{a,b})$, thus the equality (3) holds, a contradiction. 

\textbf{Case 2.} $\gamma^{\infty}(G_{a,b}) = \gamma(G) - 1$. Suppose that $\gamma^{\infty}(G \setminus Q) = \gamma(G) - 1$. Then it is easy to see that $\gamma(G \setminus Q) = \gamma(G) - 1$. Since $G$ is critical, $\gamma(G \setminus Q) = \theta(G \setminus Q)$ and we are done. Otherwise $\gamma^{\infty}(G \setminus Q) = \gamma^{\infty}(G)$. By Lemma~\ref{l7}(b), there exists a set $D'$ such that the condition~(\ref{eq5}) holds. Again, $D'$ does not dominate $\{a,b\}$, a contradiction.\end{proof}

\begin{corollary}\label{cor1}

Let $vu \in E(G)$ and $|N_1(v) \cap N_1(u)| \leq 1$. Then the following is true: 

(a) If a set $S \subseteq V(G_{v,u})$ is $\theta$-independent in $G_{v,u}$, then it does not dominate the set $N_1[v] \setminus N_1[u]$ in $G_u$.

(b) No vertex $w \in V(G_{v,u})$ dominates the set $N_1[v] \setminus N_1[u]$ in $G_u$.

\end{corollary}

\begin{proof}
We prove the statement (a) (the statement (b) is its special case with $|S| = 1$, which is often used in what follows). By Lemmas~\ref{l13} and~\ref{l14},
$$\gamma(G) - 2 = \gamma(G_{v,u}) = \theta(G_{v,u}) < \gamma(G_u) = \theta(G_u) = \gamma(G) - 1.$$

Note that $V(G_u) \setminus V(G_{v,u}) = N_1[v] \setminus N_1[u]$. Suppose that some set $S \subseteq V(G_{v,u})$ is $\theta$-independent in $G_{v,u}$ and dominates $N_1[v] \setminus N_1[u]$ in $G_u$. By Lemma~\ref{l8}, there exists a set $S^* \in DS(G_{v,u})$ such that $S \subseteq S^*$. Then $S^* \in DS(G_u)$ and, therefore, $\gamma(G_{v,u}) \geq \gamma(G_u)$, a contradiction.
\end{proof}

\textbf{Remark~7.} The main idea of the proof of Theorem~\ref{thm3} can now be roughly formulated as follows. We use the following approaches to show that a given graph $G'$ is not critical:
\begin{itemize}
\item Find an independent subset $X \subseteq V(G')$ (usually $|X| = 1$)  such that the subgraph $G'_X$ satisfies the condition of Lemma~\ref{l9} or Lemma~\ref{l10};

\item For a vertex $v \in V(G')$ such that $deg(v) = \delta(G')$, find a vertex $u \in N_1(v)$ and a set $S \subseteq V(G'_{v,u})$ such that $S$ is $\theta$-independent in $G'_{v,u}$ and dominates $N_1[v] \setminus N_1[u]$.
\end{itemize}

\subsection{Separating cycles}\label{s42}

The following lemma allows us to select, for a given separating cycle $C$ of a critical graph, a vertex in $int(C)$ not adjacent to $C$.

\begin{lemma}\label{l15}
Let $\delta(G) \geq 4$ and $C$ be a separating cycle of $G$. Moreover, $C = u_1u_2u_3$ or $C = u_1u_2u_3u_4$ and $u_4$ is not adjacent to $int(C)$. Then the following holds:

(a) if $\delta(G) = 5$, then $int(C) \nsubseteq N_1[C]$; 

(b) if $\delta(G) = 4$ and $int(C) \subseteq N_1[C]$, then there exists a vertex $w \in int(C)$ such that $G[N_1(w)] \cong C_4$. 
\end{lemma}

\begin{proof}
If $int(C) \nsubseteq N_1[C]$, then there is nothing to prove. Suppose that $int(C) \subseteq N_1[C]$. The subgraph $G[int(C)]$ is outerplanar.  Since $\delta(G) \geq 4$, we have $|int(G)| > 1$, thus by Lemma~\ref{l2}, there exist two vertices in $G[int(C)]$ with at least $\delta(G) - 2$ neighbors in $\{u_1,u_2,u_3\}$. Since $G$ is planar, $\delta(G) = 4$. Assume that $C = u_1u_2u_3$ (it is not hard to modify the proof for $C = u_1u_2u_3u_4$) and consider five cases.

\textbf{Case 1.} There exists a vertex $x \in int(C)$ such that for some $1 \leq i < j \leq 3$ the following holds: (i) $x$ dominates $\{u_i,u_j\}$; (ii) $int(xu_iu_j) \neq \emptyset$;  (iii) no vertex from $int(xu_iu_j)$ dominates $\{u_i,u_j\}$. 

Assume by symmetry that $(i,j) = (1,2)$. Since $int(C) \subseteq N_1[C]$, every vertex from $int(u_1u_2x)$ is adjacent to either $u_1$ or $u_2$. Then the subgraph $H = G[int(u_1u_2x) \cup \{x\}]$ is outerplanar and there exists a vertex $x' \in V(H) \setminus \{x\}$ with at least 2 neighbors not from $H$. Therefore, $x' \in int(xu_1u_2)$ dominates $\{u_1,u_2\}$, a contradiction.

%By Lemma~\ref{l2}, there exists a vertex $x' \in int(u_1u_2x)$ adjacent to at least $\delta(G) - 2$ vertices from $\{u_1,u_2,x\}$ and not adjacent to $x$. Then $x'$ dominates $\{u_1,u_2\}$, a contradiction.

\textbf{Case 2.} There exist vertices $x,x' \in int(C)$ with two common neighbors in $C$. Assume by symmetry that $u_1x,u_2x \in E(G)$ and $x' \in int(u_1u_2x)$. We may also assume that no vertex from $int(u_1u_2x) \setminus \{x'\}$ dominates $\{u_1,u_2\}$, hence $int(u_1u_2x') = \emptyset$ by Case~1. Moreover, $int(u_1x'u_2x) \neq \emptyset$, since $deg(x') \geq 4$. We use the notation $$U_i = N_1[u_i] \cap int(u_1x'u_2x), \ i \in \{1,2\}.$$ If, say, $U_2 = \emptyset$, then for every vertex $u'_1 \in U_1$ we have $N_1[u_1'] \subseteq N_1[u_1]$, a contradiction by Lemma~\ref{l12}. Thus we may assume that $U_1,U_2 \neq \emptyset$. Moreover, $xx' \notin E(G)$ (otherwise for every vertex $u''_1 \in U_1$ we have $N_1[u_1''] \subseteq N_1[u_1]$, again a contradiction by Lemma~\ref{l12}). Since $deg(x') \geq 4$, the following subcases are possible:

\textbf{Case 2.1.} The vertex $x'$ has at least two neighbors in $U_i$ for some $i \in \{1,2\}$. Assume by symmetry that $i = 1$. Then there exist vertices $u'_1,u''_1 \in U_1 \cap N_1[x']$ such that $u'_1 \in int(u_1x'u''_1)$ and thus $N_1[u'_1] \subseteq N_1[u_1]$, a contradiction.

\textbf{Case 2.2.} There exist vertices $u'_1 \in U_1 \cap N_1(x)$ and $u'_2 \in U_2 \cap N_1(x)$. Note that $deg(x') = 4$ and $int(u_1u'_1x)$, $int(u_2u'_2x) = \emptyset$. If $u'_1u'_2 \in E(G)$, then $N_1[x'] \cong C_4$, as required. Suppose that $u'_1u'_2 \notin E(G)$ (see Fig.~2). Since $deg(u'_1) \geq 4$, the following subcases are possible:

\textbf{Case 2.2.1.} There exist vertices $y_1,y_2 \in N_1[u'_1] \cap (U_1 \cup \{x\})$. Then either $y_1 \in int(u_1u'_1y_2)$ or $y_2 \in int(u_1u'_1y_1)$. Therefore, $N_1[y_i] \subseteq N_1[u_1]$ for some $i \in \{1,2\}$, a contradiction.

\textbf{Case 2.2.2.} There exists a vertex $z \in  N_1[u'_1] \cap (U_2 \cup \{x\})$. Then $u'_2 \in int(u_2x'u'_1z)$. Since $u'_1u'_2 \notin E(G)$, we have  $N_1[u'_2] \subseteq N_1[u_2]$, a contradiction.

\begin{center}
\includegraphics[scale=0.45]{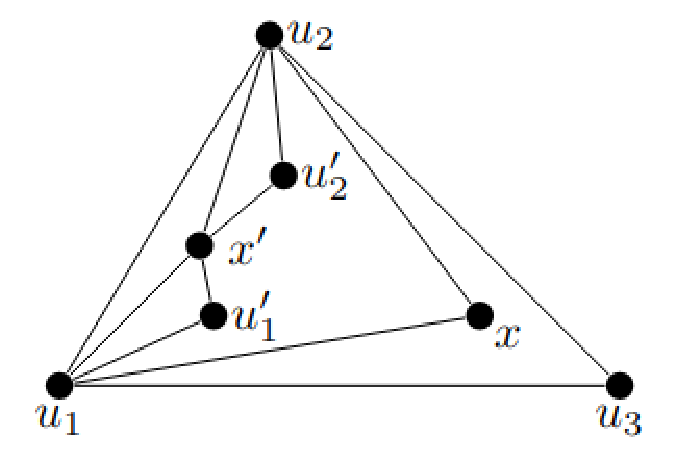}

Fig.~2. Illustration for the proof of Lemma~\ref{l15}, Case~2.2.
\end{center}

\textbf{Case 3.} There exist vertices $x,y_1,y_2 \in int(C)$ such that for some $1 \leq i < j \leq 3$ the following holds: (i) $u_ix, u_jx, u_iy_1,u_jy_2,y_1y_2 \in E(G)$;  (ii) $x \in int(u_iy_1y_2u_j)$.

Assume by symmetry that $(i,j) = (1,2)$. Cases~1 and~2 imply that $int(u_1u_2x) = \emptyset$. Note that the vertices $y_1,y_2$ may or may not have more than one neighbor in $C$. We use the notation $U_i = N_1[u_i] \cap int(u_1xu_2y_2y_1)$, where $i \in \{1,2\}$.

\textbf{Case 3.1.} $int(u_1xu_2y_2y_1) = \emptyset$.  Then $deg(x) = 4$ and $N_1[x] \cong C_4$, as required.

\textbf{Case 3.2.} $int(u_1xu_2y_2y_1) \neq \emptyset$ and $xy_i \in E(G)$ for some $i \in \{1,2\}$. Assume by symmetry that $i = 1$. Then $int(u_1xy_1) = \emptyset$, hence $U_1 = \emptyset$ and $U_2 \neq \emptyset$. If there exists a vertex $u'_2 \in U_2$ nonadjacent to $y_1$, then $N_1[u'_2] \subseteq N_1[u_2]$, a contradiction by Lemma~\ref{l12}. Otherwise $y_1$ dominates $U_2$, thus at most one vertex from $U_2$ is adjacent to $x$. Since $deg(x) = 4$, there exists such a vertex $u^*_2 \in U_2 \cap N_1(x)$, then $u_1u_2u^*_2y_1$ is a 4-cycle and $G[N_1(x)] \cong C_4$, as required.

\textbf{Case 3.3.} $int(u_1xu_2y_2y_1) \neq \emptyset$ and $xy_1,xy_2 \notin E(G)$. It is not hard to check, using the argument from Case~2, that $deg(x) = 4$ and $G[N_1(x)] \cong C_4$, as required. 

\textbf{Case 4.} There exist at most two vertices in $int(C)$ with two or three neighbors in $C$. By Lemma~\ref{l2}, there are exactly two such vertices, denote them by $x$ and $y$. Case~2 implies that neither $x$ nor $y$ has three neighbors in $C$, thus we may assume that $u_1x,u_1y,u_2x,u_3y \in E(G)$ and $u_3x,u_2y \notin E(G)$. Cases~1 and~2 imply that $int(u_1u_2x), int(u_1u_3y) = \emptyset$. We use the notation $U_i = N_1(u_i) \cap int(u_1xu_2u_3y)$, where $i \in [1;3]$. 

\textbf{Case 4.1.} $xy \in E(G)$. Lemma~\ref{l12} implies that $int(u_1xy) = \emptyset$ and thus $U_1 = \emptyset$. By Case~3, $x$ is not adjacent to $U_3$ and $y$ is not adjacent to $U_2$. Thus there exist vertices $z_2 \in U_2$ and $z_3 \in U_3$ such that $xz_2, yz_3 \in E(G)$. By Lemma~\ref{l12}, such vertices are unique.

\textbf{Case 4.1.1.} $z_2z_3 \in E(G)$. The subgraph $G[int(u_2xyu_3)]$ is outerplanar. Since $\delta(G) = 4$ and no vertex in $int(u_2xyu_3)$ dominates $\{u_2,u_3\}$, we have $|int(u_2xyu_3)| \geq 4$.  By Lemma~2, there exists a vertex $z' \in int(u_2xyu_3)$ nonadjacent to $z_2$ with at least two neighbors in $\{x,y,u_2,u_3\}$. Note that $int(u_2xz_2), int(xyz_3z_2), int(u_3z_3y) = \emptyset$, then $z'$ dominates $\{u_2,u_3\}$, a contradiction. 

\textbf{Case 4.1.2.} $z_2z_3 \notin E(G)$. If $z_2$ has a neighbor in $U_3$ then it is easy to see that $N_1[z_3] \subseteq N_1[u_3]$, a contradiction by Lemma~\ref{l12}. Otherwise there exist vertices $u'_2,u''_2 \in N_1[z_2] \cap N_1[u_2]$, thus either $N_1[u'_2] \subseteq N_1[u_2]$ or $N_1[u''_2] \subseteq N_1[u_2]$, again a contradiction by Lemma~\ref{l12}.

\textbf{Case 4.2.} $xy \notin E(G)$. Case~3 implies that $x$ is not adjacent to $U_3$. It is easy to check, using Lemma~\ref{l12}, that $x$ has at most one neighbor in $U_2$. Thus there exists a vertex $z \in N_1[x] \cap U_1$. Case~3 implies that $z$ is not adjacent to $U_2 \cup U_3$. If $zy \in E(G)$, then $int(u_1xz) = int(u_1yz) = \emptyset$, thus $deg(z) = 3$, a contradiction.  If $zy \notin E(G)$, then $z$ has at least two neighbors in $N_1[u_1]$, a contradiction by Lemma~\ref{l12}.
 
\textbf{Case 5.} There exist at least three vertices in~$int(C)$ with at most two neighbors in $C$. The previous cases imply there exist exactly three vertices $x_1,x_2,x_3 \in int(u_1u_2u_3)$ with exactly two neighbors in $C$. Assume by symmetry that $u_ix_i \notin E(G)$, $i \in [1;3]$. Case~3 implies that the set $\{x_1,x_2,x_3\}$ is independent. Cases~1 and~2 imply that $int(u_1u_2x_3)$, $int(u_2u_3x_1)$, $int(u_1u_3x_2) = \emptyset$. We use the notation $U_i = N_1[u_i] \cap int(u_1x_3u_2x_1u_3x_2)$, where $i \in [1;3]$.

\textbf{Case 5.1.}  For some $1 \leq i < j \leq 3$ the vertices $x_i$ and $x_j$ have a common neighbor $z$. Assume by symmetry that $(i,j) = (1,2)$. Case~3 implies that $z$ is not adjacent to $\{u_1,u_2\}$, hence $zu_3 \in E(G)$ and $zx_3 \notin E(G)$. Then $z$ is not adjacent to $U_1 \cup U_2$. Lemma~\ref{l12} implies that $int(u_3x_2z) = int(u_3x_1z) = \emptyset$, hence $deg(z) = 3$, a contradiction.

\textbf{Case 5.2.} The vertices $x_1,x_2,x_3$ have no common neighors in $int(C)$. Note that by Case~3, $x_1$ is not adjacent to $U_1$. Lemma~\ref{l12} implies that $x_1$ have at most one neighbor in $U_2$, thus there exist a vertex $z_1 \in N_1[x_1] \cap U_3$. Case~3 implies that $z_1$ is not adjacent to $U_1 \cup U_2 \cup \{x_3\}$. Since $deg(z_1) \geq 4$, there exist vertices $u'_3,u''_3 \in N_1[z_1] \cap U_3$. Then either $N_1[u'_3] \subseteq N_1[u_3]$ or $N_1[u''_3] \subseteq N_1[u_3]$, a contradiction by Lemma~\ref{l12}. 
\end{proof}

\subsection{Cutvertices and separating edges}\label{s43}

\begin{lemma}\label{l16}
A critical graph $G$ has no cutvertices.
\end{lemma}

\begin{proof}
Suppose that $G$ has a cutvertex $v$. Let $G_1$ be a component of $G \setminus v$ and $G_2$ be the union of all its remaining components. Let $G'_i = G[V(G_i) \cup \{v\}]$, where $i \in \{1,2\}$. Consider three cases:

\textbf{Case 1.} $\gamma(G_1) = \theta(G_1)$ and $\gamma(G_2) = \theta(G_2)$. Let $\gamma(G_1) = p$ and $\gamma(G_2) = q$. Since $G$ is critical, $ \theta(G) = p + q + 1$ and $\gamma(G) = \gamma^{\infty}(G) = p + q$. Consider a vertex $w \in V(G_1) \cap N_1[v]$. By Lemma~\ref{l8}, there exists a set $D_1 \in DS(G_1)$ such that $w \in D_1$, thus $\gamma(G'_1) = p$. Since $\theta(G) > p + q$, we have $\theta(G'_1) = p + 1$. Moreover, since $G$ is critical, we have $\gamma^{\infty}(G'_1) = p + 1$. Likewise, $\gamma(G_2) = q$ and $\gamma^{\infty}(G'_2) = \theta(G'_2) = q + 1$. Thus the graphs $G_1,G'_1,G_2,G'_2$ are maximum-demand, but the graph $G$ is not maximum-demand, a contradiction by Lemma~\ref{l4}.

\textbf{Case 2.} $\gamma(G_1) = \theta(G_1)$ and $\gamma(G_2) < \theta(G_2)$. Since $G$ is critical, $\gamma(G_2) < \gamma^{\infty}(G_2)$. Consider a vertex $w \in V(G_1) \cap N_1[v]$. By Lemma~\ref{l8}, there exists a set $D_1 \in DS(G_1)$ such that $w \in D_1$. Then for every set $D_2 \in MDS(G_2)$ we have $D_1 \cup D_2 \in MDS(G)$. Therefore, $$\gamma(G) \leq \gamma(G_1) + \gamma(G_2) < \gamma^{\infty}(G_1) + \gamma^{\infty}(G_2) = \gamma^{\infty}(G_1 \cup G_2) \leq \gamma^{\infty}(G)$$
and $G$ is not critical, a contradiction.

\textbf{Case 3.} $\gamma(G_1) < \theta(G_1)$ and $\gamma(G_2) < \theta(G_2)$. Clearly, $\gamma(G) \leq \gamma(G_1) + \gamma(G_2) + 1$. Note that for both $i \in \{1,2\}$ the subgraph $G_i$ is not critical, thus $\gamma(G_i) < \gamma^{\infty}(G_i)$. Therefore, 
$$\gamma(G) \leq \gamma(G_1) + \gamma(G_2) + 2 < \gamma^{\infty}(G_1) + \gamma^{\infty}(G_2) \leq \gamma^{\infty}(G).$$
Again, $G$ is not critical, a contradiction.
\end{proof}

\begin{lemma}\label{l17}
If $\delta(G) = 5$, then $G$ has no separating edges.
\end{lemma}

\begin{proof}
Suppose for a contradiction there exists a separating edge $ab \in E(G)$. Let $G_1$ be a component of the subgraph $G[V(G) \setminus \{a,b\}]$ and $G_2$ be the union of all its remaining components. For $i \in \{1,2\}$ let $$G^a_i = G[V(G_i) \cup \{a\}], \ G^b_i = G[V(G_i) \cup \{b\}], \  G^{a,b}_i = G[V(G_i) \cup \{a,b\}].$$ Suppose that $\gamma(G_1) < \gamma^{\infty}(G_1)$ and $\gamma(G_2) < \gamma^{\infty}(G_2)$, then  $$\gamma(G) \leq \gamma(G_1) + \gamma(G_2) + 1 \leq \gamma^{\infty}(G_1) + \gamma^{\infty}(G_2) - 1 \leq  \gamma^{\infty}(G) - 1.$$ Thus $G$ is not critical, a contradiction. Therefore, there are two remaining variants:

\textbf{Variant 1.} Either $\gamma(G_1) < \gamma^{\infty}(G_1)$ or $\gamma(G_2) < \gamma^{\infty}(G_2)$. Assume by symmetry that $\gamma(G_1) < \gamma^{\infty}(G_1)$ and $\gamma(G_2) = \theta(G_2)$. Then, since $\gamma(G) = \gamma^{\infty}(G)$, we have $\gamma(G_1) = \gamma^{\infty}(G_1) - 1$. Therefore, there exist integers $p,q \geq 1$ such that $$\gamma(G_1) = p, \ \gamma^{\infty}(G_1) = \theta(G_1) = p + 1;$$ $$\gamma(G_2) = \gamma^{\infty}(G_2) = \theta(G_2) = q;$$ $$\gamma(G) = \gamma^{\infty}(G) = p + q + 1, \ \theta(G) = p + q + 2.$$

Since $G$ has no cutvertices, both sets $A_2 = N_1(a) \cap V(G_2)$ and $B_2 = N_1(b) \cap V(G_2)$ are nonempty. If there exists a vertex $w \in A_2 \cap B_2$, then by Lemma~\ref{l8} there exists a set $D$ such that $w \in D \in MDS(G_2)$. Then $D \in MDS(G^{a,b}_2)$ and $$\gamma(G) \leq \gamma(G_1) + \gamma(G^{a,b}_2) \leq \gamma(G_1) + \gamma(G_2) = p + q,$$ a contradiction. Therefore, $A_2 \cap B_2 = \emptyset$. We now consider four cases.

\textbf{Case 1.} There exists a set $D_1 \in MEDS(G)$ and an integer $s \in \{1,2\}$ such that $|\{a,b\} \cap D_1| = s $ and $|D_1 \cap V(G_2)| \leq q-s$. Assume by symmetry that $a \in D_1$. 

\textbf{Case 1.1.} $\gamma^{\infty}(G_2 \setminus A_2) < \gamma^{\infty}(G_2)$. Note that $\gamma(G^{a,b}_2) \geq \gamma(G) - \gamma(G_1) = q + 1$, then $$\gamma(G_2 \setminus A_2) \leq \gamma^{\infty}(G_2 \setminus A_2) \leq \gamma(G^{a,b}_2)- 2.$$ However, for every set $D \in MDS(G_2 \setminus A_2)$ we have $D \cup \{a\} \in DS(G^{a,b}_2)$. Hence $\gamma(G \setminus A_2) \geq \gamma(G^{a,b}_2) - 1$, a contradiction.

\textbf{Case 1.2.} $\gamma^{\infty}(G_2 \setminus A_2) = \gamma^{\infty}(G_2)$. Let $D'_1 = D_1 \cap (V(G^b_2))$. Note that $$N_1[V(G_2) \setminus A_2] \subseteq V(G^b_2) \ \text{and} \ |D'_1| < \gamma^{\infty}(G_2 \setminus A_2).$$ Then there exists an AS $\mathfrak{B}$ with the vertices from $V(G_2) \setminus A_2$ such that the strategy $(G_2^b,D_1',\mathfrak{B})$ is winning. Therefore, the strategy $(G,D_1,\mathfrak{B})$ is also winning, a contradiction. 

\textbf{Case 2.} There exists a set $D_2 \in MEDS(G)$ such that $a,b \notin D_2$ and $|D_2 \cap V(G_2)| < q$. Then there exists an AS $\mathfrak{B}$ such that the strategy $(G_2,D_2 \cap V(G_2),\mathfrak{B})$ is winning. Clearly, the strategy $(G,D_2,\mathfrak{B})$ is also winning, a contradiction.

\textbf{Case 3.} There exists a set $D_3 \in MEDS(G)$ such that $\{a,b\} \subseteq D$ and $D_3 \cap V(G_2) = q - 1$. Then $|D_3 \cap V(G_1)| = p$ and $D_3 \cap V(G_1) \notin EDS(G_1)$. Hence one can attack the vertices of $V(G_1)$ in such a way that, at some point, one of the guards located on a vertex from $\{a,b\}$ moves to a vertex from $V(G_1)$. Therefore, there exists a set $D^*_3 \in MEDS(G)$ such that $|D^*_3 \cap V(G_1)| = p + 1$ and $|D^*_3 \cap V(G_2)| = q - 1$. Since $|D^*_3 \cap \{a,b\}| = 1$, Case~1 implies that $D^*_3 \notin MEDS(G)$, a contradiction.

\textbf{Case 4.} For every set $D \in MEDS(G)$ we have $|D \cap V(G_2)| \geq q$. Lemma~\ref{l7}(a) implies that $\gamma^{\infty}(G_1^{a,b}) \leq p + 1$. If $\gamma(G_1^{a,b}) \leq p$, then $\gamma(G) \leq p + q$, a contradiction. If $\gamma(G_1^{a,b}) = p + 1$, then, since $G$ is critical, $\theta(G_1^{a,b}) = p + 1$. Hence  $\theta(G) \leq p + q + 1$, again a contradiction. This completes the proof of Variant~1.

\textbf{Variant 2.} There exist integers $p,q \geq 1$ such that $$\gamma(G_1) = \theta(G_1) = p, \ \gamma(G_2) = \theta(G_2) = q.$$ Clearly, $\theta(G) \leq p + q + 1$. Since $G$ is critical, it  is easy to see that

$$\gamma(G) = \gamma^{\infty}(G) = p + q, \ \theta(G) = p + q + 1.$$
Lemma~\ref{l6} implies that $\gamma^{\infty}(G_1 \cup G_2) \leq \gamma^{\infty}(G \setminus a) \leq \gamma^{\infty}(G)$. Therefore, $\gamma^{\infty}(G \setminus a) = p + q$. Likewise, $\gamma^{\infty}(G \setminus b) = p + q$. It is not hard to see that $$\gamma(G \setminus a) \geq \min(\gamma(G),\gamma(G_1 \cup G_2)),$$ then $\gamma(G \setminus a) = p + q$. Likewise, $\gamma(G \setminus b) = p + q$. Then, since $G$ is critical, $$ \theta(G \setminus a) = \theta(G \setminus b) = p + q.$$

\begin{claim}\label{cl1}
For some $i \in \{1,2\}$ there exists a family $\mathcal{C}' \in MCP(G_i)$ and a clique $W \in \mathcal{C}'$ such that the set $W \cup \{a\}$ is a clique in $G$.
\end{claim}

\begin{proof}
Consider a family $\mathcal{C} \in MCP(G \setminus b)$ and its clique $W' \ni a$. If $W' = \{a\}$, then $\theta(G \setminus b) = \theta(G)$, a contradiction. Thus $|W'| \geq 2$. Since $V(G_1) \cap V(G_2) = \emptyset$, there exists $i \in \{1,2\}$ such that $W' \setminus \{a\} \subseteq V(G_i)$ (assume by symmetry that $i = 1$). Note that $\mathcal{C}$ has exactly $\theta(G_1)$ cliques with the vertices from $G^a_1$. Remove from $\mathcal{C}$ all cliques with the vertices from $V(G_2)$  and replace the clique $W'$ with~$W = W' \setminus \{a\}$ to obtain a required family $\mathcal{C}' \in MCP(G_1)$. 
\end{proof}

For the rest of the proof of Variant~2 we assume that Claim~\ref{cl1} holds with $i = 1$.

\begin{claim}\label{cl2}
There exists a family $\mathcal{C}'' \in MCP(G_1)$ such that $W \notin \mathcal{C}''$.
\end{claim}

\begin{proof}
Suppose for a contradiction that $W$ belongs to every MCP of $G_1$. It follows from Claim~\ref{cl1}  that for some $j \in \{1,2\}$ there exists a family $\mathcal{C}^* \in MCP(G_j)$ with a clique $Q \in \mathcal{C}^*$ such that the set $Q \cup \{b\}$ is also a clique. If $j = 2$, then  $\theta(G) \leq |\mathcal{C}'| + |\mathcal{C}^*| = p + q$, a contradiction. Suppose that $j = 1$. Since $W \in \mathcal{C}^*$, either $W = Q$ or $W \cap Q = \emptyset$. If $W = Q$, then $W \cup \{a,b\}$ is a clique, thus $\theta(G_1) = \theta(G^{a,b}_1)$ and $\theta(G) = p + q$, a contradiction. If $W \cap Q = \emptyset$, then $(W \cup \{a\}) \cap (Q \cup \{b\}) = \emptyset$, and $\theta(G_1) = \theta(G^{a,b}_1)$, again a contradiction.
\end{proof}

For the rest of the proof of Variant~2, consider families $\mathcal{C}',\mathcal{C}'' \in MCP(G_1)$ and a clique $W \subseteq V(G_1)$ such that $W \in \mathcal{C}'$, $W \notin \mathcal{C}''$ and $W \cup \{a\}$ is a clique in $G$. Since $\gamma(G_1) = \theta(G_1)$, the clique $W$ is maximal by inclusion in $V(G_1)$.

\begin{claim}\label{cl3}
$|W| = 3$.
\end{claim}

\begin{proof}
Since $G$ is planar and $W \cup \{a\}$ is a clique, $|W| \leq 3$. Note that $\gamma(G_1) = \theta(G_1)$ and $\delta(G_1) \geq 3$. Remind that $W$ is maximal by inclusion in $G_1$. If $W = \{x\}$, then $x$ is isolated in $G_1$, hence $\delta(G) < 5$, a contradiction. Suppose that $W = \{x_1,x_2\}$.  By Lemma~\ref{l1}, there exist nonadjacent vertices $x_1' \in N_1(x_1) \cap V(G_1)$ and $x_2' \in N_1(x_2) \cap V(G_1)$. Then the set $\{x'_1,x'_2\} \in IS(G_1)$ dominates $W$, a contradiction by Lemma~\ref{l9}.
\end{proof}

Let $W = \{x_1,x_2,x_3\}$. Assume by symmetry that $x_3 \in int(ax_1x_2)$, thus $bx_3 \notin E(G)$.

\begin{claim}\label{cl4}
The vertices $x_1$ and $x_2$ have no common neighbors in $V(G_1) \setminus \{x_3\}$. 
\end{claim}

\begin{proof}
Suppose for a contradiction there exists a vertex $y \in V(G_1) \setminus \{x_3\}$ adjacent to both $x_1$ and $x_2$. The clique $\{x_1,x_2,x_3\}$ is maximal by inclusion in $G_1$, thus $yx_3 \notin E(G)$. Let $X_3 = N_1[x_3] \setminus \{x_1,x_2,a\}$. Suppose there exists a vertex $z \in X_3$ such that the set $\{y,z\}$ is $\theta$-independent in $\mathcal{C'}$. Then $\{y,z\}$ dominates $W$ in $G_1$, a contradiction by Lemma~\ref{l9}. Therefore, $\{y\} \cup X_3$ is a clique. Since $deg_G(x_3) \geq 5$ and $x_3b \notin E(G)$, we have $|X_3| = 2$. Let $X_3 = \{x'_3,x''_3\}$, then some clique from $\mathcal{C}'$ contains $\{y,x'_3,x''_3\}$. 

\textbf{Case 4.1.} $\{x'_3,x''_3,y\} \in \mathcal{C}'$. Consider a vertex $z \in N_1[y] \setminus \{x_1,x_2,x'_3,x''_3\}$ (such a vertex exists, since $ya,yb \notin E(G)$). The set $\{z,x_3\} \in IS(G_1)$ dominates $\{x'_3,x''_3,y\}$, a contradiction by Lemma~\ref{l9}. 

\textbf{Case 4.2.} There exists a vertex $y' \in N_1(y)$ such that $\{x'_3,x''_3,y,y'\} \in \mathcal{C}'$. Since $G$ is planar, $y' \in int(x'_3x''_3y)$.  Consider a clique $W'$ such that $x_3 \in W' \in \mathcal{C}''$. Since $W \neq W'$, there exists $i \in \{1,2\}$ such that $x_i \notin W'$. Note that $W' \subseteq \{x_1,x_2,x_3,x'_3,x''_3\}$ and $x_iy' \notin E(G)$, thus the set $\{x_i,y'\} \in IS(G_1)$ dominates $W'$, a contradiction by Lemma~\ref{l9}. 
\end{proof}

\begin{claim}\label{cl5}
Neither $x_1$ nor $x_2$ has a common neighbor with $x_3$ in $V(G_1) \setminus \{x_1,x_2\}$.
\end{claim}

\begin{proof}
Suppose for a contradiction that for some $i \in \{1,2\}$ the vertices $x_i$ and $x_3$ have a common neighbor $y \in V(G_1) \setminus \{x_1,x_2\}$. Assume by symmetry that $i = 1$. Let $X_2 = N_1[x_2] \setminus \{x_1,x_3,a,b\}$. Lemma~\ref{l9} implies that the set $X_2 \cup \{y\}$ is a clique, thus $|X_2| \leq 2$.  Note that $X_2$ is nonempty, since $deg_{G_1}(x_2) \geq 3$. Hence $y$ is a unique common neighbor of $x_1$ and $x_3$ in $V(G_1) \setminus \{x_2\}$.

\textbf{Case 1.} $X_2 = \{x'_2\}$ for some $x'_2 \in V(G_1)$. Consider the clique $W'$ such that $x_2 \in W' \in \mathcal{C}''$. Claim~\ref{cl4} implies that $W' = \{x_2,x'_2\}$ or $W' = \{x_2,x'_2,x_3\}$. The set $\{x_1,y\}$ dominates $W$, thus by Lemma~\ref{l9} it is not $\theta$-independent in $\mathcal{C}''$ and there exists a clique $W''$ such that $\{x_1,y\} \subseteq W'' \in \mathcal{C}''$. Since $G$ is planar and $|N_1(x'_2) \setminus \{x_2,x_3,y\}| \geq 2$, there exists a vertex $z \in N_1(x'_2) \setminus \{x_2,x_3,y\}$ such that $z \notin W''$. Therefore, the set $\{x_1,z\}$ is $\theta$-independent in $\mathcal{C}''$ and dominates $W$, a contradiction by Lemma~\ref{l9}.

\textbf{Case 2.} $X_2 = \{x'_2,x''_2\}$ for some $x'_2,x''_2 \in V(G_1)$. Assume by symmetry that $x'_2 \in int(x_1x_2x''_2y)$, thus $x_1x''_2,x_3x'_2 \notin E(G)$. By Lemma~\ref{l9}, both sets $\{y,x'_2\}$ and $\{y,x''_2\}$ are not $\theta$-independent in $\mathcal{C}'$, hence there exists a clique $W'$ such that $\{y,x'_2,x''_2\} \subseteq W' \in \mathcal{C}'$. 

\textbf{Case 2.1.} $W' = \{y,x'_2,x''_2\}$. Consider a vertex $z \in N_1(y) \setminus \{x_1,x_3,x'_2,x''_2\}$. The set $\{x_2,z\} \in IS(G_1)$ dominates $\{y,x'_2,x''_2\}$, a contradiction by Lemma~\ref{l9}. 

\textbf{Case 2.2.} There exists a vertex $y' \in V(G_1)$ such that $\{y,y',x'_2,x''_2\} \in \mathcal{C}'$ (see Fig.~3). Then $y' \in int(x'_2x''_2y)$. Consider the clique $W''$ such that $x_2 \in W'' \in \mathcal{C}''$. Claim~\ref{cl4} implies that $W'' \subseteq \{x'_2,x''_2,x_1,x_2,x_3\}$. There exists $i \in \{1,3\}$ such that $x_i \notin W''$. Hence the set $\{x_i,y'\} \in IS(G_1)$ dominates $W''$, a contradiction by Lemma~\ref{l9}.
\end{proof}

\begin{center}
\includegraphics[scale=0.45]{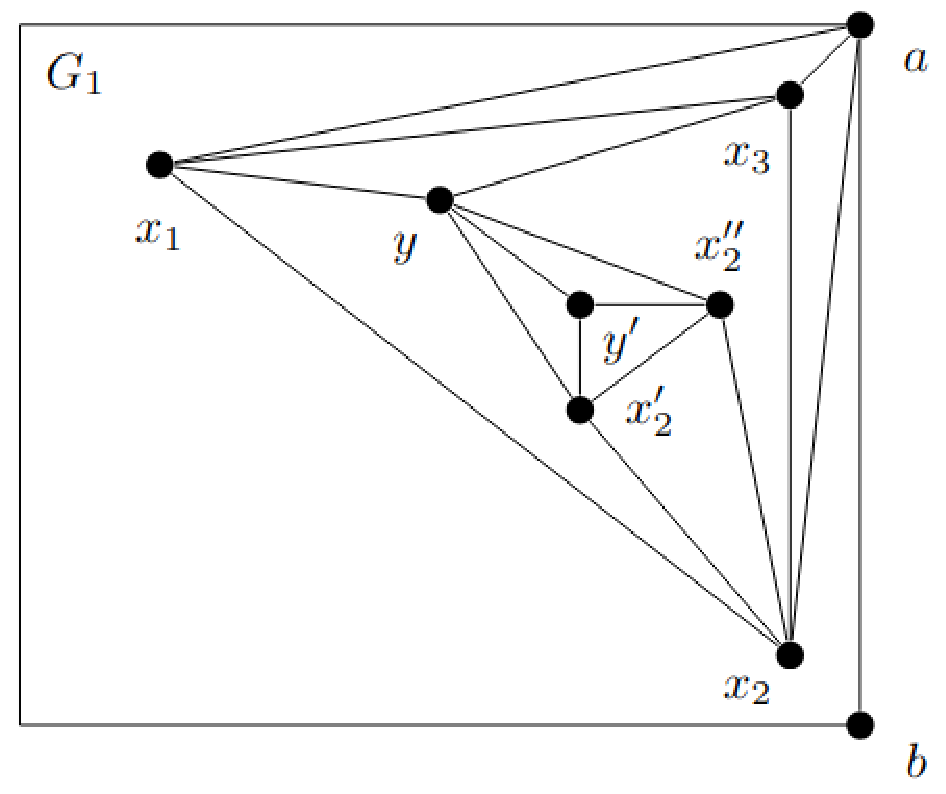}

Fig. 3. Illustration for the proof of Claim~\ref{cl5}, Case~2.2
\end{center}

\begin{claim}\label{cl6}
There exist vertices $y_1,y_2,y_3 \in V(G_1)$ such that $\{x_3,y_1,y_2,y_3\} \in \mathcal{C}''$. 
\end{claim}

\begin{proof}
Suppose there are no such vertices. Claims~\ref{cl4},~\ref{cl5} and Lemma~\ref{l9} imply that for all $1 \leq i < j \leq 3$ the set $\{x_i,x_j\}$ is $\theta$-independent in $\mathcal{C}''$. Then the following cases are possible:

\textbf{Case 1.} There exists a vertex $y \in V(G_1)$ such that $\{x_3,y\} \in \mathcal{C}''$. Consider a vertex $y' \in N_1[y] \setminus \{x_1,x_2,x_3\}$. Claims~\ref{cl4},~\ref{cl5} imply that $y'x_i \notin E(G)$ for some $i \in \{1,2\}$. Then the set $\{y',x_i\} \in IS(G_1)$ dominates $\{x_3,y\}$, a contradiction by Lemma~\ref{l9}.

\textbf{Case 2.} There exist vertices $y_1,y_2 \in V(G)$ such that $\{x_3,y_1,y_2\} \in \mathcal{C}''$. Claims~\ref{cl4},~\ref{cl5} imply that $y_ix_j \notin E(G)$ for all $i,j \in \{1,2\}$. Consider the sets $Y_i = N_1(y_i) \setminus N_1(x_3)$, $i \in \{1,2\}$. Note that $|Y_1|,|Y_2| \geq 3$.

\textbf{Case 2.1.} There exists a vertex $w \in Y_1 \cap Y_2$. Then for some $i \in \{1,2\}$ the set $\{w,x_i\}$ is independent and dominates $\{x_3,y_1,y_2\} \in \mathcal{C}''$, a contradiction by Lemma~\ref{l9}.

In cases~2.2 and~2.3 we assume that $Y_1 \cap Y_2 = \emptyset$.

\textbf{Case 2.2.} For some $i \in \{1,2\}$ there exists a vertex $w_i \in Y_i$ not adjacent to $\{x_1,x_2\}$. Assume by symmetry  that $i = 1$. Since $G$ is planar and $|Y_2| \geq 3$, the set $\{w_1\} \cup Y_2$ is not a clique. Thus there exists a vertex $w_2 \in Y_2$ such that the set $\{w_1,w_2\}$ is $\theta$-independent in $\mathcal{C}''$. Claims~\ref{cl4},~\ref{cl5} imply that there exists $j \in \{1,2\}$ such that $w_2x_j \notin E(G)$. Then the set $\{x_j,w_1,w_2\}$ is $\theta$-independent in $\mathcal{C}''$ and dominates $\{y_1,y_2,x_3\} \in \mathcal{C}''$, a contradiction.

\textbf{Case 2.3.} Every vertex from $Y_1 \cup Y_2$ is adjacent to either $x_1$ or $x_2$. Consider a subset $Z \subseteq Y_1 \cup Y_2$ of vertices $z$ such that the set $\{x_1,z\}$ is $\theta$-independent in $\mathcal{C}''$. Note that $|Z| \geq |Y_1| + |Y_2| - 2$. Assume that $|Y_1 \cap Z| \leq |Y_2 \cap Z|$. If $|Y_1 \cap Z| \geq 2$, then by Lemma~\ref{l1} there exist nonadjacent vertices $w_1 \in Y_1$ and $w_2 \in Y_2$, thus the set $\{x_1,w_1,w_2\}$ is $\theta$-independent in $\mathcal{C}''$ and dominates $\{y_1,y_2,x_3\}$, a contradiction. Otherwise $|Y_1 \cap Z| = 1$ and $|Y_2 \cap Z| \geq 3$. Let $ Y_1 \cap Z = \{w_1\}$. There exists a vertex $w_2 \in Y_2 \cap Z$ such that the set $\{w_1,w_2\}$ is $\theta$-independent in $\mathcal{C}''$. Thus $\{x_1,w_1,w_2\}$ is also $\theta$-independent in $\mathcal{C}''$ and dominates $\{y_1,y_2,x_3\}$, a contradiction.
\end{proof}

We are now ready to finish the proof of Variant~2. Assume by symmetry that $y_3 \in int(y_1y_2x_3)$, then $x_1y_3, x_2y_3 \notin E(G)$. Let $X_i = N_1[x_i] \setminus \{x_1,x_2,x_3,a,b\}$, $i \in \{1,2\}$. Assume that $|X_1| \geq |X_2| \geq 1$.

\textbf{Case 1.} $|X_2| = 1$. Let $X_2 = \{x'_2\}$. Consider the clique $W'$ such that $x_2 \in W' \in \mathcal{C}''$. Since the set $\{x_1,x_2,x_3\}$ is $\theta$-independent in $\mathcal{C}''$, we have $W' = \{x_2,x'_2\}$. Consider a vertex $z \in N_1(x'_2) \setminus \{x_1,x_2,x_3\}$. Claim~\ref{cl5} implies that there exists $i \in \{1,3\}$ such that the set $\{x_i,z\}$ is $\theta$-independent and dominates $W'$, a contradiction by Lemma~\ref{l9}.

\textbf{Case 2.} $|X_2| \geq 2$. By Lemma~\ref{l1} there exists a pair of nonadjacent vertices $x'_1 \in X_1$ and $x'_2 \in X_2$. Claims~\ref{cl4} and~\ref{cl5} imply that $x'_i \neq y_j$ for all $i,j \in \{1,2\}$. Hence the set $\{x'_1,x'_2,y_3\} \in IS(G)$ dominates the clique $\{x_1,x_2,x_3\} \in \mathcal{C}'$, a contradiction by Lemma~\ref{l9}. The proof of Lemma~\ref{l17} is complete.
\end{proof}

\section{Critical graphs with $\delta \leq 4$}

\subsection{Case $\delta \leq 3$}

\begin{lemma}\label{l18}
If a graph $G$ is critical, then $\delta(G) \geq 3$.
\end{lemma}

\begin{proof}
Suppose that there exists a vertex $v \in V(G)$ such that $deg(v) \leq 2$. By Lemma~\ref{l12}, the set $N_1[v]$ is not a clique, then  $deg(v) = 2$ and the neighbors of $v$ (denoted by $u_1$ and $u_2$) are nonadjacent. Corollary~\ref{cor1}(b) implies that $N_1(u_1) = N_1(u_2)$ (if, for example, there exists a vertex $u'_1 \in N_1(u_1) \setminus N_1(u_2)$, then $\gamma(G_{u_2}) \leq \gamma(G_{v,u_2})$, a contradiction by lemmas~\ref{l13} and~\ref{l14}). Let $deg(u_1) = 1 + k$ and $N_1[u_1] \setminus \{v\} = \{w_1,\ldots,w_k\}$. Assume that the vertices $w_1,\ldots,w_k$ are in clockwise order around $u_1$, then $w_iw_j \notin E(G)$ for all $1 \leq i \leq j + 2 \leq k$.

\textbf{Case 1.} $k = 1$. By Lemma~\ref{l5} we have $|V(G)| > 4$. Then $w_1$ is a cutvertex, a contradiction by Lemma~\ref{l16}.

\textbf{Case 2.} $k = 2$. Suppose that $w_1w_2 \notin E(G)$. Lemma~\ref{l10}(c) implies that every vertex $x \in V(G) \setminus N_1[v]$ has a neighbor in $\{w_1,w_2\}$ (otherwise $\gamma(G_x) < \theta(G_x)$, a contradiction by Lemma~\ref{l13}). Thus $\{v,w_1,w_2\} \in DS(G)$ and $\gamma(G) \leq 3$. By Lemma~\ref{l13}, $\theta(G_{u_1,u_2}) \leq 1$. Since $G$ is $K_5$-free, we have $|V(G)| \leq 9$, a contradiction by Lemma~\ref{l5}. 

Suppose now that $w_1w_2 \in E(G)$. By Lemma~\ref{l13}, $\gamma(G_v) = \theta(G_v) = \gamma(G) - 1$ and $\gamma(G_{u_1,u_2}) = \theta(G_{u_1,u_2}) = \gamma(G) - 2$. However, for every family $\mathcal{C} \in MCP(G_{u_1,u_2})$ we have $\mathcal{C} \cup \{\{v,u_1\},\{u_2,w_1,w_2\}\} \in CP(G)$, thus $\theta(G) \leq \theta(G_{u_1,u_2}) + 2 = \gamma(G)$, a contradiction.

\textbf{Case 3.} $k = 3$. By Lemma~\ref{l10}(c), every vertex from $V(G) \setminus N_2[v]$ has a neighbor in $\{w_1,w_2,w_3\}$, thus $\{v,w_1,w_2,w_3\} \in DS(G)$ and $\gamma(G) \leq 4$. If $\gamma(G) \leq 3$, then $\gamma(G_{u_1,u_2}) = \theta(G_{u_1,u_2}) \leq 1$ and $|V(G)| \leq 10$, a contradiction by Lemma~\ref{l5}. Suppose that $\gamma(G) = 4$, then there exists a vertex $x_2 \in N_1(w_2)$ not adjacent to $\{w_1,w_3\}$ (if there is no such vertex, then $\{u_1,w_1,w_3\} \in DS(G)$, a contradiction). Since $w_1w_3 \notin E(G)$, we have $\gamma(G_{x_2}) < \theta(G_{x_2})$ by Lemma~\ref{l10}(c), a contradiction.

\textbf{Case 4.} $k \geq 4$. If $N_3(v) = \emptyset$, then $\gamma(G) = |\{v,u_1\}| < |\{v,w_1,w_3\}| \leq \alpha(G)$, hence~$G$ is not critical, a contradiction. Suppose that $N_3(v) \neq \emptyset$. Lemma~\ref{l10}(c) implies that for every vertex $x \in N_3(v)$ the set $\{w_1,\ldots,w_k\} \setminus N_1[x]$ is a clique. Since $G$ is $K_{3,3}$-free, every vertex from $N_3(v)$ has at most two neighbors in $\{w_1,\ldots,w_k\}$. Therefore, $k = 4$ and every vertex from $N_3(v)$  dominates either $\{w_1,w_2\}$ or $\{w_3,w_4\}$. Assume by symmetry that there exists a vertex $x \in N_3(v)$ that dominates $\{w_1,w_2\}$. Then $w_3w_4 \in E(G)$ by Lemma~\ref{l10}(c) and thus $N_1[w_4] \subseteq N_1[w_3]$, a contradiction by Lemma~\ref{l12}.

\end{proof}

\begin{lemma}\label{l19}
If a graph $G$ is critical, then $\delta(G) \geq 4$.
\end{lemma}

\begin{proof}
Suppose for a contradiction that $\delta(G) \leq 3$, then $\delta(G) = 3$ by the previous lemma. Let $v \in V(G)$ be a vertex of degree 3 and $N_1(v) = \{u_1,u_2,u_3\}$. Consider three cases depending on the structure of the subgraph $G[N_1(v)]$.

\textbf{Case 1.} $G[N_1(v)] \cong 3K_1$. By Lemma~\ref{l13},

\begin{equation}\label{eq6}
\gamma(G) - 2 \leq \gamma(G_{v,u_i}) < \gamma(G_{u_i}) = \gamma(G) - 1 \text{ for all } i \in \{1,2,3\}.
\end{equation}

Corollary~\ref{cor1}(b) implies that no vertex from $N_2(v)$ has exactly two neighbors in $N_1(v)$. 

\textbf{Case 1.1.} There exists a vertex $w \in N_2(v)$ that dominates $N_1(v)$ (since $G$ is planar, such a vertex is unique). Then there exist pairwise distinct vertices $u'_i \in N_1(u_i) \setminus \{v,w\}$, where $i \in \{1,2,3\}$. Assume by symmetry that $u_2 \in int(vu_1wu_3)$ and $u'_2 \in int(vu_2wu_3)$. Then the set $\{u'_1,u'_2\} \in IS(G_{v,u_3})$ dominates $\{u_1,u_2\}$, a contradiction by Corollary~\ref{cor1}(a). 

\textbf{Case 1.2.} Every vertex from $N_2(v)$ has a unique neighbor in $N_1(v)$. By Lemma~\ref{l1} there exist nonadjacent vertices $u_1'' \in N_1(u_1) \setminus \{v\}$ and $u_2'' \in N_1(u_2) \setminus \{v\}$. Again, the set $\{u_1'',u_2''\} \in IS(G_{v,u_3})$ dominates $\{u_1,u_2\}$, a contradiction by Corollary~\ref{cor1}(a).

\textbf{Case 2.} $G[N_1(v)] \cong K_2 \cup K_1$. Assume that $u_1u_2 \in E(G)$ and $u_2u_3,u_1u_3 \notin E(G)$. Note that by Lemma~\ref{l14} the property (\ref{eq6}) still holds. Moreover, every neighbor of $u_3$ is adjacent to both $u_1$ and $u_2$ (if, for example, there exists a vertex $x \in N_2(v)$ such that $xu_3 \in E(G)$ and $xu_1 \notin E(G)$, then $x$ dominates $N_1[v] \setminus N_1[u_1] = \{u_3\}$ in $G_{u_1}$, a contradiction by Corollary~\ref{cor1}(b)). Since $G$ is $K_{3,3}$-free, we have $deg(u_3) \leq 2$, a contradiction.

\textbf{Case 3.} $G[N_1(v)] \cong P_3$ or $G[N_1(v)] \cong C_3$. Clearly, for some $i \in \{1,2,3\}$ we have $N_1[v] \subseteq N_1[u_i]$, a contradiction by Lemma~\ref{l12}.
\end{proof}

\subsection{Case $\delta = 4$}

Throughout this subsection we consider a planar graph $G$ with a vertex $v$ such that $deg(v) = \delta(G) = 4$. Let $N_1(v) = \{u_1,u_2,u_3,u_4\}$ and $S \subseteq N_1(v)$. Call a vertex $x \in N_2(v)$ \textit{$S$-vertex}, if $N_1(x) \cap N_1(v) = S$. We use the notation $U_i =  N_1(u_i) \setminus N_1(v)$. Since $G$ is critical, the sets $U_i$ are nonempty by Lemma~\ref{l12}.

It is easy to check that precisely one of the following possibilities holds:

\begin{enumerate}
\item $G[N_1(v)]$ has a triangle (Lemma~\ref{l20});
\item $G[N_1(v)]$ has a vertex of degree 3 (Lemma~\ref{l21});
\item $G[N_1(v)] \cong C_4$ (Lemma~\ref{l22});
\item $G[N_1(v)] \cong 4K_1$ (Lemma~\ref{l23}); 
\item $G[N_1(v)] \cong P_3 \cup K_1$ (Lemma~\ref{l24}); 
\item $G[N_1(v)] \cong K_2 \cup 2K_1$ (Lemma~\ref{l25}); 
\item $G[N_1(v)] \cong 2K_2$ (Lemma~\ref{l26}); 
\item $G[N_1(v)] \cong P_4$. (Lemma~\ref{l27});
\end{enumerate}

\begin{lemma}\label{l20}
If $G[N_1(v)]$ has a triangle, then $G$ is not critical.
\end{lemma}

\begin{proof} We may assume that $G[N_1(v)]$ has a cycle $u_1u_2u_3$ such that $u_3 \in ext(vu_1u_2)$ and $u_4 \in int(vu_1u_2)$. Select a vertex $x \in U_4$. Since $x \in int(vu_1u_2)$, the vertex $u_3$ is not adjacent to $\{u_4,x\}$. Lemmas~\ref{l13} and~\ref{l14} imply that $$\gamma(G) - 2 = \gamma(G_{v,u_3}) = \theta(G_{v,u_3}) < \gamma(G_{u_3}) = \theta(G_{u_3}) = \gamma(G) - 1.$$ Note that $V(G_{u_3}) = V(G_{v,u_3}) \cup \{u_4\}$. Since $\gamma(G_{v,u_3}) = \theta(G_{v,u_3})$, Lemma~\ref{l8} implies that there exists a set $D \in MDS(G_{v,u_3})$ such that $x \in D$.  Hence $D \in DS(G_{u_3})$ and $\gamma(G_{u_3}) \leq \gamma(G_{v,u_3})$, a contradiction.
\end{proof}

\begin{lemma}\label{l21}
If $G[N_1(v)]$ has a vertex of degree 3, then $G$ is not critical.
\end{lemma}

\begin{proof}
There exists $i \in [1;4]$, such that $N_1[v] \subseteq N_1[u_i]$, a contradiction by Lemma~\ref{l12}.
\end{proof}

\begin{lemma}\label{l22}
If $G[N_1(v)] \cong C_4$, then $G$ is not critical.
\end{lemma}

\begin{proof}
We may assume that $G$ has a cycle $u_1u_2u_3u_4$, then $u_1u_3,u_2u_4 \notin E(G)$ by Lemma~\ref{l20}. First, we prove the following fact.

\begin{claim*} For both $i \in \{1,2\}$ we have $U_i \subseteq U_{i+2}$ or $U_{i+2} \subseteq U_i$.
\end{claim*}

\begin{proof}
Suppose for a contradiction that, for example, there exist vertices $u'_1 \in U_1 \setminus U_3$ and $u_3' \in U_3 \setminus U_1$.  If, moreover, $$\gamma(G \setminus \{u_1,u_2,v,u_4\}) =  \theta(G \setminus \{u_1,u_2,v,u_4\}) = \gamma(G) - 1,$$ then by Lemma~\ref{l8} there exists a set $D \in MDS(G \setminus \{u_1,u_2,v,u_4\})$ such that $\{u'_1,u_3\} \subseteq D$. Since $D \in MDS(G)$, we obtain a contradiction. Thus by Lemma~\ref{l14} we have $\gamma^{\infty}(G_{u_1,v}) = \gamma^{\infty}(G) - 2$. Lemma~\ref{l13} implies that $$\gamma(G) - 2 = \gamma(G_{u_1,v}) = \theta(G_{u_1,v}) < \gamma(G_{u_1}) = \theta(G_{u_1}) = \gamma(G) - 1.$$ By Lemma~\ref{l8}, there exists a set  $D' \in MDS(G_{u_1,v})$ such that $u'_3 \in D'$. Clearly $D' \in DS(G_{u_1})$, hence $\gamma(G_{u_1,v}) \geq \gamma(G_{u_1})$, a contradiction. 
\end{proof}

We are now ready to finish the proof of Lemma~\ref{l22}. Assume that $U_1 \subseteq U_3$ and $U_2 \subseteq U_4$. Let $w \in N_2(v)$ be a common neighbor of $u_1$ and $u_3$. Since $G$ is planar and $\delta(G) = 4$, $w$ is the unique common neighbor of $u_2$ and $u_4$, thus $deg(u_1) = deg(u_2) = 4$.  Consider an induced $(n-6)$-vertex subgraph $H = G[V(G) \setminus (N_1[v] \cup w)]$. Lemmas~\ref{l6} and~\ref{l13} imply that $\gamma^{\infty}(H) \leq \gamma^{\infty}(G_v) = \gamma^{\infty}(G) - 1$. Consider two cases:

\textbf{Case 1.} $\gamma^{\infty}(H) \leq \gamma^{\infty}(G) - 2$, then $\gamma(H) \leq \gamma(G) - 2$. For every set $D \in MDS(H)$ we have $D \cup \{v,w\} \in DS(G)$, thus $\gamma(H) = \gamma^{\infty}(H) = \gamma(G) - 2$. Since $G$ is critical, $$\gamma(H) = \theta(H) \leq \theta(G) - 3.$$ Consider a family $\mathcal{C} \in MCP(H)$. Clearly, $$\mathcal{C} \cup \{\{u_1,u_2,v\} \cup \{u_3,u_4,w\}\} \in CP(G).$$  Therefore, $\theta(G) \leq \theta(H) + 2$, a contradiction.

\textbf{Case 2.} $\gamma^{\infty}(H) = \gamma^{\infty}(G) - 1$. By Lemma~\ref{l11}, $\gamma^{\infty}(G \setminus \{u_1,u_2,v\}) = \gamma^{\infty}(G)$, then by Lemma~\ref{l7}(b) there exists a set $D'' \in MEDS(G)$ such that $$D'' \cap \{v,u_1,u_2\} = \emptyset, \ |D'' \cap \{v,u_1,u_2,u_3,u_4,w\}| \leq 1.$$

It is easy to see that $D''$ does not dominate $\{u_1,u_2,v\}$, a contradiction.
\end{proof}

\begin{lemma}\label{l23}
If $G[N_1(v)] \cong 4K_1$, then $G$ is not critical.
\end{lemma}

\begin{proof}
Consider three cases: 

\textbf{Case 1.} There exists a vertex $w \in N_2(v)$ with at least three neighbors in $N_1(v)$. Corollary~\ref{cor1}(b) implies that $w$ is a $N_1(v)$-vertex. Since $G$ is planar, such a vertex is unique. We may assume that the vertices $u_1,u_2,u_3,u_4$ are in clockwise order around $v$, then $u_2 \in int(vu_1wu_3)$ and $u_4 \in ext(vu_1wu_3)$. Hence $int(vu_1wu_2) \cup int(vu_2wu_3) \neq \emptyset$. By Lemma~\ref{l15}, there exists a vertex $x \in int(vu_1wu_2) \cup int(vu_2wu_3)$ not adjacent to $N_1[v] \cup w$. Then $\gamma(G_x) < \theta(G_x)$ by Lemma~\ref{l10}(e), a contradiction.

\textbf{Case 2.} There exists a vertex $x \notin N_1[v]$ with at most one neighbor in $N_1(v)$. We may assume that $x$ is not adjacent to $\{u_2,u_3,u_4\}$. Since $\gamma(G_x) = \theta(G_x)$, Lemma~\ref{l10}(a,e) implies that exactly one vertex from $N_1(v)$ is pendant in $G_x$. Then for some $k \in [1;4]$ the vertex $x$ dominates the set $U_k$ and $xu_k \notin E(G)$ (see Fig.~4). Since $|U_k| \geq 3$, there exist vertices $y_1,y_2,y_3 \in U_k$ such that  $int(u_ky_1xy_3) \cap U_k = \{y_2\}$. Since $deg(y_2) \geq 4$, Lemma~\ref{l15} implies that there exists a vertex $z \in int(u_ky_1xy_2) \cup int(u_ky_2xy_3)$ not adjacent to $N_2[v]$. Then $\gamma(G_z) < \theta(G_z)$ by Lemma~\ref{l10}(e), again a contradiction.

\begin{center}
\includegraphics[scale=0.45]{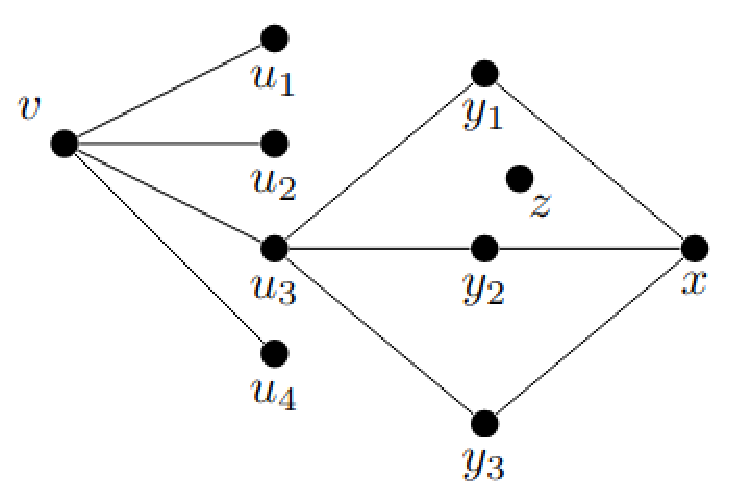}

Fig.~4. Illustration for the proof of Lemma~\ref{l23}, Case~2 (here $k = 3$).
\end{center}

\textbf{Case 3.} Every vertex from $V(G) \setminus N_1[v]$ has exactly two neighbors in $N_1(v)$. Since $G$ is planar, there exists a pair of vertices from $N_1(v)$ (say, $(u_1,u_2)$) with no common neighbors in $N_2(v)$. Then it is easy to see that $\{v,u_3,u_4\} \in DS(G)$. Note that $N_1(v) \in IS(G)$, hence $\gamma(G) < \alpha(G) \leq \gamma^{\infty}(G)$, a contradiction.
\end{proof}

\begin{lemma}\label{l24}
If  $N_1(v) \cong P_3 \cup K_1$,  then $G$ is not critical.
\end{lemma}

\begin{proof}
Assume by symmetry that $E(G[N_1(v)]) = \{u_1u_2,u_2u_3\}$ and $u_3 \in ext(vu_1u_2)$. By Corollary~\ref{cor1}(b), $G$ has neither $\{u_1,u_4\}$-vertices nor $\{u_3,u_4\}$-vertices. Consider two cases:

\textbf{Case 1.} There exists a vertex $w \in N_2(v)$ with at least three neighbors in $N_1(v)$. Corollary~\ref{cor1}(b) implies that  $w$ is a $\{u_1,u_3,u_4\}$-vertex or a $\{u_1,u_2,u_3,u_4\}$-vertex. Since $G$ is planar, such a vertex is unique and $w \in ext(vu_1u_2u_3)$. Hence $u_4 \in ext(vu_1u_2u_3)$ and we may assume that the vertices $u_1,\ldots,u_4$ are in clockwise order around $v$. 

\textbf{Case 1.1.} There exists a $\{u_1,u_3\}$-vertex $x$. Note that $x \in int(vu_1wu_3)$, then $u_2w \notin E(G)$ and $U_2 \subseteq int(vu_1xu_3)$. Consider a vertex $u_2' \in U_2$. Corollary~\ref{cor1}(b) implies that $u_2'$ is not a $\{u_1,u_2,u_3\}$-vertex. If $u'_2u_1,u'_2u_3 \notin E(G)$, then $\gamma(G_{u_2'}) < \theta(G_{u_2'})$ by Lemma~\ref{l10}(e), a contradiction. Otherwise $u'_2$ is a $\{u_2,u_i\}$-vertex for some $i \in \{1,3\}$, assume by symmetry that $i = 1$. If there exists a vertex $u'_4 \in U_4 \setminus U_3$, then the set $\{u'_2,u'_4\} \in IS(G_{v,u_3})$ dominates the set $\{u_1,u_2,u_4\}$, a contradiction by Corollary~\ref{cor1}(a). Otherwise $U_4 \subseteq U_3$, hence the vertices $u_3$ and $u_4$ have at least three common neighbors in $G_{u'_2}$ and $\gamma(G_{u'_2}) < \theta(G_{u'_2})$ by Lemma~\ref{l10}(d), again a contradiction.

\textbf{Case 1.2.} There are no $\{u_1,u_3\}$-vertices. Since $\delta(G) \geq 4$, both sets $U_1 \setminus U_3$ and $U_3 \setminus U_1$ are nonempty. Consider a vertex $u_4' \in U_4$. Corollary~\ref{cor1}(b) implies that $u'_4$ is not adjacent to $\{u_1,u_3\}$. Assume by symmetry that $u'_4 \in int(vu_3wu_4)$, then there exists a vertex $u_1' \in U_1 \setminus U_3$, nonadjacent to $u_4'$. The set $\{u_1',u_4'\} \in IS(G_{v,u_3})$  dominates the set $\{u_1,u_4\}$, a contradiction by Corollary~\ref{cor1}(a).

\textbf{Case 2.} Every vertex from $N_2(v)$ has at most two neighbors in $N_1(v)$. Recall that $G$ has neither $\{u_1,u_4\}$-vertices nor $\{u_3,u_4\}$-vertices. Then every vertex from $U_4$ is a $\{u_2,u_4\}$-vertex or a $\{u_4\}$-vertex. Consider three subcases:

\textbf{Case 2.1.} All vertices from $U_4$ are $\{u_2,u_4\}$-vertices. Since $|U_4| \geq 3$, there exist $\{u_2,u_4\}$-vertices $x_1,x_2,x_3$ such that $ int(u_2x_1u_4x_3) \cap U_4 = \{x_2\}$. Since $x_1,x_3 \notin U_1 \cup U_3$, we have $N_1[x_2] \cap (U_1 \cup U_3) = \emptyset$. 

\textbf{Case 2.1.1.} $U_1 \neq U_3$. We may assume that there exists a vertex $u'_1 \in U_1 \setminus U_3$, then the set $\{u'_1,x_2\} \in IS(G_{v,u_3})$ dominates $\{u_1,u_2,u_4\}$, a contradiction by Corollary~\ref{cor1}(a).

\textbf{Case 2.1.2.} $U_1 = U_3$. Then $u_4 \in int(vu_1u_2u_3)$. Consider a vertex $y \in U_1 \cap U_3$. Since $y \in ext(vu_1u_2u_3)$ and $U_4 \in int(vu_1u_2u_3)$, the vertices $u_2$ and $u_4$ have at least $|U_4| \geq 3$ common neighbors in $G_{y}$, thus $\gamma(G_{y}) < \theta(G_{y})$ by Lemma~\ref{l10}(d), a contradiction.
 
\textbf{Case 2.2.} There exists a $\{u_2,u_4\}$-vertex $x$ and a $\{u_4\}$-vertex $u'_4$. 

\textbf{Case 2.2.1.} $u_4 \in ext(vu_1u_2u_3)$. Assume by symmetry that $u'_4 \in int(vu_2xu_4)$. Consider a vertex $u'_1 \in U_1$. The set $\{u'_1,u'_4\} \in IS(G_{v,u_3})$ dominates $\{u_1,u_4\}$, a contradiction by Corollary~\ref{cor1}(a). 

\textbf{Case 2.2.2.} $u_4 \in int(vu_1u_2u_3)$. Assume by symmetry that $u_4 \in int(vu_2u_3)$.

\textbf{Case 2.2.2.1.} There exists a vertex $u'_1 \in U_1 \setminus U_3$. Since $u'_1 \in ext(vu_2u_3)$, we have $xu'_1 \notin E(G)$. The set $\{u'_1,x\} \in IS(G_{v,u_3})$ dominates $\{u_1,u_4\}$, a contradiction by Corollary~\ref{cor1}(a). 

\textbf{Case 2.2.2.2.} There exists a vertex $u'_3 \in U_3 \setminus U_1$. Remind that $u'_3 \notin U_4$ by Corollary~\ref{cor1}(b). Since $|U_4| \geq 3$, there exists a vertex $u^*_4 \in U_4$ such that the set $\{u'_3,u^*_4\}$ is $\theta$-independent in $G_{v,u_1}$. Since $\{u'_3,u^*_4\}$ dominates $\{u_3,u_4\}$, we obtain a contradiction by Corollary~\ref{cor1}(a). 

\textbf{Case 2.2.2.3.} $U_1 = U_3$. By Lemma~\ref{l10}, the equality $\gamma(G_x) = \theta(G_x)$ implies that $|U_1| = 2$ and the set $\{u_1\} \cup U_1$ is a clique. Let $U_1 = \{y,z\}$. Lemma~\ref{l13} implies that $\gamma(G_{u_1,u_3}) = \gamma(G) - 2$ and $\theta(G_{u_1,u_3}) \leq \theta(G) - 3$.  Note that $$V(G) = V(G_{u_1,u_3}) \cup \{v,u_1,u_2,u_3,y,z\}.$$ Consider a family $\mathcal{C} \in MCP(G_{u_1,u_3})$. Clearly, $$\mathcal{C} \cup \{\{v,u_2,u_3\},\{u_1,y,z\}\} \in CP(G).$$ Hence $\theta(G) \leq \theta(G_{u_1,u_3}) + 2$, a contradiction.

\textbf{Case 2.3.} All elements of $U_4$ are $\{u_4\}$-vertices. 

\textbf{Case 2.3.1.} $U_1 \neq U_3$. Assume by symmetry that there exists a vertex $u_1' \in U_1 \setminus U_3$. Since $|U_4| \geq 3$, the set $\{u_1'\} \cup U_4$ is not a clique, thus there exists a vertex $u'_4 \in U_4$ such that the set $\{u'_1,u'_4\}$ is $\theta$-independent in $G_{v,u_3}$, a contradiction by Corollary~\ref{cor1}(a). 

\textbf{Case 2.3.2.} $U_1 = U_3$. There exist $\{u_1,u_3\}$-vertices $y_1, y_2$ such that $y_2 \in ext(vu_1y_1u_3)$.  Consider a vertex $u_2' \in U_2$.  The set $\{u'_2,y_2\} \in IS(G_{v,u_4})$ dominates $\{u_1,u_2\}$, again a contradiction by Corollary~\ref{cor1}(a).

\end{proof}

\begin{lemma}\label{l25}
If $N_1(v) \cong K_2 \cup 2K_1$, then $G$ is not critical.
\end{lemma}
\begin{proof}

Assume by symmetry that $u_1u_2 \in E(G)$ and $u_3 \in ext(vu_1u_2)$. Corollary~\ref{cor1}(b) implies that every vertex from $U_3 \cap U_4$ is a $N_1(v)$-vertex and no vertex from $N_2(v)$ has exactly three neighbors in $N_1(v)$. Consider two cases:

\textbf{Case 1.} There exists a vertex $w \in U_1 \cap U_2$.

\textbf{Case 1.1.} $w$ is a $N_1(v)$-vertex. Since $G$ is planar, such a vertex is unique. Moreover, $w,u_4 \in ext(vu_1u_2)$. We may assume that the vertices $u_1,\ldots,u_4$ are in clockwise order around $v$, then $U_2 \cup U_3 \subseteq int(vu_1wu_4)$. 

\textbf{Case 1.1.1.} There are no $\{u_1,u_4\}$-vertices. Note that  $U_3 \cap U_4 = \{w\}$, hence $|U_3 \setminus \{w\}|,|U_4 \setminus \{w\}| \geq 2$. By Lemma~\ref{l1}, there exist nonadjacent vertices $u_3' \in U_3 \setminus U_1$ and $u_4' \in U_4 \setminus U_1$. Then the set $\{u_3',u_4'\} \in IS(G_{v,u_1})$ dominates the set $\{u_3,u_4\}$,  a contradiction by Corollary~\ref{cor1}(a).

\textbf{Case 1.1.2.} There exists a $\{u_1,u_4\}$-vertex $x \in V(G)$ (see Fig.~5). It is easy to check, using Corollary~\ref{cor1}(a), that $U_2 = U_3$. Since $|U_3 \setminus \{w\}| \geq 2$, there exists a $\{u_2,u_3\}$-vertex $y$ nonadjacent to $w$. Hence $y$ is not adjacent to $U_1 \cup U_4$ and $\gamma(G_y) < \theta(G_y)$ by Lemma~\ref{l10}(d), a contradiction. 

\begin{center}
\includegraphics[scale=0.45]{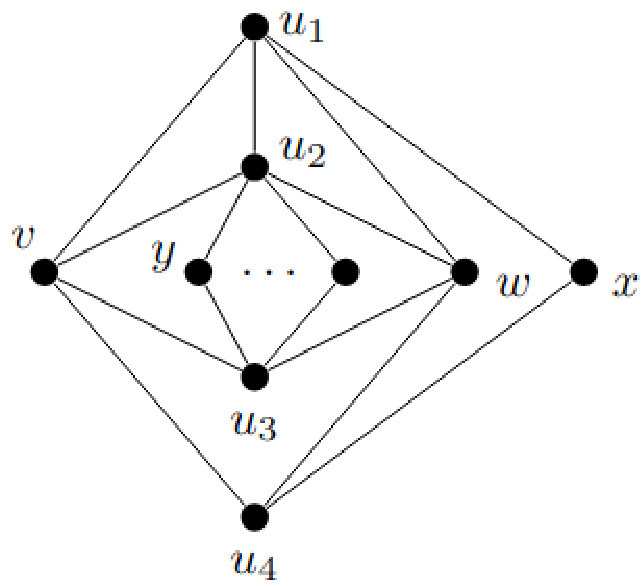}

Fig.~5. Illustration for the proof of Lemma~\ref{l25}, Case 1.1.2.
\end{center}

\textbf{Case 1.2.} $w$ is not a $N_1(v)$-vertex, then it is a $\{u_1,u_2\}$-vertex. Then $U_3 \cap U_4 = \emptyset$. Since $|U_3| \geq 3$, the set $\{w\} \cup U_3$ is not a clique, hence there exists a vertex $u'_3 \in U_3$ such that the set $\{w,u'_3\}$ is $\theta$-independent in $G_{v,u_4}$ and dominates $\{u_1,u_2,u_3\}$, a contradiction by Corollary~\ref{cor1}(a).

\textbf{Case 2.} $U_1 \cap U_2 = \emptyset$. By Lemma~\ref{l1}, there exist nonadjacent vertices $u'_1 \in U_1$ and $u'_2 \in U_2$. Corollary~\ref{cor1}(b) implies that these vertices have at most two neighbors in $N_1(v)$. The following subcases are possible: 

\textbf{Case 2.1.} One of the vertices $u_3$ and $u_4$ (say, $u_3$), is not adjacent to $\{u'_1,u'_2\}$. If $u'_1u_4 \in E(G)$ or $u'_2u_4 \in E(G)$, then $\{u'_1,u'_2\} \in IS(G_{v,u_3})$ dominates $\{u_1,u_2,u_4\}$, a contradiction by Corollary~\ref{cor1}(a). Therefore, $u'_iu_{j+2} \notin E(G)$, for all $i,j \in \{1,2\}$. By Lemma~\ref{l10}(e), in both $G_{u_1'}$ and $G_{u_2'}$ exactly one of the vertices $u_3$ or $u_4$ is pendant. Since $|U_3|,|U_4| \geq 3$, if $u_3$ is pendant in $G_{u_1'}$ then $u_4$ is pendant in $G_{u_2'}$, and vice versa. Therefore, both vertices $u_3$ and $u_4$ are pendant in $G_{u'_1,u'_2}$, a contradiction by Lemma~\ref{l10}(a).

In the subcases~2.2--2.4 we may assume that both $u_3$ and $u_4$ have exactly one neighbor in $\{u'_1,u'_2\}$.

\textbf{Case 2.2.} $u_4 \in int(vu_1u_2$). Assume by symmetry that $u'_1u_4, u'_2u_3 \in E(G)$. 

\textbf{Case 2.2.1.} $U_4 \subseteq U_1$. Since $u'_2 \in ext(vu_1u_2)$ and $U_4 \subseteq int(vu_1u_2)$, we have $\gamma(G_{u'_2}) < \theta(G_{u'_2})$ by Lemma~\ref{l10}(d), a contradiction.

\textbf{Case 2.2.2.} There exists a vertex $u'_4 \in U_4 \setminus U_1$. Since $u'_4 \in int(vu_1u_2)$ and $u'_2 \in ext(vu_1u_2)$, we have $u'_2u'_4 \notin E(G)$. Hence the set $\{u'_2,u'_4\} \in IS(G_{v,u_1})$ dominates the set$\{u_3,u_4\}$, a contradiction by Corollary~\ref{cor1}(a).

In the subcases~2.3--2.4  we assume that $u_3,u_4 \in ext(vu_1u_2)$. We may also assume that the vertices $u_1,u_2,u_3,u_4$ are located clockwise around $v$. Then the only possible configuration is $u'_1u_4, u'_2u_3 \in E(G)$.

\textbf{Case 2.3.} $U_1 \cap U_4 = \{u'_1\}$. Corollary~\ref{cor1}(a) implies that every $\{u_2,u_3\}$-vertex dominates $U_4 \setminus \{u'_1\}$, thus $u'_2$ is a unique such vertex in $G$. Consider vertices $u'_3 \in U_3 \setminus \{u'_2\}$ and $u'_4 \in U_4 \setminus \{u'_1\}$. Clearly, $u'_3 \in int(vu_2u'_2u'_4u_4)$ and $u'_1u'_3 \notin E(G)$ (see Fig.~6). Thus the set $\{u'_1,u'_3\} \in IS(G_{v,u_2})$ dominates the set $\{u_1,u_3,u_4\}$, a contradiction by Corollary~\ref{cor1}(a).

\begin{center}
\includegraphics[scale=0.45]{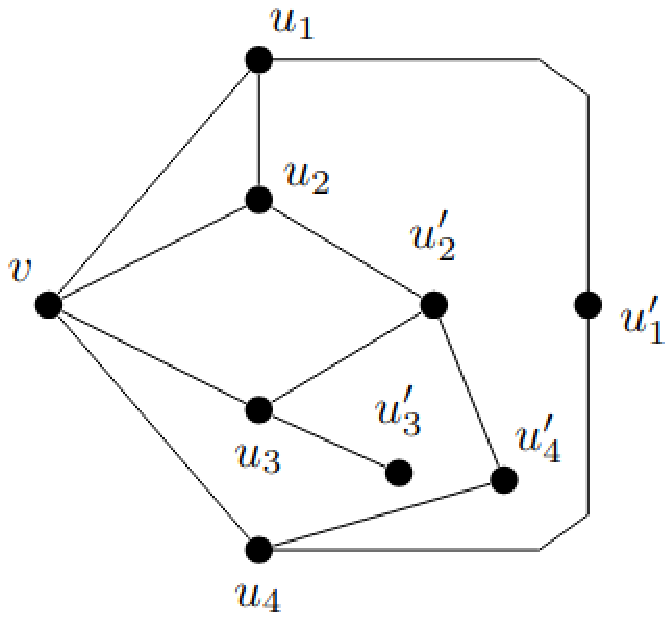}

Fig.~6. Illustration for the proof of Lemma~\ref{l25}, Case 2.3.
\end{center}

\textbf{Case 2.4.} There exists a vertex $x \in (U_1 \cap U_4) \setminus \{u'_1\}$. We may assume that $x \in int(vu_1u'_1u_4)$ (otherwise rename the vertices $x$ and $u'_1$). Since $\gamma(G_{u'_1}) = \theta(G_{u'_1})$, by Lemma~\ref{l10}(b) there exists a vertex $y \in (U_2 \cap U_3) \setminus \{u'_2\}$, assume that $y \in int(vu_2u'_2u_3)$. If there exists a vertex $z \in U_4 \setminus U_1$, then the independent set $\{y,z\}$ dominates $\{u_2,u_3,u_4\}$ in $G_{u_1}$, a contradiction by Corollary~\ref{cor1}(a). Otherwise $U_4 \subseteq U_1$, hence $\gamma(G_{y}) < \theta(G_{y})$ by Lemma~\ref{l10}(d), a contradiction.

\end{proof}

\begin{lemma}\label{l26}
If $N(v) \cong 2K_2$, then $G$ is not critical.
\end{lemma}
\begin{proof}

Assume by symmetry that $u_1u_2,u_3u_4 \in E(G)$. We may also assume that $u_3,u_4 \in ext(vu_1u_2)$ and the vertices $u_1,\ldots,u_4$ are in clockwise order around $v$. Corollary~\ref{cor1}(b) implies that no vertex from $N_2(v)$ has exactly three neighbors in $N_1(v)$. Consider two cases:

\textbf{Case 1.} There exists a $N_1(v)$-vertex $w$. 

\textbf{Case 1.1.} There exists a $\{u_2,u_3\}$-vertex $y$.  It is easy to check, using Corollary~\ref{cor1}(a), that $U_1 = U_4$. Lemma~\ref{l10}(b) implies that there exists a $\{u_1,u_4\}$-vertex $x \in ext(vu_1wu_4)$. Likewise, it is easy to check, using Corollary~\ref{cor1}(a), that $U_2 = U_3$. 

\textbf{Case 1.1.1.} $min(deg(u_1),deg(u_2)) = 4$. Assume by symmetry that $deg(u_1) = 4$. Then $wx \in E(G)$ (otherwise $N_1[u_1] \cong P_3 \cup K_1$, a contradiction by Lemma~\ref{l24}). Note that $deg(x) \geq 4$ and $x$ is not a cutvertex by Lemma~\ref{l16}. Then by Lemma~\ref{l15} there exists a vertex $z \in int(u_1xw)\cup int(u_4xw)$ not adjacent to $N_1[v] \cup \{w,x,y\}$. Consider a family $\mathcal{C} \in MCP(G_z)$ and its clique $W \ni v$. Recall that $W$ is maximal by inclusion, thus $|W| = 3$ and either $W = \{v,u_1,u_2\}$ or $W = \{v,u_3,u_4\}$ (assume that $W = \{v,u_1,u_2\}$). If $\{w,u_3\}$ or $\{w,u_4\}$ is $\theta$-independent in $\mathcal{C}$, we obtain a contradiction by Lemma~\ref{l9}. Otherwise, since $U_3 \cap U_4 = \{w\}$, we have $\{w,u_3,u_4\} \in \mathcal{C}$. Then the set $\{x,v\} \in IS(G_z)$ dominates $\{w,u_3,u_4\}$, a contradiction by Lemma~\ref{l9}.

\textbf{Case 1.1.2.} $min(deg(u_1),deg(u_2)) \geq 5$. Then $G$ has at least two $\{u_1,u_4\}$-vertices and at least two $\{u_2,u_3\}$-vertices. Hence there exists a vertex $z \in U_2 \cap U_3$ not adjacent to $U_1$. Since $|U_1| \geq 3$, we have $\gamma(G_{z}) < \theta(G_{z})$ by Lemma~\ref{l10}(d), a contradiction. 

\textbf{Case 1.2.} There are no $\{u_2,u_3\}$-vertices. Assume by symmetry that there are no $\{u_1,u_4\}$-vertices.

\textbf{Case 1.2.1.} $deg(w) = 4$. By Lemma~\ref{l13}, $\gamma(G_{v,w}) = \gamma(G) - 2$ and $\theta(G_{v,w}) \leq \theta(G) - 3$. Consider a family $\mathcal{C} \in MCP(G_{v,w})$. Clearly, $$\mathcal{C} \cup \{\{v,u_1,u_2\} \cup \{u_3,u_4,w\}\} \in CP(G).$$ Hence $\theta(G) \leq \theta(G_{v,w}) + 2$, a contradiction.

\textbf{Case 1.2.2.} $deg(w) \geq 5$. 

\textbf{Case 1.2.2.1.}
$int(vu_2wu_3) \cup ext(vu_1wu_4) \neq \emptyset$. Assume that $int(vu_2wu_3) \neq \emptyset$. By Lemma~\ref{l16}, $w$ is not a cutvertex, thus for some $i \in \{2,3\}$ there exists a vertex $u'_i \in U_i \cap int(vu_2wu_3)$. Assume that $i = 2$, consider  a vertex $u'_1 \in U_1$. Remind that $u'_2u_3 \notin E(G)$, then the set $\{u'_1,u'_2\} \in IS(G_{v,u_4})$ dominates $\{u_1,u_2\}$, a contradiction. 
 
\textbf{Case 1.2.2.2.} $int(vu_2wu_3) \cup ext(vu_1wu_4) = \emptyset$. For some $j \in \{1,3\}$ there exists a vertex $x \in N_1[w] \cap int(u_ju_{j+1}w)$. Assume that $j = 3$. By Lemma~\ref{l15} the set $int(vu_1u_2) \cup int(wu_1u_2)$ contains a vertex $y$ not adjacent to $N_1[v] \cup \{w,x\}$. Consider a family $\mathcal{C} \in MCP(G_y)$ with a clique $W \ni v$. Assume by symmetry that $W = \{v,u_1,u_2\}$. Note that $U_3 \cap U_4 = \{w\}$. If $\{w,u_3\}$ or $\{w,u_4\}$ is $\theta$-independent in $\mathcal{C}$, we obtain a contradiction by Lemma~\ref{l9}. Otherwise  $\{w,u_3,u_4\} \in \mathcal{C}$ and the set $\{v,x\} \in IS(G_x)$ dominates $\{w,u_3,u_4\}$ in $G_x$, a contradiction. 

\textbf{Case 2.} All vertices of $N_2(v)$ have at most two neighbors in $N_1(v)$. Corollary~\ref{cor1}(b) implies that there are no $\{u_i,u_{i+1}\}$-vertices, where $i \in \{1,3\}$. By Lemma~\ref{l1} there exist nonadjacent vertices $u'_1 \in U_1$ and $u'_2 \in U_2$. Then by Corollary~\ref{cor1}(a), each of the vertices $u_3$ and $u_4$ is adjacent to either $u'_1$ or $u'_2$. Since both $u'_1$ and $u'_2$ have at most two neighbors in $N_1(v)$, the only possible configuration is $u'_2u_3, u'_1u_4 \in E(G)$.

\textbf{Case 2.1.} $U_1 \neq U_4$ or $U_2 \neq U_3$. We may assume that there exists a vertex $u''_1 \in U_1 \setminus U_4$. Then $u''_1u'_2$ by Corollary~\ref{cor1}(a). Consider a vertex $u''_2 \in U_2 \setminus \{u'_2\}$, then $u''_1u''_2 \in E(G)$ by Corollary~\ref{cor1}(a) and $u''_2u_3, u''_2u'_1 \notin E(G)$. Hence $\{u'_1,u''_2\} \in IS(G_{v,u_3})$ dominates $\{u_1,u_2,u_4\}$, a contradiction.

\textbf{Case 2.2.} $U_1 = U_4$ and $U_2 = U_3$. Consider the subgraph $H = G[V(G) \setminus N_2[v]]$. Note that $G_{u_1,u_3} \cong H$ and $G_{u_1,u_2} \cong K_2 \cup H$, thus $\gamma(G_{u_1,u_3}) = \gamma(G_{u_1,u_2}) - 1$.  However, Lemmas~\ref{l13} and~\ref{l14} imply that $\gamma(G_{u_1,u_2}) = \gamma(G_{u_1,u_3}) = \gamma(G) - 2$, a contradiction. 
\end{proof}

\begin{lemma}\label{l27}
If $N(v) \cong P_4$, then $G$ is not critical.
\end{lemma}

\begin{proof}

Assume that $u_1u_2,u_2u_3,u_3u_4 \in E(G)$ and $u_1 \in ext(vu_2u_3)$. Corollary~\ref{cor1}(b) implies that $G$ has no $S$-vertices, where $$S \in \{ \{u_1,u_2\}, \{u_3,u_4\}, \{u_1,u_2,u_3\}, \{u_2,u_3,u_4\}\}.$$ 

In the following cases (except subcases~3.3 and~4.3) we assume that the vertices $u_1$ and $u_4$ have a common neighbor in $N_2(v)$. Thus we may also assume that $u_4 \in ext(vu_2u_3)$ and the vertices $u_1,\ldots,u_4$ are in clockwise order around $v$.

\textbf{Case 1.} There exists a $N_1(v)$-vertex $w$. Lemma~\ref{l22} implies that $deg(u_2) \geq 5$, then $int(vu_1wu_3) \setminus \{u_2\} \neq \emptyset$. If $deg(w) = 4$, then for every family $\mathcal{C} \in MCP(G_{v,w})$ we have $$\mathcal{C} \cup \{\{v,u_1,u_2\} \cup \{w,u_3,u_4\}\} \in CP(G),$$ thus  $\theta(G_{v,w}) \geq \theta(G) - 2$, a contradiction by Lemma~\ref{l13}. Hence $deg(w) = 5$ and the set $Y = N_1(w) \setminus N_1(v)$ is nonempty. Lemma~\ref{l15} implies that $int(vu_1wu_3) \setminus (N_2[v] \cup Y) \neq \emptyset$, since one of the cycles $vu_1u_2$, $vu_2u_3$, $u_1u_2w$, $u_2u_3w$ is separating.

\textbf{Case 1.1.} There exists a vertex $x \in V(G) \setminus (N_2[v] \cup Y)$ that does not dominate the set $Y$. Select a vertex $z \in Y \setminus N_1(x)$, consider a family $\mathcal{C} \in MCP(G_x)$ and its clique $W \ni v$. If $W = \{v,u_2,u_3\}$, then $\{u_1,u_4\} \in IS(G_x)$ dominates $W$, a contradiction. Otherwise we may assume that $W = \{v,u_1,u_2\}$. If $\{w,u_3\}$ or $\{w,u_4\}$ is $\theta$-independent in $G_x$, we obtain a contradiction by Lemma~\ref{l9}. Otherwise $\{w,u_3,u_4\} \in \mathcal{C}$ (since $G$ has no $\{u_3,u_4\}$-vertices, this clique is maximal by inclusion). The set $\{v,z\} \in IS(G_x)$ dominates $\{w,u_3,u_4\}$, a contradiction by Lemma~\ref{l9}.

\textbf{Case 1.2.} Every vertex from $V(G) \setminus (N_2[v] \cup Y)$ dominates $Y$. Then  all vertices from $Y$ either belong to the set $int(u_iu_{i+1}w)$ for some $i \in \{1,2,3\}$, or belong to the set $int(vu_1wu_4)$. Lemma~\ref{l15} implies that $int(vu_iu_{i+1}) = \emptyset$ for all $i \in \{1,2,3\}$. Note that $deg(u_2)$, $deg(u_3) \geq 5$ by Lemma~\ref{l22}, hence the only possible configuration is $int(u_1u_2w) = int(u_3u_4w) = \emptyset$ and $Y \subseteq int(u_2u_3w)$. Moreover, every vertex from $ext(vu_1wu_4)$ is adjacent to $u_1$ or $u_4$ (otherwise it does not belong to $N_2[v] \cup Y$ and does not dominate $Y$, a contradiction). Corollary~\ref{cor1}(a) implies that $U_1 = U_4$ (if, say, there exists a vertex $u'_1 \in U_1 \setminus U_4$, then for any vertex $u'_2 \in U_2$ the set $\{u'_1,u'_2\} \in IS(G_{v,u_4})$ dominates $\{u_1,u_2\}$, a contradiction). Since $U_1 \setminus \{w\} \neq \emptyset$, there exists a vertex $y \in ext(vu_1wu_4)$ such that $ext(vu_1yu_4) = \emptyset$. Clearly, $deg(y) \leq 3$, a contradiction.

\textbf{Case 2.} There exists a vertex $w \in N_2(v)$ with exactly three neighbors from $N_1(v)$. Assume by symmetry that $w$ is a $\{u_1,u_3,u_4\}$-vertex. Since $G$ is planar and has no $\{u_1,u_2,u_3\}$-vertices, all vertices from $N_2(v) \setminus \{w\}$ have at most two neighbors in $N_1(v)$.

\textbf{Case 2.1.} There are no $\{u_1,u_4\}$-vertices. It is easy to verify, using Corollary~\ref{cor1}(a), that the set $(U_1 \cup U_2) \setminus \{w\}$ is a clique on at most three vertices. Remind that $G$ has no $\{u_1,u_2\}$-vertices, hence $min(deg(u_1),deg(u_2)) = 4$. If $deg(u_1) = 4$, then the subgraph $G[N_1[u_1]]$ has at most two edges, a contradiction by Lemmas~\ref{l20}--\ref{l26}. Otherwise $deg(u_2) = 4$. Let $U_2 = \{u'_2\}$. Lemmas~\ref{l20}--\ref{l26} imply that $u'_2u_3 \in E(G)$. Consider a vertex $u'_4 \in U_4 \setminus \{w\}$. The set $\{u'_2,u'_4\} \in IS(G_{v,u_1})$ dominates $\{u_3,u_4\}$, a contradiction by Corollary~\ref{cor1}(a).

\textbf{Case 2.2.} There exists a $\{u_1,u_4\}$-vertex. Consider vertices $x,y \in U_1 \cap U_4$ such that $ext(vu_1xu_4) \cap U_1 \cap U_4 = \{y\}$ (if there exists a unique $\{u_1,u_4\}$-vertex $y$, then $x = w$). It is easy to see that $U_1 \setminus U_4 \in int(vu_1wu_4)$ by Corollary~\ref{cor1}(a). Hence by Lemma~\ref{l15} there exists a vertex $z \in int(u_1xu_4y) \cup ext(vu_1yu_4)$ not adjacent to $N_1[v] \cup \{x,y\}$. Consider a family $\mathcal{C} \in MCP(G_z)$ and its clique $W \ni v$. It is easy to check that one of the sets $\{u_1,u_4\}$, $\{u_2,y\}$, $\{u_3,y\} \in IS(G_z)$ dominates $W$, a contradiction. 

\textbf{Case 3.} Every vertex from $N_2(v)$ has at most two neighbors in $N_1(v)$. Moreover, there exists a $\{u_2,u_3\}$-vertex $w$. Recall that by Corollary~\ref{cor1}(b) $G$ has neither $\{u_1,u_2\}$-vertices nor $\{u_3,u_4\}$-vertices.

\textbf{Case 3.1.} There exist two distinct $\{u_1,u_4\}$-vertices $x$ and $y$ nonadjacent to $w$. We may assume that $U_1 \cap U_4 \cap ext(vu_1xu_4) = \{y\}$. Corollary~\ref{cor1}(a) implies that $w$ dominates both sets $U_1 \setminus U_4$ and $U_4 \setminus U_1$, hence these sets belong to $int(u_1u_2wu_3u_4x)$. By Lemma~\ref{l16}, $y$ is not a cutvertex, then $ext(vu_1yu_4) = \emptyset$ and $int(u_1xu_4y) \neq \emptyset$. By Lemma~\ref{l15}, since $U_1 \cap int(u_1xu_4y) = \emptyset$, there exists a vertex $z \in int(u_1xu_4y)$ not adjacent to $N_1[v] \cup \{x,w\}$. Consider a family $\mathcal{C} \in MCP(G_z)$ and its clique $W \ni v$. It is easy to check that one of the sets $\{u_1,u_4\}$, $\{u_2,x\}$, $\{u_3,x\} \in IS(G_z)$ dominates $W$, a contradiction. 

\textbf{Case 3.2.} There exists a unique $\{u_1,u_4\}$-vertex $x$ nonadjacent to $w$. Remind that $G$ has no $\{u_1,u_2\}$-vertices. It is easy to check, using Lemmas~\ref{l20}--\ref{l26}, that $deg(u_1) \geq 5$ and $|U_1 \setminus \{x\}| \geq 2$. Corollary~\ref{cor1}(a) implies that $(U_1 \setminus \{x\}) \cup \{w\}$ is a clique, then $|U_1 \setminus \{x\}| = 2$. Let $U_1 = \{u'_1,u''_1,x\}$, assume by symmetry that $u'_1 = int(u_1u''_1wu_2)$ (see Fig.~7). Since $deg(u'_1) \geq 4$, by Lemma~\ref{l15} the set $int(u_1u_2wu''_1) \setminus \{u'_1\}$ contains a vertex $z$ not adjacent to $N_1[v] \cup \{w,x\}$. Again, consider a family $\mathcal{C} \in MCP(G_z)$ and its clique $W \ni v$. One of the sets $\{u_1,u_4\}$, $\{u_2,x\}$, $\{u_3,x\} \in IS(G_z)$ dominates $W$, a contradiction.

\begin{center}
\includegraphics[scale=0.45]{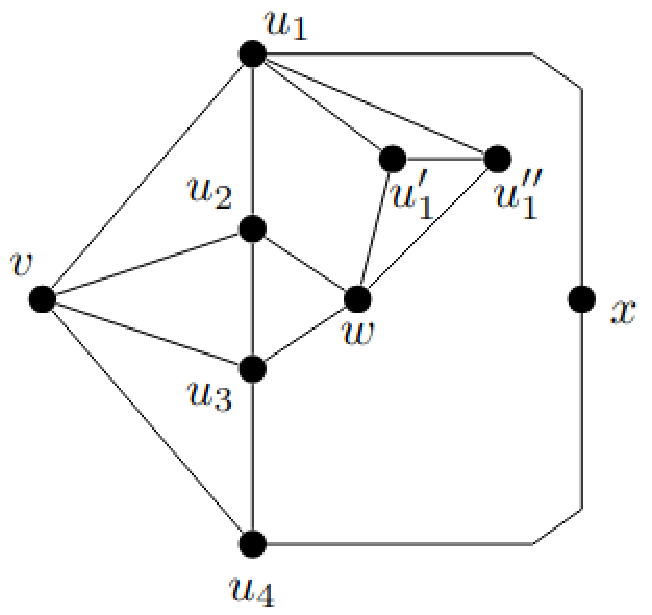}

Fig.~7. Illustration for the proof of Lemma~\ref{l27}, Case 3.2.
\end{center}

\textbf{Case 3.3.} There are no $\{u_1,u_4\}$-vertices nonadjacent to $w$. If $w$ dominates $U_1 \cup U_4$, then the vertices $u_1$ and $u_4$ are pendant in $G_w$, thus $\gamma(G_w) < \theta(G_w)$ by Lemma~\ref{l10}(a), a contradiction. Otherwise there exists a vertex $y \in U_1 \cup U_4$ such that $wy \notin E(G)$. Assume by symmetry that $y \in U_1$, then the set $\{w,y\}$ is $\theta$-independent in $G_{v,u_4}$ and dominates the set $\{u_1,u_2,u_3\}$, a contradiction by Corollary~\ref{cor1}(a).

\textbf{Case 4.} Every vertex from $N_2(v)$ has at most two neighbors in $N_1(v)$. Moreover, $G$ has no $\{u_2,u_3\}$-vertices. If $u_4 \in ext(vu_1u_2u_3)$, then we assume that for some $i \in \{1,2\}$ the graph $G$ does not have $\{u_i,u_{i+2}\}$-vertices (say, $i = 2$).

\textbf{Case 4.1.} There exist a pair of distinct $\{u_1,u_4\}$-vertices $x$ and $y$. We may assume that $ext(vu_1xu_4) \cap U_1 \cap U_4 = \{y\}$. 

\textbf{Case 4.1.1.} There exists a vertex $u'_1 \in U_1 \cap int(u_1xu_4y)$. Note that $u'_1 \notin U_4$. Consider a vertex $u'_2 \in U_2$. Since $u'_2 \in ext(u_1xu_4y)$, the set $\{u'_1,u'_2\} \in IS(G_{v,u_4})$ dominates $\{u_1,u_2\}$, a contradiction by Corollary~\ref{cor1}(a).

\textbf{Case 4.1.2.} $U_1 \cap int(u_1xu_4y) = \emptyset$. Since $deg(y) \geq 4$, by Lemma~\ref{l15} there exists a vertex $z \in int(u_1xu_4y) \cup ext(vu_1yu_4)$ not adjacent to $N_1[v] \cup \{x\}$. Consider a family $\mathcal{C} \in MCP(G_z)$ and its clique $W \ni v$. It is easy to check that one of the sets $\{u_1,u_4\}$, $\{u_2,x\}$, $\{u_3,x\} \in IS(G_z)$ dominates $W$, thus $\gamma(G_z) < \theta(G_z)$ by Lemma~\ref{l9}, a contradiction.

\textbf{Case 4.2.} There exists a unique $\{u_1,u_4\}$-vertex. Note that $G$ has neither $\{u_1,u_2\}$-vertices nor $\{u_2,u_3\}$-vertices. Then Lemmas~\ref{l20}--\ref{l26} imply that $deg(u_1)$, $deg(u_2) \geq 5$ and $|U_1 \setminus U_4|, |U_2 \setminus U_4| \geq 2$. By Lemma~\ref{l1} there exist nonadjacent vertices $u_1' \in U_1 \cap V(G_{u_4})$ and $u_2' \in U_2 \cap V(G_{u_4})$, a contradiction by Corollary~\ref{cor1}(a).

\textbf{Case 4.3.} There are no $\{u_1,u_4\}$-veritces. If $u_4 \in ext(vu_1u_2u_3)$, we apply the argument from~Case~4.2. Assume that $u_4 \in int(vu_2u_3)$. 

\textbf{Case 4.3.1.} There exists a vertex $u'_2 \in U_2 \setminus U_4$. Lemmas~\ref{l20}--\ref{l26} imply that $deg(u_1) \geq 5$, hence $|U_1| \geq 3$ and there exists a vertex $u'_1 \in U_1$ such that the set $\{u'_1,u'_2\}$ is $\theta$-independent in $G_{v,u_4}$, a contradiction by Corollary~\ref{cor1}(a). 

\textbf{Case 4.3.2.} $U_2 \subseteq U_4$ and there exists a vertex $u'_3 \in U_3 \setminus U_1$. Note that $deg(u_2) \geq 5$, then there exists a $\{u_2,u_4\}$-vertex $u'_2$ nonadjacent to $u'_3$. Hence the set $\{u'_2,u'_3\} \in V(G_{v,u_1})$ dominates $\{u_3,u_4\}$, a contradiction by Corollary~\ref{cor1}(a).

\textbf{Case 4.3.3.} $U_2 \subseteq U_4$ and $U_3 \subseteq U_1$. Since $|U_3| \geq 2$, there exist vertices $x,y \in U_1 \cap U_3$ such that $int(u_1u_2u_3y) \cap U_3 = \{x\}$. Note that $U_2 \in int(vu_2u_3)$, then by Lemma~\ref{l15} there exists a vertex $z \in int(u_1u_2u_3x) \cup int(u_1xu_3y)$ not adjacent to $N_1[v] \cup \{x,y\}$. Consider a family $\mathcal{C} \in MCP(G_z)$ and its clique $W \ni v$. Note that $|U_1|,|U_4| \geq 3$ by Lemmas~\ref{l20}--\ref{l26}, thus there exist vertices $u^*_1 \in U_1$ and  $u^*_4 \in U_4$ such that both sets $\{u^*_1,u_3\}$ and $\{u^*_4,u_2\}$ are $\theta$-independent in $\mathcal{C}$. It is easy to check that one of the sets $\{u^*_1,u_3\}$, $\{u^*_4,u_2\}$, $\{u_1,u_4\}$ dominates $W$, a contradiction by Lemma~\ref{l9}.
\end{proof}

\section{Critical graphs with $\delta = 5$}

Throughout this section we consider a critical graph $G$ and its vertex $v$ such that $deg(v) = \delta(G) = 5$. Let $N_1(v) = \{u_1,\ldots,u_5\}$ and $S \subseteq N_1(v)$. Call a vertex $x \in N_2(v)$ \textit{$S$-vertex}, if $N_1(x) \cap N_1(v) = S$. We use the notation $U_i =  N_1(u_i) \setminus N_1(v)$.

\subsection{Additional properties}

\begin{lemma}\label{l28}
Suppose that $u_1u_2 \in E(G)$ and there exist four pairwise distinct vertices $u^1_1,u^2_1 \in U_1$ and $u^1_2,u^2_2 \in U_2$, such that $u^1_1u^1_2 \notin E(G)$ and $u^2_1u^2_2 \notin E(G)$. Then there exist $i,j \in \{1,2\}$ and $k \in [3;5]$ such that the set $\{u^i_1,u^j_2,u_k\}$ is independent.
\end{lemma}

\begin{proof}

Assume by symmetry that $u_3 \in ext(vu_1u_2)$. Consider two cases:

\textbf{Case 1.} $u_4,u_5 \in ext(vu_1u_2)$. We may assume that the vertices $u_1,\ldots,u_5$ are in clockwise order around $v$, and the vertices $u^1_2,u^2_2$ are in clockwise order around $u_2$. 

\textbf{Case 1.1.} There exists $i \in \{4,5\}$ such that $u_iu^2_2 \in E(G)$. Then $u_3 \in int(vu_2u^2_2u_i)$ and $u^1_1,u^1_2 \in ext(vu_2u^2_2u_i)$. Thus the set $\{u^1_1,u^1_2,u_3\}$ is independent, as required.

\textbf{Case 1.2.} $u_4u^2_2,u_5u^2_2 \notin E(G)$. If both sets $\{u^2_1,u^2_2,u_4\}$ and $\{u^2_1,u^2_2,u_5\}$ are not independent, then $u_4u^2_1, u_5u^2_1 \in E(G)$ and $u^1_2 \in int(vu_1u_1^2u_4)$. If, moreover, the set $\{u^1_1,u^1_2,u_5\}$ is not independent, then $u^1_1u_5 \in E(G)$. Hence $u^1_1 \in ext(vu_1u^2_1u_5)$ and $u^1_1u_3, u^1_1u_4 \notin E(G)$. If both sets $\{u^1_1,u^1_2,u_i\}, i \in \{3,4\}$ are not independent, then $u^1_2u_3,u^1_2u_4 \in E(G)$. Therefore, $u^2_2 \in int(vu_2u^1_2u_3)$ and the set $\{u^1_1,u^2_2,u_4\}$ is independent, as required (see Fig.~8).

\begin{center}
\includegraphics[scale=0.45]{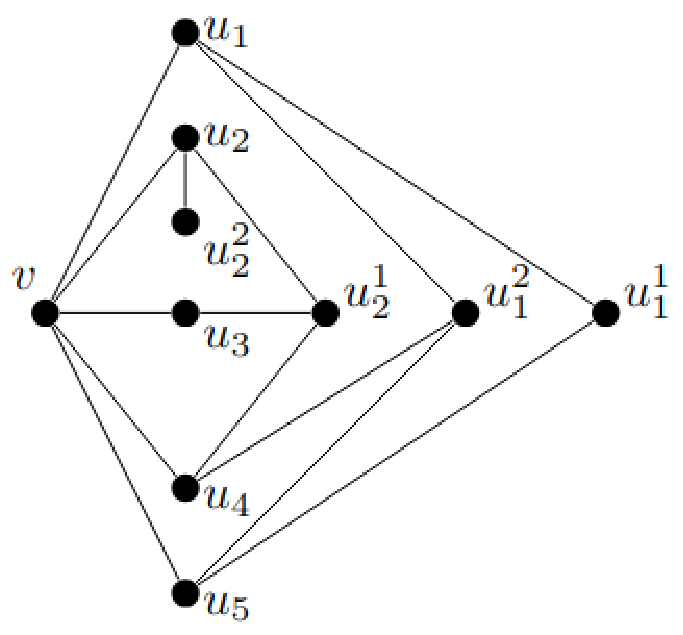}

Fig.~8. Illustration for the proof of Lemma~\ref{l28}, Case 1.2.
\end{center}

\textbf{Case 2.}  $\{u_4,u_5\} \cap int(vu_1u_2) \neq \emptyset$. Assume by symmetry that $u_3,u_4 \in ext(vu_1u_2)$ and $u_5 \in int(vu_1u_2)$. Suppose that for all $i \in [3;5]$ and $j \in \{1,2\}$ the set $\{u_1^j,u_2^j,u_i\}$ is not independent, then there exist $a,b \in \{1,2\}$ such that $u^a_1,u^b_2 \in ext(vu_1u_2)$ and $u^{3-a}_1,u^{3-b}_2 \in int(vu_1u_2)$. Both $u_1^a$ and $u_2^b$ dominate $\{u_3,u_4\}$ (if, say, $u_1^au_3 \notin E(G)$, then the set $\{u_1^a,u_2^{3-b},u_3\}$ is independent, a contradiction). This contradicts the planarity of $G$.
\end{proof}

\begin{lemma}\label{l29}
If the subgraph $G[N_1(v)]$ is triangle-free, than $N_4(v) = \emptyset$. 
\end{lemma}

\begin{proof}
Suppose there exists a vertex $x \in N_4(v)$. Note that $N_2[v] \subseteq V(G_x)$. Consider a family $\mathcal{C} \in MCP(G_x)$ and its clique $W \ni v$. Since $G[N_1(v)]$ is triangle-free and $W$ is maximal by inclusion, we have $|W| \in \{2,3\}$. Suppose that $|W| = 2$ (say, $W = \{v,u_1\}$). Consider a vertex $u'_1 \in U_1$ (it exists by Lemma~\ref{l12}). The set $\{u'_1,u_2,u_3,u_4\}$ is not a clique, then for some $i \in [2;4]$ the set $\{u'_1,u_i\}$ is $\theta$-independent in $\mathcal{C}$ and dominates $W$. Thus $\gamma(G_x) < \theta(G_x)$ by Lemma~\ref{l9}, a contradiction.

We now suppose that $|W| = 3$ (say, $W = \{v,u_1,u_2\}$). Consider three cases:

\textbf{Case 1.} Both vertices $u_1$ and $u_2$ are pendant in $G[N_1(v)]$. If there exists a vertex $w \in U_1 \cap U_2$, then the set $\{w,u_3,u_4,u_5\}$ is not a clique and for some $i \in [3;5]$ the set $\{w,u_i\}$ is $\theta$-independent in $\mathcal{C}$ and dominates~$W$, a contradiction by Lemma~\ref{l9}. Suppose that $U_1 \cap U_2 = \emptyset$. Since $|U_1|, |U_2| \geq 3$, by Lemma~\ref{l1} there exist vertices $u_1^{1},u_1^{2} \in U_1$ and $u_2^{1},u_2^{2} \in U_2$ such that $u_1^1u_2^1,u_1^2u_2^2 \notin E(G)$. By Lemma~\ref{l28}, there exist integers $i,j \in [1;2]$ and $k \in [2;5]$ such that the set $\{u_1^{i},u_2^{j},u_k\}$ is independent. This set dominates $\{v,u_1,u_2\}$ in $G_x$, again a contradiction by Lemma~\ref{l9}.

\textbf{Case 2.} Either $u_1$ or $u_2$ is not pendant in $G[N_1(v)]$. Assume that $u_2u_3 \in E(G)$. Since $|U_1| \geq 3$, the set $\{u_3\} \cup U_1$ is not a clique. Then there exists a vertex $u_1' \in U_1$ such that the set $\{u_1',u_3\}$ is $\theta$-independent in $\mathcal{C}$ and dominates $W$, a contradiction by Lemma~\ref{l9}.

\textbf{Case 3.} Both vertices $u_1$ and $u_2$ are not pendant in $G[N_1(v)]$. Since $G[N_1(v)]$ is triangle-free, we may assume  that $u_1u_5,u_2u_3 \in E(G)$ and $u_1u_3,u_2u_5 \notin E(G)$. Note that $U_1,U_2 \neq \emptyset$. If there exists a vertex $w \in U_1 \cap U_2$, then the set $\{w,u_3,u_4,u_5\}$ is not a clique and for some $i \in [3;5]$ the set $\{w,u_i\}$ is $\theta$-independent in $\mathcal{C}$ and covers $W$, a contradiction by Lemma~\ref{l9}. Otherwise consider some vertices $u'_1 \in U_1$ and $u'_2 \in U_2$. If none of the sets $\{u'_1,u_3\}$, $\{u'_2,u_5\}$, $\{u_3,u_5\}$ is $\theta$-independent in $\mathcal{C}$, then $u'_1u_3, u'_2u_5, u_3u_5 \in E(G)$ and the induced subgraph $G[\{v,u_1,u_2,u_3,u_5,u'_1,u'_2\}]$ is nonplanar, a contradiction.
\end{proof}

\begin{corollary}\label{cor2}
Suppose there exists a vertex $x \in V(G) \setminus N_1[v]$ and an integer $k \in [1;5]$ such that $x$ dominates $U_k$. Then $|U_k| \leq 2$.
\end{corollary}

\begin{proof}
Suppose for a contradiction that $|U_k| \geq 3$. Then there exist vertices $y_1,y_2,y_3 \in U_k$ such that $U_k \cap int(u_ky_1xy_3) = y_2$. Since $deg(y_2) \geq 5$, the set $int(u_ky_1xy_2) \cup int(u_ky_2xy_3)$ is nonepmty, assume by symmetry that $int(u_ky_1xy_2) \neq \emptyset$. If $u_kx \in E(G)$ and $y_1 \in ext(u_ky_2x)$, then one of the edges $y_1x$ and $y_2x$ is separating, a contradiction by Lemma~\ref{l17}. Otherwise if $y_1y_2 \in E(G)$ ($y_1y_2 \notin E(G)$), let $C = y_1y_2x$ ($C = u_ky_1xy_2$). By Lemma~\ref{l15}, there exists a vertex $z \in int(C)$ not adjacent to the vertices from $C$. It is easy to see that $z \notin N_3[v]$, hence $N_4(v) \neq \emptyset$, a contradiction by the previous lemma.
\end{proof}

\subsection{Case $\delta = 5$}

In this subsection we finish the proof of Theorem~\ref{thm3}. It is easy to check that precisely one of the following possibilities holds:

\begin{enumerate}
\item $G[N_1(v)]$ has a cycle $C_3$ (Lemma~\ref{l30})
\item $G[N_1(v)]$ has a cycle $C_4$ (Lemma~\ref{l31})
\item $G[N_1(v)]$ has a cycle $C_5$ (Lemma~\ref{l32})
\item $G[N_1(v)]$ is acyclic and $\Delta(G[N_1(v)]) \geq 3$ (Lemma~\ref{l33})
\item $G[N_1(v)] \cong 5K_1$ (Lemma~\ref{l34})
\item $G[N_1(v)] \cong K_2 \cup 3K_1$ (Lemma~\ref{l35})
\item $G[N_1(v)] \cong 2K_2 \cup K_1$ (Lemma~\ref{l36})
\item $G[N_1(v)] \cong P_3 \cup 2K_1$ (Lemma~\ref{l37})
\item $G[N_1(v)] \cong P_3 \cup K_2$ (Lemma~\ref{l38})
\item $G[N_1(v)] \cong P_5$ (Lemma~\ref{l39})
\item $G[N_1(v)] \cong P_4 \cup K_1$ (Lemma~\ref{l40})
\end{enumerate}

\begin{lemma}\label{l30}
If $G[N_1(v)]$ has a cycle $C_3$, then $G$ is not critical.
\end{lemma}

\begin{proof}

We may assume that $G[N_1(v)]$ has a cycle $u_1u_2u_3$ and $v,u_4,u_5 \in int(u_1u_2u_3)$. Assume by symmetry that $u_4 \in int(vu_1u_2)$.

\textbf{Case 1.} $u_5 \in int(vu_1u_3)$ or $u_5 \in int(vu_2u_3)$. Assume that $u_5 \in int(vu_1u_3)$.

\textbf{Case 1.1.} $u_4u_1,u_4u_2 \in E(G)$ or $u_5u_1,u_5u_3 \in E(G)$ (say, $u_4u_1,u_4u_2 \in E(G)$). Lemmas~\ref{l13} and~\ref{l14} imply that 
$$\gamma(G) - 2 = \gamma(G_{v,u_4}) = \theta(G_{v,u_4}) < \gamma(G_{u_4}) = \theta(G_{u_4}) = \gamma(G) - 1.$$ 

Note that $V(G_{u_4}) = V(G_{v,u_4}) \cup \{u_3,u_5\}$ and $U_4 \cap (U_3 \cup U_5) = \emptyset$. Since $|U_3|,|U_5| \geq 2$, Lemma~\ref{l1} implies that there exists a vertex $x \in V(G_{v,u_4})$ that dominates $\{u_3,u_5\}$ or there exist nonadjacent vertices $u'_3 \in U_3 \cap V(G_{v,u_4})$ and $u'_5 \in U_5 \cap V(G_{v,u_4})$. Then $\gamma(G_{u_4}) \leq \gamma(G_{v,u_4})$ by Lemma~\ref{l8}, a contradiction.

\textbf{Case 1.2.} $u_4u_2,u_5u_3 \in E(G)$. Case~1.1 implies that $u_4u_1,u_5u_1 \notin E(G)$, then the set $N_1[u_1] \cap N_1[v]$ is a clique and Lemmas~\ref{l13} and~\ref{l14} imply that

$$\gamma(G) - 2 = \gamma(G_{v,u_1}) = \theta(G_{v,u_1}) < \gamma(G_{u_1}) = \theta(G_{u_1}) = \gamma(G) - 1.$$ 

\textbf{Case 1.2.1.} Neither $U_4 \subseteq U_1$ nor $U_5 \subseteq U_1$. There exist vertices $u_4' \in U_4 \setminus U_1$ and $u_5' \in U_5 \setminus U_1$. Since $u'_4 \in int(vu_1u_2)$ and $u'_5 \in int(vu_1u_3)$, we have $u_4'u_5' \notin E(G)$. Note that $u'_4,u'_5 \in V(G_{v,u_1})$, thus  $\gamma(G_{u_1}) \leq \gamma(G_{v,u_1})$ by Lemma~\ref{l8}, a contradiction.

\textbf{Case 1.2.2.} There exists $i \in \{4,5\}$ such that $U_i \subseteq U_{1}$. Assume by symmetry that $i = 4$. Since $|U_4| \geq 3$, there exists a vertex $u'_4 \in U_1 \cap U_4$ nonadjacent to $U_2 \cup U_5$. Consider a family $\mathcal{C} \in MCP(G_{u'_4})$, its clique $W \ni v$ and some vertices $u'_2 \in U_2$, $u'_5 \in U_5$. If $W = \{v,u_2,u_3\}$, then $\{u'_2,u_5\}$ dominates $W$. Otherwise $W = \{v,u_3,u_5\}$ and $\{u_2,u'_5\}$ dominates $W$, a contradiction.

\textbf{Case 1.3.} There exists $i \in \{4,5\}$ such that $u_iu_{i-2} \notin E(G)$. Assume by symmetry that $i = 4$. Lemmas~\ref{l13} and~\ref{l14} imply that 
$$\gamma(G) - 2 = \gamma(G_{v,u_2}) = \theta(G_{v,u_2}) < \gamma(G_{u_2}) = \theta(G_{u_2}) = \gamma(G) - 1.$$ 

\textbf{Case 1.3.1.} There exists a vertex $u'_4 \in U_4 \setminus U_2$. Consider a vertex $u'_5 \in U_5$, then the set $\{u'_4,u'_5\} \in IS(G_{v,u_2})$  dominates $\{u_4,u_5\}$. Thus $\gamma(G_{v,u_2}) \geq \gamma(G_{u_2})$ by Lemma~\ref{l8}, a contradiction. 

\textbf{Case 1.3.2.} $U_4 \subseteq U_2$. Since $|U_4| \geq 3$, there exists a vertex $x \in U_4 \cap U_2$ not adjacent to $U_1$. Moreover, since $vu_2$ is not a separating edge, there exists a vertex $u'_1 \in U_1 \cap int(vu_1u_2)$. Consider a family $\mathcal{C} \in MCP(G_x)$ and its clique $W \ni v$. 

\textbf{Case 1.3.2.1.} $W = \{v,u_3,u_5\}$ or $W = \{v,u_5\}$. Since $|U_5| \geq 3$, there exists a vertex $u''_5 \in U_5$ such that the set $\{u_1,u''_5\}$ is $\theta$-independent in $\mathcal{C}$ and dominates $W$, a contradiction.

\textbf{Case 1.3.2.2. }$W = \{v,u_1,u_5\}$. Again, since $|U_5| \geq 3$, there exists a vertex $u'''_5 \in U_5$ such that the set $\{u_3,u'''_5\}$ is $\theta$-independent in $\mathcal{C}$ and dominates $W$, a contradiction. 

\textbf{Case 1.3.2.3.} $W = \{v,u_1,u_3\}$. If there exists a vertex $u'_3 \in U_3$ nonadjacent to $u_5$, then the set $\{u_1',u'_3,u_5\} \in IS(G_{x})$ dominates $W$, a contradiction. Otherwise $U_3 \subseteq U_5$. By Lemma~\ref{l17} the edge $vu_3$ is not separating, thus there exists a vertex $u''_1 \in U_1 \cap int(vu_1u_3)$. Consider a vertex $z \in U_3 \cap U_5$ nonadjacent to $u''_1$ (since $|U_3| \geq 2$, such a vertex exists). Since $|U_4| \geq 3$ and $U_4 \subseteq ext(vu_1u_3)$, we have $\gamma(G_{u''_1,z}) < \theta(G_{u''_1,z})$ by Lemma~\ref{l10}(d), a contradiction.

\textbf{Case 2.} $u_4,u_5 \in int(vu_1u_2)$. Lemmas~\ref{l13} and~\ref{l14} imply that  $$\gamma(G) - 2 = \gamma(G_{v,u_3}) = \theta(G_{v,u_3}) < \gamma(G_{u_3}) = \theta(G_{u_3}) = \gamma(G) - 1.$$ If the vertices $u_4$ and $u_5$ have a common neighbor in $N_2(v)$, then by Lemma~\ref{l8} $\gamma(G_{u_3}) \leq \gamma(G_{v,u_3})$, a contradiction. Suppose no such vertex exists. By Lemma~\ref{l12}, neither $u_1$ nor $u_2$ dominates $\{u_4,u_5\}$. Thus $|U_4| + |U_5| \geq 4$ and the set $U_4 \cup U_5$ is not a clique. Hence there exist vertices $u'_4 \in U_4 \cap V(G_{v,u_3})$ and $u'_5 \in U_5 \cap V(G_{v,u_3})$ such that the set $\{u'_4,u'_5\}$ is $\theta$-independent in $G_{v,u_3}$, thus $\gamma(G_{v,u_3}) \geq \gamma(G_{u_3})$, a contradiction.
\end{proof}

\begin{lemma}\label{l31}
If $G[N_1(v)]$ has an induced cycle $C_4$, then $G$ is not critical.
\end{lemma}

\begin{proof}
We may assume that $G[N_1(v)]$ has a cycle $u_1u_2u_3u_4$ and $u_5 \in int(vu_1u_2)$. By Lemma~\ref{l30} we may also assume that $u_2u_5 \notin E(G)$.

\textbf{Case 1.} $u_1u_5 \notin E(G)$. By Lemma~\ref{l17}, the edges $vu_1$ and $vu_2$ are not separating, thus there exist vertices $u'_1 \in U_1 \cap int(vu_1u_2) $ and $u'_2 \in U_2 \cap int(vu_1u_2)$ (these vertices may or may not coincide). The set $\{u'_1,u_3\} \in IS(G)$ dominates the set $\{v,u_1,u_2,u_4\}$, hence by Lemma~\ref{l8} the equality $\gamma(G \setminus \{v,u_1,u_2,u_4\}) = \theta(G \setminus \{v,u_1,u_2,u_4\}) = \gamma(G) - 1$ is not possible. Lemmas~\ref{l13} and~\ref{l14} imply that $\gamma(G_{v,u_1}) < \gamma(G_{u_1})$. Therefore, $U_i \subseteq U_1$ for some $i \in \{3,5\}$. Likewise, $\gamma(G_{v,u_2}) < \gamma(G_{u_2})$ and, therefore, $U_j \subseteq U_2$ for some $i \in \{4,5\}$. Since $|U_i| \geq 2$ for all $i \in [3;5]$, we may assume by symmetry that $U_5 \subseteq U_1$ and $U_4 \subseteq U_2$ (see Fig.~9). Then there exists a vertex $x \in U_2 \cap U_4$ not adjacent to $U_3 \setminus U_4$. If $xu_1 \in E(G)$ ($xu_1 \notin E(G)$), then $\gamma(G_x) < \theta(G_x)$ by Lemma~\ref{l10}(b) (Lemma~\ref{l10}(e)), a contradiction. 

\begin{center}
\includegraphics[scale=0.45]{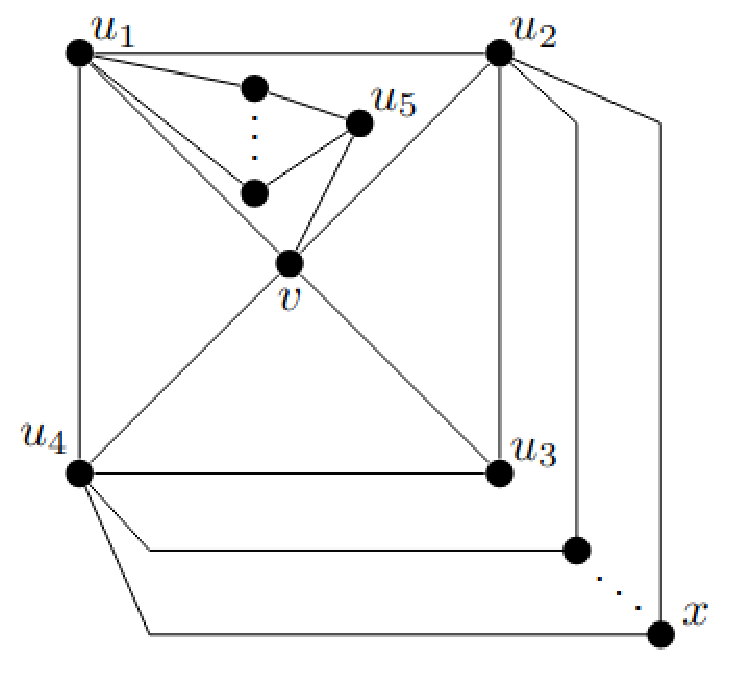}

Fig.~9. Illustration for the proof of Lemma~\ref{l30}, Case 1.
\end{center}

\textbf{Case 2.} $u_1u_5 \in E(G)$. By Lemma~\ref{l17}, there exists a vertex $u'_2 \in U_2 \cap int(vu_1u_2)$. The set $\{u'_2,u_4\} \in IS(G)$ dominates  $\{v,u_1,u_2,u_3\}$, thus $\gamma(G_{v,u_2}) < \gamma(G_{u_2})$ and $U_j \subseteq U_2$ for some $i \in \{4,5\}$.

\textbf{Case 2.1.} $U_4 \subseteq U_2$.  There exist vertices $u'_4,u''_4 \in U_4$ such that $U_4 \cap int(u_4u'_4u_2u_3) = \{u''_4\}$. Pick a vertex $u'_5 \in U_5 \setminus U_1$ (it exists by Lemma~\ref{l12} and may conicide with $u'_2$).  By Lemma~\ref{l17}, $int(vu_3u_4) = int(vu_2u_3) = \emptyset$, hence by Lemma~\ref{l15} there exists a vertex $z \in int(u_2u_3u_4u''_4)$ not adjacent to $N_1[v] \cup \{u'_2,u''_4,u'_5\}$. Consider a family $\mathcal{C} \in MCP(G_z)$ and its clique $W \ni v$.  If $W = \{v,u_1,u_2\}$ or $W = \{v,u_2,u_3\}$ then the set $\{u'_2,u_4\} \in IS(G_z)$ dominates it. If $W = \{v,u_1,u_5\}$ then the set $\{u_4,u'_5\} \in IS(G_z)$ dominates it. If $W = \{v,u_1,u_4\}$ then the set $\{u_3,u_5\} \in IS(G_z)$ dominates it. Finally, if $W = \{v,u_3,u_4\}$ then, since $u_1u''_4 \notin E(G)$, one of the sets $\{u_1,u_2\}$ and $\{u_2,u''_4\}$ is $\theta$-independent in $\mathcal{C}$ and dominates $W$, a contradiction.

\textbf{Case 2.2.} $U_5 \subseteq U_2$. Since $|U_5| \geq 3$, there exist vertices $y_1,y_2,y_3 \in U_5$ such that $int(u_5y_1u_2y_3) \cap U_5 = \{y_2\}$. We may assume that  $int(u_5y_1u_2y_2) \neq \emptyset$. If, moreover, $ int(u_5y_1u_2y_2) \cap U_2  = \emptyset$, then Lemma~\ref{l15} implies that $N_4(v) \neq \emptyset$, a contradiction by Lemma~\ref{l29}. Otherwise there exists a vertex $z \in U_2 \cap int(u_5y_1u_2y_2)$. Note that $z$ is not adjacent to $U_3 \cup U_4$. Consider a family $\mathcal{C} \in MCP(G_z)$ and its clique $W \ni v$. If $W = \{v,u_1,u_4\}$ then the set $\{u_3,u_5\} \in IS(G_z)$ dominates it. If $W = \{v,u_1,u_5\}$ then the set $\{y_3,u_4\} \in IS(G_z)$ dominates it. If $W = \{v,u_3,u_4\}$ then by Lemma~\ref{l1} there exists a vertex $w \in U_3 \cap U_4$ and the set $\{w,u_5\} \in IS(G_z)$ dominates $W$ or there exist nonadjacent vertces $u'_3 \in U_3, u'_4 \in U_4$ and the set $\{u_5,u'_3,u'_4\} \in IS(G_x)$ dominates $W$, a contradiction.

\end{proof}

\begin{lemma}\label{l32}
If $G[N_1(v)] \cong C_5$, then $G$ is not critical.
\end{lemma}

\begin{proof}
Assume that $G$ has a cycle $C = u_1u_2u_3u_4u_5$ and $v \in int(C)$. Lemmas~\ref{l30} and~\ref{l31} imply that $C$ is chordless. If $N_3(v) \neq \emptyset$, consider a vertex $x \in N_3[v]$ and a family $\mathcal{C} \in MCP(G_x)$ with a clique $W \ni v$. Since $W$ is maximal by inclusion, we have $|W| = 3$. Assume by symmetry that $W = \{v,u_1,u_2\}$, then the independent set $\{u_3,u_5\}$ dominates $W$, a contradiction. We now assume that $N_3(v) = \emptyset$. Consider four cases:

\textbf{Case 1.} There exists a vertex $w \in N_2(v)$, adjacent to three consecutive vertices of $C$ (assume by symmetry that $wu_1,wu_2,wu_3 \in E(G)$). Since $deg(u_2) \geq 5$, Lemma~\ref{l15} implies that there exists a vertex $x \in N_3(v) \cap int(vu_1wu_3)$, a contradiction.

\textbf{Case 2.} There exists a vertex $w \in N_2(v)$, adjacent to exactly three vertices of $C$ that are not consecutive. We may assume that $w$ is a $\{u_1,u_3,u_5\}$-vertex. Then $U_2 \cap U_4 = \emptyset$ and by Lemma~\ref{l10}(b) $w$ dominates either $U_2$ or $U_4$ (say, $U_2$). Corollary~\ref{cor2} implies that $|U_2| = 2$. Let $U_2 = \{u'_2,u''_2\}$, assume by symmetry that $u'_2 \in int(u_1u_2u''_2w)$. Since $deg(u'_2) \geq 5$, Lemma~\ref{l15} implies that there exists a vertex $x \in N_3(v) \cap int(u_1u_2u''_2w)$, a contradiction.

\textbf{Case 3.} There exists a vertex $w \in N_2(v)$, adjacent to a pair of nonadjacent vertices from $N_1(v)$. Assume by symmetry that $w$ is a $\{u_1,u_3\}$-vertex. 

\textbf{Case 3.1.}  $w$ dominates $U_2$. We use the argument from Case~2.

\textbf{Case 3.2.} There exists a $\{u_2\}$-vertex $u_2'$ nonadjacent to $w$. Consider a family $\mathcal{C} \in MCP(G_{u'_2})$ and its clique $W \ni v$.  It is easy to check that one of the independent sets $\{u_4,w\}$, $\{u_1,u_3\}$ and $\{u_5,w\}$ dominates $W$, a contradiction.

\textbf{Case 3.3.} For some $i \in \{1,3\}$ there exists a $\{u_2,u_i\}$-vertex $x$  nonadjacent to $w$. Assume by symmetry that $i = 3$. Consider a family $\mathcal{C} \in MCP(G_{x})$ and its clique $W \ni v$. If $W = \{v,u_1,u_5\}$ then the set $\{u_4,w\} \in IS(G_x)$ dominates it. Otherwise $W = \{v,u_4,u_5\}$. If there exists a vertex $u'_4 \in U_4 \setminus U_1$ then the set $\{u'_4,u_1\} \in IS(G_x)$ dominates $W$, a contradiction. Thus we assume that $U_4 \subseteq U_1$. By Lemma~\ref{l9}, it remains to consider the case when for every vertex $u'_4 \in U_4$ the set $\{u'_4,u_1\}$ is not $\theta$-independent in $G_x$, then $|U_4| = 2$. Let $U_4 = \{y_1,y_2\}$, assume that $y_2w \notin E(G)$ (see Fig.~10). If $y_2$ is adjacent to $U_5$, we rename the vertices and use the argument from Case~2. Otherwise  Lemma~\ref{l10}(b) implies that $\gamma(G_{\{y_2,w\}}) < \theta(G_{\{y_2,w\}})$, a contradiction.

\begin{center}
\includegraphics[scale=0.45]{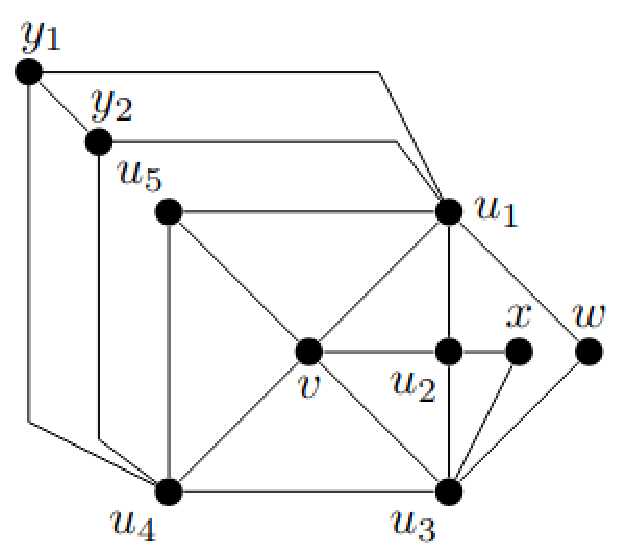}

Fig.~10. Illustration for the proof of Lemma~\ref{l32}, Case 3.3.
\end{center}

\textbf{Case 4.} Every vertex from $N_2(v)$ is adjacent to one vertex or two consecutive vertices from $N_1(v)$. By Lemma~\ref{l14}, it is sufficient to consider the following cases: 

\textbf{Case 4.1.} $\gamma(G \setminus \{u_1,u_2,u_3,v\}) = \theta(G \setminus \{u_1,u_2,u_3,v\}) = \gamma(G) - 1.$ If there exists a $\{u_1,u_2\}$-vertex  $x$, then the set $\{x,u_4\} \in IS(G)$ dominates $\{v,u_1,u_2,u_3\}$ in $G$, a contradiction. Otherwise by Lemma~\ref{l1} there exist nonadjacent vertices $u_1' \in U_1$ and $u_2' \in U_2$. By the previous cases, $u_1'u_4,u_2'u_4 \notin E(G)$, hence the set $\{u_1',u_2',u_4\} \in IS(G)$ dominates $\{v,u_1,u_2,u_3\}$ in $G$, again a contradiction.

\textbf{Case 4.2.} $ \gamma(G_{v,u_2}) = \theta(G_{v,u_2}) = \gamma(G) - 2.$ Since $|U_4|, |U_5| \geq 2$, by Lemma~\ref{l1} there exists a  $\{u_4,u_5\}$-vertex $x$ that dominates $\{u_4,u_5\}$ in $G_{u_2}$, or there exist vertices $u_4' \in U_4$ and $u_5' \in U_5$ such that the set is $\{u'_4,u'_5\} \in IS(G)$ dominates $\{u_4,u_5\}$ in $G_{u_2}$. Hence $\gamma(G_{v,u_2}) \geq \gamma(G_{u_2})$, a contradiction.
\end{proof}

\begin{lemma}\label{l33}
If $G[N_1(v)]$ is acyclic and $\Delta(G[N_1(v)]) \geq 3$, then $G$ is not critical.
\end{lemma}

\begin{proof}
We may assume that $u_1$ has the maximum degree in $(G[N_1(v)])$ and $u_1u_2$, $u_1u_3$, $u_1u_4 \in E(G)$. Lemma~\ref{l12} implies that $u_1u_5 \notin E(G)$. By Lemma~\ref{l30}, $G[N_1(v)]$ is triangle-free, thus $u_2u_3,u_3u_4,u_2u_4 \notin E(G)$. We may also assume that $u_2,u_5 \in int(vu_1u_3)$ and $u_4 \in ext(vu_1u_3)$. Lemma~\ref{l17} implies that $u_1u_4$ is not a separating edge, thus $ext(vu_1u_4) = \emptyset$ and $int(vu_3u_1u_4) \neq \emptyset$. By Lemma~\ref{l15}, there exists a vertex $x \in int(vu_3u_1u_4)$ not adjacent to $N_1[v] \cup U_2 \cup U_5$.  Consider a family $\mathcal{C} \in MCP(G_x)$ and a clique $W \ni v$. Remind that $W$ is maximal by inclusion, hence either $u_1 \in W$ or $u_5 \in W$.

\textbf{Case 1.} $W = \{v,u_5\}$. Consider a vertex $u_5' \in U_5$. The set $\{u_5',u_4\} \in IS(G_x)$ dominates $W$, a contradiction. 

\textbf{Case 2.}  $W = \{v,u_i,u_5\}$, where $i \in \{2,3\}$. Let $j = 5 - i$.  If $u_ju_5 \notin E(G)$ then $|U_5| \geq 3$ and there exists a vertex $u'_5 \in U_5$ such that the set $\{u_1,u'_5\}$ is $\theta$-independent in $\mathcal{C}$ and dominates $W$, a contradiction. Suppose that $u_ju_5 \in E(G)$. If the vertices $u_1$ and $u_j$ are $\theta$-independent in $\mathcal{C}$, we obtain a contradiction by Lemma~\ref{l9}. Otherwise, since $|U_5| \geq 2$ and the set $\{u_1,u_j\} \cup U_5$ is not a clique, there exists a vertex $u'_5 \in U_5$ such that the set $\{u'_5,u_1\}$ is $\theta$-independent in $\mathcal{C}$ and dominates $W$, again a contradiction. 

\textbf{Case 3.} $W = \{v,u_1,u_2\}$. Consider a vertex $u_2' \in U_2$. The set $\{u_2',u_4\} \in IS(G_x)$ dominates $W$, a contradiction. 

\textbf{Case 4.} $W = \{v,u_1,u_3\}$. By Corollary~\ref{cor2}, there exists a vertex $u'_3 \in U_3$ such that $u'_3x \notin E(G)$. Then one of the sets $\{u'_3,u_2\}$ and $\{u'_3,u_4\}$ is independent in $G_x$ and dominates $W$, a contradiction.

\textbf{Case 5.} $W = \{v,u_1,u_4\}$. By Corollary~\ref{cor2}, there exists a vertex $u_4' \in U_4$ such that $u_4'x \notin E(G)$. Then the set $\{u_2,u_4'\} \in IS(G_x)$ dominates $W$, a contradiction.  
\end{proof}

\begin{lemma}\label{l34}
If $N_1(v) \cong 5K_1$, then $G$ is not critical.
\end{lemma}

\begin{proof}

Consider two cases:

\textbf{Case 1.} $N_3(v) \neq \emptyset$. Select a vertex $x \in N_3(v)$. Consider a family $\mathcal{C} \in MCP(G_x)$ and its clique $W \ni v$. Assume by symmetry that $W = \{x,u_1\}$. By Lemma~\ref{l10}(e), $u_1$ is pendant in $G_x$. Thus $x$ dominates $U_1$, a contradiction by Corollary~\ref{cor2}. 

\textbf{Case 2.} $N_3(v) = \emptyset$. If every vertex from $N_2(v)$ has at least three neighbors in $N_1(v)$, then $\{v,u_1,u_2,u_3\} \in DS(G)$, thus $\gamma(G) < |N_1(v)| \leq \alpha(G)\leq \gamma^{\infty}(G)$, a contradiction. Otherwise there exists a vertex $x \in N_2(v)$, with at most two neighbors in $N_1(v)$. By Lemma~\ref{l10}(e), exactly one vertex from $N_1[v]$ is pendant in $G_x$. Thus for some $i \in [1;5]$ the vertex $x$ dominates the set $U_i$, again a contradiction by Corollary~\ref{cor2}.
\end{proof}

\begin{lemma}\label{l35}
If $G[N_1(v)] \cong K_2 \cup 3K_1$, then $G$ is not critical.
\end{lemma}

\begin{proof}
We may assume that $u_1u_2 \in E(G)$ and $u_3 \in ext(vu_1u_2)$. Consider two cases:

\textbf{Case 1.} There exists a vertex $w \in N_2(v)$ that dominates $N_1(v)$. Then $w,u_4,u_5 \in ext(vu_1u_2)$, we may assume that the vertices $u_1,u_2,\ldots,u_5$ are in clockwise order around $v$. Since $deg(u_4) \geq 5$, by Lemma~\ref{l15} the set $int(vu_3wu_5)$ contains a vertex $x$ not adjacent to $N_1[v] \cup \{w\}$. Consider a family $\mathcal{C} \in MCP(G_x)$ with a clique $W \ni v$. It is easy to check that for some $i \in [1;5]$ the set $\{u_i,w\}$ is $\theta$-independent in $\mathcal{C}$ and dominates $W$, a contradiction.

\textbf{Case 2.} No vertex from $N_2(v)$ dominates $N_1(v)$. Corollary~\ref{cor1}(b) implies that every vertex from $N_2(v)$ has at most three neighbors in $N_1(v)$ in $G$. Moreover, there are no $\{u_3,u_4,u_5\}$-vertices.

\textbf{Case 2.1.} There exists a vertex $w \in U_1 \cup U_2$ with at most two neighbors in $N_1(v)$. By Lemma~\ref{l10}(e), the equality $\gamma(G_w) = \theta(G_w)$ implies that for some $i \in [1;5]$ the vertex $w$ dominates the set $U_i$, a contradiction by Corollary~\ref{cor2}.

\textbf{Case 2.2.} Every vertex from $U_1 \cup U_2$ has exactly three neighbors in $N_1(v)$. Moreover, $u_4,u_5 \in ext(vu_1u_2)$. Since $|U_1|,|U_2| \geq 3$ and every vertex from $U_1 \cap U_2$ has a neighbor in $\{u_3,u_4,u_5\}$, it is easy to see that $|U_1 \cap U_2| \leq 1$ and thus $|U_1 \cup U_2| \geq 5$. Hence there exist two vertices from $U_1 \cup U_2$ with two common neighbors in $\{u_3,u_4,u_5\}$. Clearly, the subgraph $G[N_1[v] \cup U_1 \cup U_2]$ is nonplanar, a contradiction.

\textbf{Case 2.3.} Every vertex from $U_1 \cup U_2$ has exactly three  neighbors in $N_1(v)$. Moreover, $\{u_4,u_5\} \cap int(vu_1u_2) \neq \emptyset$. Assume by symmetry that $u_4 \in ext(vu_1u_2)$ and $u_5 \in int(vu_1u_2)$. If there exists a $\{u_1,u_2,u_5\}$-vertex $x$, then $\gamma(G_x) < \theta(G_x)$ by Lemma~\ref{l10}(d), a contradiction. Otherwise $U_1 \cup U_2 \in ext(vu_1u_2)$, thus $|U_1 \cap U_2| \leq 1$ and $|U_1 \cup U_2| \geq 5$. Since no vertex from $U_1 \cup U_2$ is adjacent to $u_5$, there exist two vertices from $U_1 \cup U_2$ with three common neighbors in $\{u_1,u_2,u_3,u_4\}$. Again, the subgraph $G[N_1[v] \cup U_1 \cup U_2]$ is nonplanar, a contradiction.
\end{proof}

\begin{lemma}\label{l36}
If $G[N_1(v)] \cong P_3 \cup 2K_1$, then $G$ is not critical.
\end{lemma}

\begin{proof}

Assume by symmetry that $E(N_1(v)) = \{u_1u_2,u_2u_3\}$ and $u_3 \in ext(vu_1u_2)$. Consider some vertices $x_1,x_2 \in U_2$. It follows from Lemma~\ref{l10}, Corollary~\ref{cor1}(b) and Corollary~\ref{cor2} that both of them have either three or five neighbors in $N_1(v)$. Consider two cases: 

\textbf{Case 1.} At least one of the vertices $u_4,u_5$ belongs to $int(vu_1u_2u_3)$. Assume by symmetry that $u_4 \in int(vu_1u_2)$.

\textbf{Case 1.1.} $u_5 \in int(vu_1u_2)$. If for some $i \in \{1,2\}$ the vertex $x_i$ is a $\{u_1,u_2,u_3\}$-vertex, then $\gamma(G_{x_i}) < \theta(G_{x_i})$ by Lemma~\ref{l10}(d), a contradiction. Otherwise $x_1,x_2 \in int(vu_1u_2)$. There exist $i,j \in [1;2]$ such that $x_iu_{j+3} \notin E(G)$. Assume by symmetry that $i = j = 1$. By Corollary~\ref{cor2}, $u_4$ is not pendant in $G_{x_1}$, thus $\gamma(G_{x_1}) < \theta(G_{x_1})$ by Lemma~\ref{l10}(b), a contradiction.

\textbf{Case 1.2.} $u_5 \in int(vu_2u_3)$. Then $u_4$ and $u_5$ have no common neighbors in $N_2(v)$, hence $G$ has no $\{u_1,u_2,u_3\}$-vertices by Lemma~\ref{l10}(b). Then either $x_1$ or $x_2$ is a $\{u_1,u_2,u_4\}$-vertex (say, $x_1$). Lemma~\ref{l10}(d) implies that $\gamma(G_{x_1}) < \theta(G_{x_1})$, a contradiction.

\textbf{Case 1.3.}  $u_4 \in int(vu_1u_2)$, $u_5 \in ext(vu_1u_2u_3)$. Assume that the vertices $u_1,u_2,u_3$,  $u_5$ are located clockwise around $v$. It is not hard to check, using Lemma~\ref{l10}(b,d) and Corollary~\ref{cor2}, that $x_1u_5, x_2u_5 \in E(G)$. We may assume that $x_1$ is a $\{u_1,u_2,u_5\}$-vertex and $x_2$ is a $\{u_2,u_3,u_5\}$-vertex. Then $x_1$ is not adjacent to $U_3 \setminus \{x_2\}$, hence $\gamma(G_{x_1}) < \theta(G_{x_1})$ by Lemma~\ref{l10}(b), a contradiction.

\textbf{Case 2.} $u_4,u_5 \in ext(vu_1u_2u_3)$. Then $x_1,x_2 \in ext(vu_1u_2u_3)$ and $G$ has no $\{u_1,u_2,u_3\}$-vertices. We may assume that the vertices $u_1,u_2,\ldots,u_5$ are in clockwise order around $v$ and the vertices $x_1$ and $x_2$ are in clockwise order around $u_2$. If $x_1u_3 \in E(G)$ ($x_2u_1 \in E(G)$) then $x_2$ ($x_1$) has at most two neighbors in $N_1(v)$, a contradiction. Thus both $x_1$ and $x_2$ have exactly three neighbors in $N_1(v)$.

\textbf{Case 2.1.} There exists $i \in \{1,2\}$ such that $x_i$ is a $\{u_2,u_4,u_5\}$-vertex. Assume by symmetry that $i = 1$. Then $u_1 \in ext(vu_2x_1u_4)$ and $u_3 \in int(vu_2x_1u_4)$.  Lemma~\ref{l10}(b) and Corollary~\ref{cor2} imply that $\gamma(G_{x_1}) < \theta(G_{x_1})$, a contradiction. 

\textbf{Case 2.2.} Neither $x_1$ nor $x_2$ is a $\{u_2,u_4,u_5\}$-vertex. Then $x_1u_1, x_2u_3 \in E(G)$. If $x_2u_5 \in E(G)$, then Lemma~\ref{l10}(b) and Corollary~\ref{cor2} imply that $\gamma(G_{x_2}) < \theta(G_{x_2})$, a contradiction. Hence $x_2u_4 \in E(G)$. If $x_1u_4 \in E(G)$ then Lemma~\ref{l10}(b) and Corollary~\ref{cor2} imply that $\gamma(G_{x_1}) < \theta(G_{x_1})$, hence $x_1u_5 \in E(G)$. If there exists a vertex $u'_1 \in U_1 \setminus U_5$, then $\{u'_1,x_2\} \in IS(G_{v,u_5})$ dominates $\{u_1,u_2,u_3,u_4\}$, a contradiction by Corollary~\ref{cor1}(a). Therefore, $U_1 \subseteq U_5$. By Case~2.1 $G$ has no $\{x_2,x_4,x_5\}$-vertices, thus $deg(u_2) = 5$. Lemma~\ref{l32} implies that $x_1x_2 \notin E(G)$. Since $x_2 \in int(vu_2x_1u_5)$, we have $N_1[x_2] \cap U_1 = \emptyset$. Since $|U_1| \geq 3$, we have $\gamma(G_{x_2}) < \theta(G_{x_2})$ by Lemma~\ref{l10}(d), a contradiction.
\end{proof}

\begin{lemma}\label{l37}
If $N_1(v) \cong 2K_2 \cup K_1$, then $G$ is not critical.
\end{lemma}

\begin{proof}
Assume by symmetry that $E(N_1(v)) = \{u_1u_2,u_3u_4\}$ and $u_3,u_4 \in ext(vu_1u_2)$. Consider three cases:

\textbf{Case 1.} There exists a vertex $w \in N_2(v)$ with at least four neighbors in $N_1(v)$. Corollary~\ref{cor1}(b) implies that $w$ dominates $N_1(v)$. Then $u_5 \in ext(vu_1u_2) \cap ext(vu_3u_4)$ and we may assume that the vertices $u_1,\ldots,u_5$ are in clockwise order around $v$. Since $deg(u_2) \geq 5$, the set $int(vu_1wu_3)$ contains a vertex is not adjacent to $N_1[v] \cup \{w\}$. Consider a family $\mathcal{C} \in MCP(G_x)$ and its clique $W \ni v$. It is easy to check that for some $i \in [1;5]$ the set $\{u_i,w\}$ is $\theta$-independent in $\mathcal{C}$ and dominates $W$, a contradiction.

\textbf{Case 2.} There exists a vertex $w \in N_2(v)$ with exactly three neighbors in $N_1(v)$. Corollary~\ref{cor1}(b) implies that $w$ is neither a $\{u_1,u_2,u_5\}$-vertex, nor a $\{u_3,u_4,u_5\}$-vertex.

\textbf{Case 2.1.} $wu_5 \notin E(G)$. Assume by symmetry that $w$ is a $\{u_1,u_2,u_3\}$-vertex. It is easy to check, using Corollary~\ref{cor1}(a), that $U_4 \setminus N_1[w] = U_5 \setminus N_1[w]$. Since $\gamma(G_w) = \theta(G_w)$, Lemma~\ref{l10}(a,d) implies that there exists a vertex $x \in U_4 \cap U_5 \cap N_1[w]$ or vertices $u'_4 \in U_4 \cap N_1[w]$ and $u'_5 \in U_5 \cap N_1[w]$. Moreover, Lemma~\ref{l10}(b) implies that there exists a vertex $y \in (U_4 \cap U_5) \setminus N_1[w]$. Consider a vertex $u'_3 \in U_3$. It is easy to check that $u'_3y \notin E(G)$. The set $\{u'_3,y\}$ dominates $\{u_3,u_4,u_5\}$ in one of the subgraphs $G_{u_1}$ and $G_{u_2}$, a contradiction by Corollary~\ref{cor1}(a).

\textbf{Case 2.2.} $wu_5 \in E(G)$. By Corollary~\ref{cor1}(b), $w$ is a $\{u_i,u_j,u_5\}$-vertex, for some $i \in [1;2]$ and $j \in [3;4]$. Then $u_5 \in ext(vu_1u_2) \cap ext(vu_3u_4)$ and we may assume that the vertices $u_1,\ldots,u_5$ are in clockwise order around $v$. If $(i,j) \neq (1,4)$, then Lemma~\ref{l10}(b) and Corollary~\ref{cor2} imply that $\gamma(G_w) < \theta(G_w)$. Suppose that $(i,j) = (1,4)$. By Lemma~\ref{l10}(b) and Corollary~\ref{cor2} there exists a vertex $x \in U_2 \cap U_3$ nonadjacent to $w$. Case~2.1 implies that $x$ is a $\{u_2,u_3\}$-vertex, then $\gamma(G_x) < \theta(G_x)$ by Lemma~\ref{l10}(e), a contradiction.

\textbf{Case 3.} Every vertex of $N_2(v)$ has at most two neighbors in $N_1(v)$.

\textbf{Case 3.1.} For some $i \in \{1,3\}$ there exists $\{u_i,u_{i+1}\}$-vertex $w$. Assume by symmetry that $i = 1$. By Corollary~\ref{cor2}, there exists a vertex $u'_5 \in U_5$ such that $wu'_5 \notin E(G)$. By Corollary~\ref{cor1}(b), for some $j \in \{3,4\}$ we have $u_ju'_5 \notin E(G)$. Hence the set $\{w,u'_5\} \in IS(G_{v,u_j})$ dominates $\{u_1,u_2,u_5\}$, a contradiction by Corollary~\ref{cor1}(a). 

\textbf{Case 3.2.} For some $i,j \in \{1,2\}$ there exists $\{u_i,u_{j+2}\}$-vertex $w$. Then $\gamma(G_w) < \theta(G_w)$ by Lemma~\ref{l10}(e) and Corollary~\ref{cor2}, a contradiction.

\textbf{Case 3.3.} For some $i \in [1;4]$ there exists $\{u_i,u_5\}$-vertex $w$ (assume that $i = 4$). Corollary~\ref{cor2} implies that there exists a vertex $u'_3 \in U_3$ nonadjacent to $w$. Case~3.2 implies that $u_1u'_3 \notin E(G)$. Then $\{u'_3,w\} \in IS(G_{v,u_1})$ dominates $\{u_3,u_4,u_5\}$, a contradiction by Corollary~\ref{cor1}(a).

\textbf{Case 3.4.} Every vertex from $N_2(v)$ has a unique neighbor in $N_1(v)$. By Lemma~\ref{l1}, for $i \in [1;4]$ there exist vertices $u'_i \in U_i$ such that $u'_1u'_2,u'_3u'_4 \notin E(G)$. Suppose there exists a vertex $u'_5 \in U_5$ such that one of the sets $\{u'_1,u'_2,u'_5\}$ and $\{u'_3,u'_4,u'_5\}$ is independent (say, $\{u'_1,u'_2,u'_5\}$). This set dominates $\{u_1,u_2,u_5\}$ in a subgraph $G_{v,u_3}$, a contradiction by Corollary~\ref{cor1}(a). Otherwise every vertex from $U_5$ has a neighbor in $\{u'_1,u'_2\}$ and a neighbor in $\{u'_3,u'_4\}$. Since $|U_5| \geq 4$, the subgraph $G[N_2[v]]$ is nonplanar, a contradiction.
\end{proof}

\begin{lemma}\label{l38}
If $G[N_1(v)] \cong P_3 \cup K_2$, then $G$ is not critical.
\end{lemma}

\begin{proof}

Assume by symmetry that $E(N_1(v)) = \{u_1u_2,u_2u_3,u_4u_5\}$ and $u_3 \in ext(vu_1u_2)$. By Corollary~\ref{cor1}(b), if a vertex from $N_2(v)$ is adjacent to exactly four vertices from $N_1(v)$, then it is a $\{u_1,u_3,u_4,u_5\}$-vertex. Moreover, $G$ has neither $\{u_1,u_4,u_5\}$-vertices nor $\{u_3,u_4,u_5\}$-vertices. Consider three cases:

\textbf{Case 1.} There exists a vertex $w \in N_2(v)$ with at least four neighbors in $N_1(v)$. Then $u_4,u_5 \in ext(vu_1u_2u_3)$, hence we may assume that the vertices $u_1,\ldots, u_5$ are in clockwise order around $v$. Since $U_4 \setminus \{w\} \neq \emptyset$, by Lemma~\ref{l15} there exists a vertex $x \in int(vu_3wu_5)$ not adjacent to the set $N_1[v] \cup U_2 \cup \{w\}$. Consider a family $\mathcal{C} \in MCP(G_x)$ and its clique $W \ni v$. Clearly, $|W| = 3$. If $u_2w \in E(G)$, then it is easy to check that for some $i \in [1;5]$ the set $\{w,u_i\}$ is $\theta$-independent and dominates $W$, a contradiction. 

Suppose that $u_2w \notin E(G)$. If $W = \{v,u_4,u_5\}$, then one of the sets $\{u_1,w\}$ and $\{u_3,w\}$ is $\theta$-independent in $\mathcal{C}$ and dominates $W$, a contradiction. Otherwise assume by symmetry that $W = \{v,u_1,u_2\}$. If the set $\{u_3,w\}$ is $\theta$-independent in $\mathcal{C}$, we obtain a contradiction by Lemma~\ref{l9}. Otherwise there exists a clique $W' \in \mathcal{C}$ such that $u_3,w \in W'$. Clearly, $u_5 \notin W'$. Since $|U_2| \geq 2$, there exists a vertex $u'_2 \in U_2 \setminus W'$. Note that $u'_2 \in int(vu_1wu_3)$, hence $u'_2u_5 \notin E(G)$. Then the set $\{u'_2,w,u_5\}$ is $\theta$-independent in $\mathcal{C}$ and dominates $W$, a contradiction.

\textbf{Case 2.} Every vertex from $N_2(v)$ has at most three neighbors in $N_1(v)$. Moreover, there exists a vertex $w \in U_2$ that is not a $\{u_2\}$-vertex.

\textbf{Case 2.1.} $w$ is a $\{u_2,u_i\}$-vertex for some $i \in \{4,5\}$. Lemma~\ref{l10}(e) and Corollary~\ref{cor2} imply that $\gamma(G_w) < \theta(G_w)$, a contradiction.

\textbf{Case 2.2.}  $w$ is a $\{u_2,u_4,u_5\}$-vertex. If $w \in ext(vu_1u_2u_3)$, then $\gamma(G_w) < \theta(G_w)$ by Lemma~\ref{l10}(b) and Corollary~\ref{cor2}, a contradiction. Hence $w \in int(vu_1u_2u_3)$, we may assume that $w \in int(vu_2u_3)$. If there exists a vertex $u'_1 \in U_1 \setminus U_3$, then $\{u'_1,w\}$ dominates $\{u_1,u_2,u_4,u_5\}$ in $G_{u_3}$, a contradiction by Corollary~\ref{cor1}(a). Otherwise $U_1 \subseteq U_3$. Since $w$ is not adjacent to $U_1$, we have $\gamma(G_w) < \theta(G_w)$ by Lemma~\ref{l10}(d), a contradiction.

\textbf{Case 2.3.} $w$ is a $\{u_i,u_2,u_j\}$-vertex, where $i \in \{1,3\}$ and $j \in \{4,5\}$. 

\textbf{Case 2.3.1.} $w \in int(vu_1u_2u_3)$. Assume that $w \in int(vu_2u_3)$, then $u_4,u_5 \in int(vu_2u_3)$ and $i = 3$. The vertices $u_{1}$ and $u_{9-j}$ have no common neighbors in $N_2(v)$, thus Lemma~\ref{l10}(b) and Corollary~\ref{cor2}  imply that $\gamma(G_w) < \theta(G_w)$, a contradiction.

\textbf{Case 2.3.2.} $w \in ext(vu_1u_2u_3)$. Then $u_4,u_5 \in ext(vu_1u_2u_3)$, we may assume that the vertices $u_1,\ldots,u_5$ are in clockwise order around $v$. Lemma~\ref{l10}(b) and Corollary~\ref{cor2} imply that either $(i,j) = (1,5)$ or $(i,j) = (3,4)$. Assume by symmetry that $(i,j) = (1,5)$. Lemma~\ref{l10}(b) implies that there exists a vertex $x \in U_3 \cap U_4$ nonadjacent to $w$.

\textbf{Case 2.3.2.1.} There exists a vertex $u'_5 \in U_5 \setminus U_1$. By Corollary~\ref{cor1}(a), $u'_5x \in E(G)$. If there exists a vertex $u'_4 \in U_4 \setminus U_3$, then $\{u'_4,w\}$ dominates $\{u_1,u_4,u_5\}$ in $G_{u_3}$, a contradiction by Corollary~\ref{cor1}(a). Otherwise $U_4 \subseteq U_3$ and there exists a vertex $x' \in (U_3 \cap U_4) \setminus \{x\}$. Therefore, $\{x',u'_5\} \in IS(G_{v,u_1})$ dominates $\{u_3,u_4,u_5\}$, again a contradiction by Corollary~\ref{cor1}(a).

\textbf{Case 2.3.2.2.} $U_5 \subseteq U_1$. Then $xu_2 \notin E(G)$ (if $xu_2 \in E(G)$, then $\gamma(G_{x}) < \theta(G_{x})$ by Lemma~\ref{l10}(d), a contradiction). Since $|U_5| \geq 3$, there exist vertices $y_1,y_2 \in U_5 \setminus \{w\}$ such that $U_5 \cap ext(vu_1y_1u_5) = \{y_2\}$. By Lemma~\ref{l17}, $ext(vu_1y_2u_5) = \emptyset$, thus by Lemma~\ref{l15} there exists a vertex $z \in int(u_1y_1u_5y_2)$ not adjacent to $N_1[v] \cup \{x,w,y_1\}$. Consider a family $\mathcal{C} \in MCP(G_z)$ and its clique $W \ni v$. It is easy to check that one of the sets $\{u_3,w\}$, $\{u_1,x\}$, $\{u_2,x,y_1\} \in IS(G_z)$ dominates $W$, a contradiction.

\textbf{Case 2.4.} $w$ is a $\{u_i,u_2\}$-vertex, $i \in \{1,3\}$. Assume that $i = 1$. By Corollary~\ref{cor2}, there exists a vertex $u'_3 \in U_3$ nonadjacent to $w$. Corollary~\ref{cor1}(b) implies that $u'_3u_j \notin E(G)$ for some $j \in \{4,5\}$. Then $\{w,u'_3\} \in IS(G_{v,u_j})$ dominates $\{u_1,u_2,u_3\}$, a contradiction by Corollary~\ref{cor1}(a).

\textbf{Case 3.} Every vertex from $N_2(v)$ has at most three neighbors in $N_1(v)$. Moreover, every vertex from $U_2$ is $\{u_2\}$-vertex.

\textbf{Case 3.1.} There exists a vertex $w \in U_1 \cap U_3$. Since $U_2 \cap (U_4 \cup U_5) = \emptyset$, the set $\{w\} \cup U_2$ is a clique by Corollary~\ref{cor1}(a). Let $U_2 = \{u'_2,u''_2\}$, assume by symmetry that $u'_2 \in int(u_1u_2u''_2w)$. By Lemma~\ref{l17}, $int(u_2u'_2u''_2) = \emptyset$. By Lemmas~\ref{l15} and~\ref{l29}, $int(u'_2u''_2w) = \emptyset$. Hence $int(u_1u_2u'_2w)$,  $int(u_2u_3wu''_2) \neq \emptyset$. By Lemma~\ref{l17}, there exist vertices $u'_1 \in U_1 \cap int(u_1u_2u'_2w)$ and $u'_3 \in U_3 \cap  int(u_2u_3wu''_2)$. Since the set $\{u'_2,u''_2\}$ is not $\theta$-independent in $G_{v,u_5}$ and $u'_1u'_3,u'_2u'_3,u''_2u'_1 \notin E(G)$, the set $\{u'_1,u'_2,u'_3\}$ is $\theta$-independent in $G_{v,u_5}$, a contradiction by Corollary~\ref{cor1}(a).

\textbf{Case 3.2.} There exists a vertex $w \in U_4 \cap U_5$. Note that $wu_1 \notin E(G)$ by Corollary~\ref{cor1}(b). Corollary~\ref{cor2} implies that there exists a vertex $u'_3 \in U_3$ nonadjacent to $w$. Then $\{u'_3,w\} \in IS(G_{v,u_1})$ dominates $\{u_3,u_4,u_5\}$, a contradiction.

\textbf{Case 3.3.} There exists $\{u_i,u_j\}$ vertex $x$ for some $i \in \{1,3\}$ and $j \in \{4,5\}$. 

\textbf{Case 3.3.1.} $u_4 \in int(vu_1u_2u_3)$. Assume by symmetry that $u_4 \in int(vu_2u_3)$ and $(i,j) = (3,4)$. By Corollary~\ref{cor2}, there exists a vertex $u'_5 \in U_5$ nonadjacent to $x$. Then the set $\{u'_5,x\} \in IS(G_{u_1})$ dominates $\{u_3,u_4,u_5\}$, a contradiction. 

In subcases~3.3.2 and~3.3.3 we may assume that $u_4,u_5 \in ext(vu_1u_2u_3)$ and the vertices $u_1,u_2,\ldots,u_5$ are in clockwise order around $v$.

\textbf{Case 3.3.2.} $(i,j) = (1,4)$ or $(i,j) = (3,5)$. Assume that $(i,j) = (3,5)$. By Corollary~\ref{cor2}, there exists a vertex $u'_4 \in U_4 \setminus U_1$ nonadjacent to $x$. Hence the set $\{u'_4,x\} \in IS(G_{v,u_1})$ dominates $\{u_3,u_4,u_5\}$, a contradiction.

\textbf{Case 3.3.3.} $(i,j) = (3,4)$ or $(i,j) = (1,5)$. Assume that $(i,j) = (3,4)$ and the number of $\{u_3,u_4\}$-vertices in $G$ is greater or equal than the number of $\{u_1,u_5\}$-vertices. Corollary~\ref{cor2} implies that there exists a vertex $y \in U_5$ nonadjacent to $x$. Corollary~\ref{cor1}(a) implies that $y$ is a $\{u_1,u_5\}$-vertex. 

\textbf{Case 3.3.3.1.} $x$ is the only $\{u_3,u_4\}$-vertex in $G$. Then there exist vertices $u'_4 \in U_4 \setminus U_3$ and $u'_5 \in U_5 \setminus U_1$. Corollary~\ref{cor1}(a) imply that $u'_4y,u'_5x \in E(G)$, this contradicts the planarity of $G$.

\textbf{Case 3.3.3.2.} There exists a $\{u_3,u_4\}$-vertex $x' \neq x$. Assume that $x' \in int(vu_3xu_4)$. By Lemma~\ref{l15}, there exists a vertex $z \in int(vu_3xu_4)$ nonadjacent to $N_1[v] \cup \{x,y\}$. Consider the family $\mathcal{C} \in MCP(G_z)$ and its clique $W \ni v$. It is easy to check that one of the sets $\{u_3,y\}$, $\{u_1,x\}$, $\{u_2,x,y\} \in IS(G_z)$ dominates $W$, a contradiction.

\textbf{Case 3.4.} All vertices in $N_2(v)$ have a unique neighbor in $N_1(v)$. By Lemma~\ref{l1}, there exist nonadjacent vertices $u'_4 \in U_4$ and $u'_5 \in U_5$. Corollary~\ref{cor1}(a) implies that every vertex from $U_1 \cup U_3$ is adjacent to $u'_4$ or $u'_5$. Then it is easy to see that $u_4,u_5 \in ext(vu_1u_2u_3)$, we may assume that the vertices $u_1,\ldots,u_5$ are in clockwise order around $v$. Corollary~\ref{cor2} implies that $u'_4$ does not dominate $U_3$, then there exists a vertex $u'_3 \in U_3$ such that $u'_3u'_5 \in E(G)$ and $u'_4 \in int(vu_3u'_3u'_5u_5)$. Hence $u'_4$ is not adjacent to $U_1$. By Corollary~\ref{cor2}, $u'_5$ does not dominate $U_1$, thus there exists a vertex $u'_1 \in U_1$ such that $u'_1u'_5 \notin E(G)$. Therefore, the set $\{u'_1,u'_4,u'_5\} \in IS(G_{v,u_3})$ dominates $\{u_1,u_4,u_5\}$, a contradiction by Corollary~\ref{cor1}(a).
\end{proof}

\begin{lemma}\label{l39}
If  $N_1(v) \cong P_5$, then $G$ is not critical.
\end{lemma}

\begin{proof}

We may assume that $E(N_1(v)) = \{u_1u_2,u_2u_3,u_3u_4,u_4u_5\}$ and $u_1 \in ext(vu_2u_3)$. Consider three cases:

\textbf{Case 1.} There exists a vertex $w \in N_2(v)$ with at least four neighbors in $N_1(v)$. Corollary~\ref{cor1}(b) implies that $u_1w,u_5w \in E(G)$. We may assume that $u_4w \in E(G)$ and the vertices $u_1,\ldots,u_5$ are located clockwise around $v$.

\textbf{Case 1.1.} $u_3w \in E(G)$. By Lemma~\ref{l15}, there exists a vertex $x \in int(vu_3wu_5)$ not adjacent to $N_1[v] \cup U_2 \cup \{w\}$. Consider a family $\mathcal{C} \in MCP(G_x)$ and its clique $W \ni v$. Suppose that none of the sets $\{u_i,w\}$, $i \in \{1,3,4,5\}$ is both $\theta$-independent in $\mathcal{C}$ and dominates $W$. Then it is easy to check that $u_2w \notin E(G)$ and $W = \{v,u_1,u_2\}$. Moreover, the set $\{u_3,w\}$ is not $\theta$-independent in $\mathcal{C}$. Note that $|U_2| \geq 2$ and the set $U_2 \cup \{u_3,w\}$ is not a clique. Then there exists a vertex $u_2' \in U_2$ such that the set $\{u_2',w\}$ is $\theta$-independent in $\mathcal{C}$. Since $u_3u_5 \notin E(G)$, the set $\{u_2',w,u_5\}$ is also $\theta$-independent in $\mathcal{C}$ and dominates $W$, a contradiction by Lemma~\ref{l9}.

\textbf{Case 1.2.} $u_3w \notin E(G)$. Note that $w$ dominates $\{u_1,u_5\}$, hence $\gamma(G_{v,u_3}) \geq \gamma(G_{u_3})$ and $\gamma(G \setminus \{v,u_2,u_3,u_4\}) < \gamma(G)$ by Lemma~\ref{l14}. Since $u_3 \in int(vu_2wu_4)$, neither $u_1$ nor $u_5$ is adjacent to $U_3$. Consider a vertex $u'_3 \in U_3$. The set $\{u_1,u'_3,u_5\} \in IS(G)$ dominates $\{v,u_2,u_3,u_4\}$, thus $\gamma(G \setminus \{v,u_2,u_3,u_4\}) \geq \gamma(G)$ by Lemma~\ref{l8}, a contradiction.

\textbf{Case 2.} There exists a vertex $w \in N_2(v)$ with exactly three neighbors in $N_1(v)$. Assume that either $wu_i,wu_{6-i} \in E(G)$ or $wu_{6-i} \notin E(G)$, where $i \in \{1,2\}$. If $u_1w,u_5w \in E(G)$, then we also assume that the vertices $u_1,u_2,\ldots,u_5$ are located clockwise around~$v$.

\textbf{Case 2.1.} $w$ is a $\{u_2,u_3,u_4\}$-vertex. By Corollary~\ref{cor2} and Lemma~\ref{l10}(b), there exists a vertex $x \in U_1 \cap U_5$ nonadjacent to $w$. Thus $u_4 \in ext(vu_1u_2u_3)$ and $u_5 \in ext(vu_1u_2u_3u_4)$.  By Lemma~\ref{l15} there exists a vertex $y \in int(vu_2wu_4)$ not adjacent to $N_1[v] \cup \{x\}$. Consider a family $\mathcal{C} \in MCP(G_y)$ and its clique $W \ni v$. It is easy to see that one of the sets $\{u_3,x\}, \{u_1,u_4\}, \{u_2,u_5\} \in IS(G_y)$ dominates $W$, a contradiction.

\textbf{Case 2.2.} $w$ is a $\{u_1,u_3,u_4\}$-vertex. Then $u_4 \in ext(vu_1u_2u_3)$ and $u_5$ is not adjacent to $U_2$. Lemma~\ref{l10}(b) implies that $w$ dominates $U_2$, then $|U_2| = 2$ by Corollary~\ref{cor2}. Let $U_2 = \{u'_2,u''_2\}$, assume that $u'_2 \in int(u_1u_2u''_2w)$. Corollary~\ref{cor1}(a) implies that $U_2 \cup \{w\}$ is a clique, thus $u'_2u''_2 \in E(G)$.  By Lemmas~\ref{l15} and~\ref{l29}, $int(u_2u'_2wu''_2) = \emptyset$, thus by Lemma~\ref{l17} there exist vertices $u'_1 \in U_1 \cap int(u_1u_2u'_2w)$ and $u'_3 \in U_3 \cap int(u_2u_3wu''_2)$. Lemmas~\ref{l30}--\ref{l38} imply that either $u_3u''_2 \in E(G)$ or $u_1u'_2 \in E(G)$. If $u_3u''_2 \in E(G)$ ($u_1u'_2 \in E(G)$), then the set $\{u'_1,u''_2\} \in IS(G_{v,u_5})$ ($\{u'_2,u'_3\} \in IS(G_{v,u_5})$) dominates $\{u_1,u_2,u_3\}$, a contradiction by Corollary~\ref{cor1}(a).

\textbf{Case 2.3.} $w$ is a $\{u_1,u_2,u_4\}$-vertex. Again, $u_4 \in ext(vu_1u_2u_3)$.

\textbf{Case 2.3.1.} $u_5 \in ext(vu_1u_2u_3u_4)$. Consider a vertex $u'_3 \in U_3$, it is easy to see that $u'_3u_1,u'_3u_5 \notin E(G)$. The set $\{u_1,u'_3,u_5\} \in IS(G)$ dominates $\{v,u_2,u_3,u_4\}$, hence $\gamma(G_{v,u_3}) < \gamma(G_{u_3})$ by Lemmas~\ref{l13} and \ref{l14}. Note that $|U_5| \geq 3$. By Corollary~\ref{cor2}, there exists a vertex $u'_5 \in U_5$ such that $\{w,u'_5\} \in IS(G_{v,u_3})$, thus $\gamma(G_{v,u_3}) \geq \gamma(G_{u_3})$ by Lemma~\ref{l8}, a contradiction.

\textbf{Case 2.3.2.} $u_5 \in int(vu_1u_2u_3u_4)$. Then $u_5 \in int(vu_3u_4)$. Lemma~\ref{l17} implies that there exists a vertex $x \in U_3 \cap int(vu_3u_4)$. Corollary~\ref{cor1}(a) implies that $xu_5 \in E(G)$, hence $xu_4 \notin E(G)$ by Corollary~\ref{cor1}(b). It is easy to check, using Corollary~\ref{cor1}(a), that $U_4 \cap int(u_1u_2u_3u_4w) = \emptyset$. By Lemma~\ref{l17}, $int(vu_1u_2) = int(vu_2u_3) = \emptyset$. Hence Lemma~\ref{l15} implies that there exists a vertex $z \in int(u_1u_2u_3u_4w)$ not adjacent to $N_1[v] \cup \{w,x\}$. Consider a family $\mathcal{C} \in MCP(G_z)$ and its clique $W \ni v$. Since $|U_5| \geq 4$, there exists a vertex $u'_5 \in U_5$ such that the set $\{u_3,u'_5\}$ is $\theta$-independent in $\mathcal{C}$. Therefore, one of the sets $\{w,u_5\}$, $\{u_1,u_4\}$, $\{u_2,u_5\}$, $\{u_3,u'_5\}$ is $\theta$-independent in $\mathcal{C}$ and dominates $W$, a contradiction.

\textbf{Case 2.4.} $w$ is a $\{u_1,u_2,u_5\}$-vertex. Clearly, $w$ dominates $\{u_1,u_5\}$ in $G_{u_3}$, thus $\gamma(G \setminus \{v,u_2,u_3,u_4\}) = \gamma(G) - 1$ by Lemmas~\ref{l8} and~\ref{l14}. Note that $U_3 \cap U_1 = \emptyset$. If there exists a vertex $u'_3 \in U_3 \setminus U_5$, then the set $\{u_1,u'_3,u_5\} \in IS(G)$ dominates $\{v,u_2,u_3,u_4\}$, a contradiction. Otherwise $U_3 \subseteq U_5$ and there exist vertices $x_1,x_2 \in U_3$ such that $x_2 \in ext(vu_3x_1u_5)$. Consider a vertex $u'_4 \in U_4$ (note that $u'_4 \neq x_1$ by Corollary~\ref{cor1}(b)). The set $\{u'_4,x_2\} \in IS(G_{v,u_1})$ dominates $\{u_3,u_4,u_5\}$, a contradiction by Corollary~\ref{cor1}(a).

\textbf{Case 2.5.} $w$ is a $\{u_1,u_3,u_5\}$-vertex. Since $G$ is $K_{3,3}$-free, such a vertex is unique. Thus all vertices from $N_2(v) \setminus \{w\}$ have at most two neighbors in $N_1(v)$.

\textbf{Case 2.5.1.} There exists a vertex $u'_3 \in U_3 \setminus (U_1 \cup U_5)$. The set $\{u_1,u'_3,u_5\} \in IS(G)$ dominates $\{v,u_2,u_3,u_4\}$, hence $\gamma(G_{v,u_3}) < \gamma(G_{u_3})$ by Lemma~\ref{l14}. Thus $G$ has no $\{u_1,u_5\}$-vertices. We may assume that $u'_3 \in int(vu_1wu_3)$. If there exists $\{u_4,u_5\}$-vertex $x$, the set $\{u'_3,x\} \in IS(G_{v,u_1})$ dominates $\{u_3,u_4,u_5\}$, a contradiction by Corollary~\ref{cor1}(a). Otherwise by Lemma~\ref{l1} there exist nonadjacent vertices $u'_4 \in U_4 \setminus U_1$ and $u'_5 \in U_5 \setminus U_1$. The set $\{u'_3,u'_4,u'_5\} \in IS(G_{v,u_1})$ dominates $\{u_3,u_4,u_5\}$, again a contradiction by Corollary~\ref{cor1}(a).

In subcases~2.5.2 and~2.5.3 we assume that $U_3 \subseteq U_1 \cup U_5$.

\textbf{Case 2.5.2.}  For some $i \in \{1,5\}$ there exist two $\{u_3,u_i\}$-vertices $x_1$ and $x_2$. Assume that $i = 1$ and $x_2 \in ext(vu_1x_1u_3)$. Consider a vertex $u'_2 \in U_2$. The set $\{u'_2,x_2\} \in IS(G_{v,u_5})$ dominates $\{u_1,u_2,u_3\}$, a contradiction. 

\textbf{Case 2.5.3.} $|(U_1 \cap U_3) \setminus \{w\}|$, $|(U_3 \cap U_5) \setminus \{w\}| \in \{0,1\}$. We may assume that there exists a $\{u_1,u_3\}$-vertex $x$. If $deg(u_3) = 5$,  then $xu_2 \in E(G)$ by Lemmas~\ref{l30}--\ref{l38}, a contradiction by Corollary~\ref{cor1}(b). If $deg(u_3) \geq 6$, then there exists a $\{u_3,u_5\}$-vertex $y$. Corollary~\ref{cor1}(a) implies that $x$ dominates $U_2$ and $y$ dominates $U_4$. Thus the vertices $u_2$ and $u_4$ are pendant in $G_{x,y}$, a contradiction by~Lemma~\ref{l10}(a).

\textbf{Case 3.} Every vertex from $N_2(v)$ has at most two neighbors in $N_1(v)$.

\textbf{Case 3.1.} $\gamma(G \setminus \{v,u_2,u_3,u_4\}) = \gamma(G) - 1$. If there exists a vertex $u'_3 \in U_3 \setminus (U_1 \cup U_5)$, then the set $\{u_1,u'_3,u_5\} \in IS(G)$ dominates $\{v,u_2,u_3,u_4\}$, a contradiction. Therefore, $U_3 \subseteq U_1 \cup U_5$. 

\textbf{Case 3.1.1.} $max(|U_1 \cap U_3|,|U_3 \cap U_5|) \leq 1$. Clearly, $deg(u_3) = 5$. There exist vertices $x_1 \in U_1 \cap U_3$ and $x_2 \in U_3 \cap U_5$. Lemmas~\ref{l30}--\ref{l38} imply that $u_2x_1 \in E(G)$ or $u_4x_2 \in E(G)$, a contradiction by Corollary~\ref{cor1}(b).

\textbf{Case 3.1.2.} For some $i \in \{1,5\}$ there exist two $\{u_3,u_i\}$-vertices $x_1$ and $x_2$. Assume that $i = 1$ and $x_2 \in ext(vu_1x_1u_3)$. If there exists a vertex $u'_2 \in U_2 \setminus U_5$, then the set $\{u'_2,x_2\} \in IS(G_{v,u_5})$ dominates $\{u_1,u_2,u_3\}$, a contradiction by Corollary~\ref{cor1}(a). Otherwise $U_2 \subseteq U_5$, thus $u_4 \in int(vu_2u_3)$ and $u_5 \in int(vu_2u_3u_4)$. Consider a vertex $u'_2 \in U_2$, note that $u'_2u_4, u'_2u_3 \notin E(G)$. Since $U_2 \cap int(u_1u_2u_3x_2) = \emptyset$, by Lemma~\ref{l15} there exists a vertex $z \in int(u_1u_2u_3x_2)$ not adjacent to $N_1[v] \cup \{x_1,x_2\}$. Consider a family $\mathcal{C} \in MCP(G_z)$ and its clique $W \ni v$. It is easy to see that one of the sets $\{x_1,u'_2,u_4\}$, $\{u_1,u_4\}$, $\{u_2,u_5\}$, $\{u_3,u'_2\} \in IS(G_z)$ dominates $W$, a contradiction. 

\textbf{Case 3.2.} $\gamma(G_{v,u_3}) < \gamma(G_{u_3})$. Clearly, $G$ has no $\{u_1,u_5\}$-vertices.

\textbf{Case 3.2.1.} There exist vertices $u'_1 \in U_1 \setminus U_3$ and $u'_5 \in U_5 \setminus U_3$. Then $u'_1u'_5 \in E(G)$ (otherwise the set $\{u'_1,u'_5\} \in IS(G_{v,u_3})$ dominates $\{u_1,u_5\}$, a contradiction). Hence $u_5 \in ext(vu_1u_2u_3u_4)$, we may assume that the vertices $u_1,\ldots,u_5$ are located clockwise around $v$. The set $(U_1 \cup U_5) \setminus U_3$ is a clique, thus $|(U_1 \cup U_5) \setminus U_3| \leq 3$ and $|(U_1 \cup U_5) \cap U_3| \geq 3$. Hence there exist vertices $x_1,x_2 \in U_i \cap U_3$ for some $i \in \{1,5\}$ (say, $i = 1$). Assume that $x_2 \in ext(vu_1x_1u_3)$. Consider a vertex $u'_2 \in U_2$, then the set $\{u'_2,x_2\} \in IS(G_{v,u_5})$ dominates $\{u_1,u_2,u_3\}$, a contradiction by Corollary~\ref{cor1}(a).

\textbf{Case 3.2.2.} $U_1 \subseteq U_3$ or $U_5 \subseteq U_3$ (assume that $U_1 \subseteq U_3$). There exist vertices $x_1,x_2,x_3 \in U_1 \cap U_3$ such that $\{x_2\} = int(u_1x_1u_3x_3) \cap U_1$. Then it is easy to check, using Corollary~\ref{cor1}(a), that $U_2 \subseteq U_5$. Hence $u_4 \in int(vu_2u_3)$ and $u_5 \in int(vu_2u_3u_4)$. Consider a vertex $u'_2 \in U_2$, note that $u'_2u_4 \notin E(G)$. By Lemma~\ref{l15}, there exists a vertex $z \in int(u_1x_1u_3x_2) \cup int(u_1x_2u_3x_3)$ not adjacent to $N_1[v] \cup \{x_1,x_2,x_3,u'_2\}$. Consider a family $\mathcal{C} \in MCP(G)$ and its clique $W \ni v$. Clearly, $\{u_3,x_1,x_2,x_3\}$ is not a clique, then there exists $j \in [1;3]$ such that the set $\{u_3,x_j\}$ is $\theta$-independent in $\mathcal{C}$. Hence one of the $\theta$-independent sets $\{u_3,x_j\}$, $\{u_1,u_4\}$, $\{u_2,u_5\}$, $\{u_3,u'_2\}$ dominates $W$, a contradiction.
\end{proof}

\begin{lemma}\label{l40}
If $N_1(v) \cong P_4 \cup K_1$, then $G$ is not critical.
\end{lemma}

\begin{proof}

Assume that $E(N_1(v)) = \{u_1u_2,u_2u_3,u_3u_4\}$ and $u_1 \in ext(vu_2u_3)$. Consider three cases:

\textbf{Case 1.} There exists a vertex $w \in N_2(v)$ with at least four neighbors in $N_1(v)$. Corollary~\ref{cor1}(b) implies that $u_1w,u_4w,u_5w \in E(G)$. Then $u_4 \in ext(vu_1u_2u_3)$ and $u_5,w \in ext(vu_1u_2u_3u_4)$. We may assume that $u_3w \in E(G)$ and the vertices $u_1,\ldots,u_5$ are in clockwise order around $v$. By Lemma~\ref{l15}, there exists a vertex $x \in int(vu_3wu_5)$ not adjacent to $N_1[v] \cup U_2 \cup \{w\}$. We use the argument from Case~1.1 of the previous lemma to obtain a contradiction.

\textbf{Case 2.} There exists a vertex $w \in N_2(v)$ such that $wu_i, wu_{i+2} \in E(G)$ for some $i \in \{1,2\}$. Assume by symmetry that $i = 2$.

\textbf{Case 2.1.} $wu_3 \in E(G)$. The set $\{u_1,w\} \in IS(G)$ dominates $\{v,u_2,u_3,u_4\}$, thus $\gamma(G_{v,u_3}) < \gamma(G_{u_3})$ by Lemmas~\ref{l13} and~\ref{l14}. Lemma~\ref{l10}(b) and Corollary~\ref{cor2} imply that there exists a vertex $x \in U_1 \cap U_5$ such that $wx \notin E(G)$. Since $\gamma(G_{v,u_3}) < \gamma(G_{u_3})$, we have $xu_3 \in E(G)$, then it is easy to see that $u_4,w \in int(vu_2u_3)$ and $x,u_5 \in ext(vu_1u_2u_3)$.  Therefore, $w$ is not adjacent to $U_1 \cup U_5$ and $\gamma(G_w) < \theta(G_w)$ by Lemma~\ref{l10}(d), a contradiction.

\textbf{Case 2.2.} $wu_3 \notin E(G)$. 

\textbf{Case 2.2.1.} $w$ does not dominate $U_3$. Lemma~\ref{l10}(b,e) implies that $wu_1 \in E(G)$ and $wu_5 \notin E(G)$. Moreover, there exists a vertex $x \in U_3 \cap U_5$ (Corollary~\ref{cor1}(b) implies that $xu_4 \notin E(G)$. Clearly, $u_5 \in int(vu_2u_3u_4)$. Two configurations are possible:

\textbf{Case 2.2.1.1.} $u_5 \in int(vu_2u_3)$. Then $U_4 \subseteq U_1$ by Corollary~\ref{cor1}(a) and there exists vertices $y_1,y_2 \in U_1 \cap U_4$ such that $ext(vu_1y_2u_4) \cap U_4 = \{y_1\}$. Then $ext(vu_1y_1u_4) = \emptyset$ and by lemma~\ref{l15} there exists a vertex $z \in int(u_1y_1u_4y_2)$ not adjacent to $N_1[v] \cup \{x,y_1\}$. Consider a family $\mathcal{C} \in MCP(G_z)$ and its clique $W \ni v$. It is easy to check that one of the independent sets $\{u_1,u_4\}$, $\{u_1,x\}$, $\{u_3,y_1\}$, $\{u_2,y_1\}$ dominates $W$; a contradiction.

\textbf{Case 2.2.1.2.} $u_5 \in int(vu_3u_4)$. Corollary~\ref{cor1} implies that $U_4 \cap int(u_2u_3u_4w) = \emptyset$. Moreover, for every vertex $u'_4 \in U_4$ either $u'_4 \in U_1$ or $u'_4x \in E(G)$. 

If $U_4 \subseteq N_1[x]$, then $deg(u_4) = 5$ by Corollary~\ref{cor2}. Moreover, $G[N_1(u_4)] \cong P_4 \cup K_1$, hence there exists a vertex $u'_4 \in U_3 \cap U_4$. Since $|U_5| \geq 4$, there exists a vertex $u'_5 \in U_5$ such that the set $\{u'_4,u'_5\}$ is $\theta$-independent in $G$, a contradiction by Corollary~\ref{cor1}(a). 

Therefore, there exists a vertex $y \in U_1 \cap U_4$. Clearly, $y \in ext(vu_1wu_4)$. By lemma~\ref{l17}, $int(vu_1u_2) = int(vu_2u_3) = \emptyset$. Hence by lemma~\ref{l15} there exists a vertex $z \in int(u_1u_2w) \cup int(u_2u_3u_4w)$ not adjacent to $N_1[v] \cup \{x,y\}$. Consider a family $\mathcal{C} \in MCP(G_z)$ and its clique $W \ni v$. It is easy to check that one of the independent sets $\{u_1,u_4\}$, $\{u_1,x\}$, $\{u_3,y\}$, $\{u_2,y\}$ dominates $W$; a contradiction.

\textbf{Case 2.2.2.} $w$ dominates $U_3$. Corollary~\ref{cor2} imply that $|U_3| = 2$. Let $U_3 = \{u'_3,u''_3\}$, assume by symmetry that $u'_3 \in int(u_2u_3u''_3w)$. Lemmas~\ref{l30}--\ref{l39} imply that $N_1[u_3] \cong P_4 \cup K_1$, thus $u'_3u''_3 \notin E(G)$ and either $u_2u'_3 \in E(G)$ or $u_4u''_3 \in E(G)$. 

\textbf{Case 2.2.2.1.} $u_2u'_3 \in E(G)$. Corollary~\ref{cor2} implies that the vertices $u_1,u_4,u_5$ are not pendant in $G_{u'_3}$, hence $\gamma(G_{u'_3}) < \theta(G_{u'_3})$ by Lemma~\ref{l10}(e), a contradiction. 

\textbf{Case 2.2.2.2.} $u_4u''_3 \in E(G)$. If $int(u_3u'_3wu''_3) \neq \emptyset$, then $N_4[v] \neq \emptyset$ by Lemma~\ref{l15}, a contradiction by Lemma~\ref{l29}. Thus $int(u_2u_3u'_3w) \neq \emptyset$, and by Lemma~\ref{l17} there exists a vertex $u'_2 \in int(u_2u_3u'_3w) \cap U_2$. If $U_1 \cap U_5 = \emptyset$, then $\gamma(G_{u'_2,u''_3}) < \theta(G_{u'_2,u''_3})$ by Lemma~\ref{l10}(b), a contradiction. Otherwise there exists a vertex $x \in U_1 \cap U_5$. Clearly, $u_5 \in ext(vu_1u_2u_3u_4)$. If there exists a vertex $y \in U_5 \setminus U_1$, then $\{u''_3,y\} \in IS(G_{v,u_1})$ dominates $\{u_3,u_4,u_5\}$, a contradiction. Hence $U_5 \subseteq U_1$ and $\gamma(G_{u'_2,u''_3}) < \theta(G_{u'_2,u''_3})$ by Lemma~\ref{l10}(d), a contradiction.

\textbf{Case 3.} For both $i \in \{1,2\}$ no vertex from $N_2(v)$ dominates $\{u_i,u_{i+2}\}$. By Lemma~\ref{l14}, it suffices to consider the following subcases:

\textbf{Case 3.1.} There exists $i \in \{1;2\}$ such that $$\gamma(G \setminus \{v,u_i,u_{i+1},u_{i+2}\}) = \theta(G \setminus \{v,u_i,u_{i+1},u_{i+2}\}) = \gamma(G) - 1.$$ 

Assume by symmetry that $i = 2$.

\textbf{Case 3.1.1.} There exists a vertex $x \in U_2 \cap U_3$. Case~2 implies that $xu_1,xu_4 \notin E(G)$. Then Lemma~\ref{l10}(e) and Corollary~\ref{cor2} imply that $\gamma(G_x) < \theta(G_x)$, a contradiction.

\textbf{Case 3.1.2.} $U_2 \cap U_3 = \emptyset$. Moreover, $|U_1 \cap U_4| \geq 2$. Then there exist vertices $w_1,w_2 \in U_1 \cap U_4$ such that $U_1 \cap ext(vu_1w_1u_4) = \{w_2\}$. 

\textbf{Case 3.1.2.1.} $u_5 \in int(vu_1u_2u_3u_4)$. Assume by symmetry that $u_5 \notin int(vu_1u_2)$. Select a vertex $u'_5 \in U_5$. 

If there exists a vertex $z \in ext(vu_1w_1u_4)$ not adjacent to $N_1[v] \cup U_5 \cup \{w_1\}$, consider a family $\mathcal{C} \in MCP(G_z)$ and its clique $W \ni v$. It is easy to check that one of the sets $\{u_3,w_1\}, \{u_2,w_1\}, \{u_1,u_4\}, \{u_1,u'_5\} \in IS(G_z)$ dominates $W$, a contradiction by Lemma~\ref{l9}. 

Suppose that all vertices from $ ext(vu_1w_1u_4)$ are adjacent to $N_1[v] \cup U_5 \cup \{w_1\}$. By Lemma~\ref{l17}, $ext(vu_1w_2u_4) = \emptyset$, then $int(u_1w_1u_4w_2) \neq \emptyset$. By Lemma~\ref{l15} there exists a vertex $u'_4 \in U_4 \cap int(u_1w_1u_4w_2)$. Consider a vertex $u'_3 \in U_3$, remind that $u'_3u_1 \notin E(G)$. The set $\{u_1,u'_3,u'_4\} \in IS(G)$ dominates $\{v,u_2,u_3,u_4\}$, hence $\gamma(G \setminus \{v,u_2,u_3,u_4\}) \geq \gamma(G)$, a contradiction.

\textbf{Case 3.1.2.2.} $u_5 \in ext(vu_1u_2u_3u_4)$. Since $U_2 \cap U_3 = \emptyset$, by Lemma~\ref{l1} there exist nonadjacent vertices $u'_2 \in U_2 \setminus U_5$ and $u'_3 \in U_3 \setminus U_5$. If $w_2u_5 \notin E(G)$, the set $\{u'_2,u'_3,w_2\} \in IS(G_{v,u_5})$ dominates $\{u_1,u_2,u_3,u_4\}$, a contradiction. Otherwise Lemma~\ref{l10}(e) implies that $\gamma(G_{u'_2,u'_3}) < \theta(G_{u'_2,u'_3})$, a contradiction.

\textbf{Case 3.1.3.} $U_2 \cap U_3 = \emptyset$. Moreover, $|U_1 \cap U_4| \leq 1$. Then $|U_i \setminus U_1| \geq 2$ for both $i \in \{3,4\}$. If there exists a vertex $x \in U_3 \cap U_4$, then $\{x,u_1\} \in IS(G)$ dominates $\{v,u_2,u_3,u_4\}$ and $\gamma(G) \leq \gamma(G \setminus \{v,u_2,u_3,u_4\})$ by Lemma~\ref{l8}, a contradiction. Otherwise by Lemma~\ref{l1} there exist nonadjacent vertices $u'_3 \in U_3 \setminus U_1$ and $u'_4 \in U_4 \setminus U_1$. The set $\{u_1,u'_3,u'_4\} \in IS(G)$ dominates $\{v,u_2,u_3,u_4\}$, again a contradiction by Lemma~\ref{l8}.

\textbf{Case 3.2.} $\gamma(G_{u_2,v}) = \gamma(G_{u_3,v}) = \gamma(G) - 2$. Since $G$ is $K_{3,3}$-free and $|U_5| \geq 4$, there exists a vertex $u'_5 \in U_5$ such that $u'_5u_i \notin E(G)$ for some $i \in \{2,3\}$. Assume by symmetry that $i = 3$. Remind that $U_1 \cap U_3 = \emptyset$. By Corollary~\ref{cor2} there exists a vertex $u'_1 \in U_1$ nonadjacent to $u'_5$, then the set $\{u'_1,u'_5\} \in IS(G_{u_3})$ dominates $\{u_1,u_5\}$. Therefore, $\gamma(G_{v,u_3}) \geq \gamma(G_{u_3})$ by Lemma~\ref{l8}, a contradiction.

\end{proof}

\section*{Acknowledgments}

The article was prepared within the framework of the Basic Research Program at the National Research University Higher School of Economics (HSE).

\end{document}